\documentclass{gtart}

\def\ifplaintex{\expandafter\ifx\csname documentclass\endcsname\relax}

\ifplaintex 
\else
\fi

\expandafter\ifx\csname beginpicture\endcsname\relax
\expandafter\ifx\csname documentclass\endcsname\relax
\input pictex \else
\input prepictex \input pictex \input postpictex \fi\fi

\def\gt{{\mathsurround=0pt\it $\cal G\mskip-2mu$eometry \&\ 
$\cal T\!\!$opology}}        %  journal title in recommended style

\def\gtp{{\mathsurround=0pt\it $\cal G\mskip-2mu$eometry \&\ 
$\cal T\!\!$opology $\cal P\!$ublications}}  % GT publications

\def\lognumber#1{\def\thelognumber{#1}}
\def\volumenumber#1{\def\thevolumenumber{#1}}
\def\papernumber#1{\def\thepapernumber{#1}}
\def\volumeyear#1{\def\thevolumeyear{#1}}

\def\pagenumbers#1#2{\def\startpage{#1}\def\finishpage{#2}}
\def\published#1{\def\publishdate{#1}}
\def\proposed#1{\def\theproposer{#1}}
\def\seconded#1{\def\theseconders{#1}}
\def\received#1{\def\receiveddate{#1}}
\def\revised#1{\def\reviseddate{#1}}
\def\accepted#1{\def\accepteddate{#1}}

\def\coverauthors#1{\def\thecoverauthors{#1}}
\def\asciiauthors#1{\def\theasciiauthors{#1}}
\def\asciiaddress#1{\def\theasciiaddress{#1}}

\long\def\asciiabstract#1{\long\def\theasciiabstract{#1}}
\def\asciikeywords#1{\def\theasciikeywords{#1}}

\let\\\par\let\thelognumber\relax
\let\thevolumenumber\relax\let\thepapernumber\relax
\let\thevolumeyear\relax\let\thesamplenumber\relax\let\startpage\relax
\let\finishpage\relax\let\publishdate\relax\let\receiveddate\relax
\let\reviseddate\relax\let\accepteddate\relax\let\theasciititle\relax
\let\theasciiauthors\relax\let\theasciiaddress\relax
\let\theasciiabstract\relax\let\theasciikeywords\relax
\let\theasciiemail\relax\let\theshortauthors\relax\let\theshorttitle\relax
\let\thecoverauthors\relax

\long\def\maketitlep{   % start of definition of \maketitlep

\count0=\startpage

\gt\hfill      %   Journal title (top left) 
\beginpicture
\setcoordinatesystem units <0.33truein, 0.33truein> point at 2.2 0.9
\setplotsymbol ({$\cal G$})
\plotsymbolspacing=9truept
\circulararc 315 degrees from 0 1 center at 0 0
\setplotsymbol ({$\cal T$})
\circulararc 315 degrees from 1 -1 center at 1 0
\endpicture
\break
{\small\ifx\thesamplenumber\relax % sample?  
Volume \else Sample
\fi\thevolumenumber\ (\thevolumeyear)
\startpage--\finishpage\nl
Published: \publishdate}
\vglue 0.5truein plus 0.4fil minus 0.1truein

{\parskip=0pt\leftskip 0pt plus 1fil\def\\{\par\smallskip}{\ifplaintex\large
\else\Large\fi\bf\thetitle}\par\medskip}   

\vglue 0pt plus 0.1fil 

{\parskip=0pt\leftskip 0pt plus 1fil\def\\{\par}{\sc\theauthors}
\par\medskip}

\vglue 0pt plus 0.1fil 

{\small\parskip=0pt\let\newline\\
{\leftskip 0pt plus 1fil\def\\{\par}{\sl\theaddress}\par}
\expandafter\ifx\theemail\relax    % email address?
\relax\else\vglue 5pt plus 0.02fil minus 2pt\def\\{\stdspace{\rm 
and}\stdspace} 
\cl{Email:\stdspace\tt\theemail}\fi
\ifx\theurl\relax                  % URL given?
\relax\else\vglue 5pt plus 0.02fil minus 2pt\def\\{\stdspace{\rm 
and}\stdspace}
\cl{URL:\stdspace\tt\theurl}\fi\par}

\vglue 7pt plus 0.3fil minus 3pt

{\bf Abstract}
\vglue 5pt plus 0.1fil minus 2pt

\theabstract

\vglue 7pt plus 0.3fil minus 3pt

{\bf AMS Classification numbers}\quad Primary:\quad \theprimaryclass

Secondary:\quad \thesecondaryclass

\vglue 5pt plus 0.3fil minus 2pt

{\bf Keywords}\quad \thekeywords

\vglue 10pt plus 0.5fil minus 5pt

{\small  Proposed: \theproposer\hfill Received: \receiveddate\nl
Seconded: \theseconders\hfill 
\ifx\reviseddate\relax                         % paper revised?
Accepted: \accepteddate                        % no
\else
Revised: \reviseddate                          % yes
\fi}
\eject
}       %  end of definition of \maketitlep

\let\maketitlepage\maketitlep
\let\maketitle\maketitlepage

\font\phead=cmsl9 scaled 950
\font=cmsl9 scaled 1050
\font\pnum=cmbx10 scaled 913
\font\lnum=cmbx10 
\font\pfoot=cmsl9 scaled 950
\font=cmsl9 scaled 1050
\ifplaintex
\headline{\vbox to 0pt{\vskip -4.5mm\line{\small\phead\ifnum
\count0=\startpage ISSN 1364-0380 (on line)
1465-3060 (printed) \hfill {\pnum\folio}\else\ifodd\count0\def\\{ }% 
\ifx\theshorttitle\relax\thetitle\else\theshorttitle\fi\hfill{\pnum\folio}
\else\def\\{ and }{\pnum\folio}\hfill\ifx\theshortauthors\relax\theauthors
\else\theshortauthors\fi\fi\fi}\vss}}
\footline{\vbox to 0pt{\vglue 0mm\line{\small\pfoot\ifnum\count0=\startpage
\copyright\ \gtp\hfill\else
\gt, Volume \thevolumenumber\ (\thevolumeyear)\hfill\fi}\vss
}}
\else
\makeatletter
\def\@oddhead{{\small\lhead\ifnum\count0=\startpage ISSN 1364-0380 (on line)
1465-3060 (printed) \hfill {\lnum\number\count0}\else\ifodd\count0
\def\\{ }\ifx\theshorttitle\relax \thetitle \else\theshorttitle\fi\hfill
{\lnum\number\count0}\else\def\\{ and }{\lnum\number\count0}
\hfill\ifx\theshortauthors\relax 
\theauthors\else\theshortauthors\fi\fi\fi}}\def\@evenhead{\@oddhead}
\def\@oddfoot{\small\lfoot\ifnum\count0=\startpage\copyright\ \gtp\hfill\else
\gt, Volume \thevolumenumber\ (\thevolumeyear)\hfill\fi}
\def\@evenfoot{\@oddfoot}
\makeatother
\fi

\newwrite\gtoutfile
\long\gdef\makeheadfile{  %%% start of definition of \makeheadfile
{\def\\{, }\def\s{ }
\immediate\openout\gtoutfile head.xxx
\immediate\write\gtoutfile{To: math@arxiv.org}
\immediate\write\gtoutfile{Subject: put or rep NNNNN:pppp}
\immediate\write\gtoutfile{--text follows this line--}
\immediate\write\gtoutfile{Proxy-for: \ifx\theasciiauthors\relax
\theauthors\else\theasciiauthors\fi\s<\ifx\theasciiemail\relax\theemail\else\theasciiemail\fi>}
\immediate\write\gtoutfile{\noexpand\\}
\immediate\write\gtoutfile{Authors: \ifx\theasciiauthors\relax
\theauthors\else\theasciiauthors\fi}
{\def\\{ }\immediate\write\gtoutfile{Title: \ifx\theasciititle\relax
\thetitle\else\theasciititle\fi}}
\immediate\write\gtoutfile{Subj-class: GT or SG or MG etc}
\immediate\write\gtoutfile{MSC-class: \theprimaryclass\ifx\thesecondaryclass\relax\else, \thesecondaryclass\fi}
\immediate\write\gtoutfile{Journal-ref: Geom. Topol. \thevolumenumber
(\thevolumeyear) \startpage-\finishpage}
\immediate\write\gtoutfile{Comments: Published by Geometry and Topology at}
\immediate\write\gtoutfile{\s\s http://www.maths.warwick.ac.uk/gt/GTVol\thevolumenumber/paper\thepapernumber.abs.html}
\immediate\write\gtoutfile{\noexpand\\}
\immediate\write\gtoutfile{}
\ifx\theasciiabstract\relax
\immediate\write\gtoutfile{\theabstract}\else
\immediate\write\gtoutfile{\theasciiabstract}\fi
\immediate\write\gtoutfile{}
\immediate\write\gtoutfile{\noexpand\\}
\immediate\write\gtoutfile{}
\immediate\closeout\gtoutfile}}  %%% end of definition of \makeheadfile

\def\maketitlepage{\maketitlep\makeheadfile}
\let\maketitle\maketitlepage

\def\ifplaintex{\expandafter\ifx\csname documentclass\endcsname\relax}

\ifplaintex 
\else
\fi

\expandafter\ifx\csname beginpicture\endcsname\relax
\expandafter\ifx\csname documentclass\endcsname\relax
\input pictex \else
\input prepictex \input pictex \input postpictex \fi\fi

\def\gt{{\mathsurround=0pt\it $\cal G\mskip-2mu$eometry \&\ 
$\cal T\!\!$opology}}        %  journal title in recommended style

\def\gtp{{\mathsurround=0pt\it $\cal G\mskip-2mu$eometry \&\ 
$\cal T\!\!$opology $\cal P\!$ublications}}  % GT publications

\def\lognumber#1{\def\thelognumber{#1}}
\def\volumenumber#1{\def\thevolumenumber{#1}}
\def\papernumber#1{\def\thepapernumber{#1}}
\def\volumeyear#1{\def\thevolumeyear{#1}}

\def\pagenumbers#1#2{\def\startpage{#1}\def\finishpage{#2}}
\def\published#1{\def\publishdate{#1}}
\def\proposed#1{\def\theproposer{#1}}
\def\seconded#1{\def\theseconders{#1}}
\def\received#1{\def\receiveddate{#1}}
\def\revised#1{\def\reviseddate{#1}}
\def\accepted#1{\def\accepteddate{#1}}

\def\coverauthors#1{\def\thecoverauthors{#1}}
\def\asciiauthors#1{\def\theasciiauthors{#1}}
\def\asciiaddress#1{\def\theasciiaddress{#1}}

\long\def\asciiabstract#1{\long\def\theasciiabstract{#1}}
\def\asciikeywords#1{\def\theasciikeywords{#1}}

\let\\\par\let\thelognumber\relax
\let\thevolumenumber\relax\let\thepapernumber\relax
\let\thevolumeyear\relax\let\thesamplenumber\relax\let\startpage\relax
\let\finishpage\relax\let\publishdate\relax\let\receiveddate\relax
\let\reviseddate\relax\let\accepteddate\relax\let\theasciititle\relax
\let\theasciiauthors\relax\let\theasciiaddress\relax
\let\theasciiabstract\relax\let\theasciikeywords\relax
\let\theasciiemail\relax\let\theshortauthors\relax\let\theshorttitle\relax
\let\thecoverauthors\relax

\long\def\maketitlep{   % start of definition of \maketitlep

\count0=\startpage

\gt\hfill      %   Journal title (top left) 
\beginpicture
\setcoordinatesystem units <0.33truein, 0.33truein> point at 2.2 0.9
\setplotsymbol ({$\cal G$})
\plotsymbolspacing=9truept
\circulararc 315 degrees from 0 1 center at 0 0
\setplotsymbol ({$\cal T$})
\circulararc 315 degrees from 1 -1 center at 1 0
\endpicture
\break
{\small\ifx\thesamplenumber\relax % sample?  
Volume \else Sample
\fi\thevolumenumber\ (\thevolumeyear)
\startpage--\finishpage\nl
Published: \publishdate}
\vglue 0.5truein plus 0.4fil minus 0.1truein

{\parskip=0pt\leftskip 0pt plus 1fil\def\\{\par\smallskip}{\ifplaintex\large
\else\Large\fi\bf\thetitle}\par\medskip}   

\vglue 0pt plus 0.1fil 

{\parskip=0pt\leftskip 0pt plus 1fil\def\\{\par}{\sc\theauthors}
\par\medskip}

\vglue 0pt plus 0.1fil 

{\small\parskip=0pt\let\newline\\
{\leftskip 0pt plus 1fil\def\\{\par}{\sl\theaddress}\par}
\expandafter\ifx\theemail\relax    % email address?
\relax\else\vglue 5pt plus 0.02fil minus 2pt\def\\{\stdspace{\rm 
and}\stdspace} 
\cl{Email:\stdspace\tt\theemail}\fi
\ifx\theurl\relax                  % URL given?
\relax\else\vglue 5pt plus 0.02fil minus 2pt\def\\{\stdspace{\rm 
and}\stdspace}
\cl{URL:\stdspace\tt\theurl}\fi\par}

\vglue 7pt plus 0.3fil minus 3pt

{\bf Abstract}
\vglue 5pt plus 0.1fil minus 2pt

\theabstract

\vglue 7pt plus 0.3fil minus 3pt

{\bf AMS Classification numbers}\quad Primary:\quad \theprimaryclass

Secondary:\quad \thesecondaryclass

\vglue 5pt plus 0.3fil minus 2pt

{\bf Keywords}\quad \thekeywords

\vglue 10pt plus 0.5fil minus 5pt

{\small  Proposed: \theproposer\hfill Received: \receiveddate\nl
Seconded: \theseconders\hfill 
\ifx\reviseddate\relax                         % paper revised?
Accepted: \accepteddate                        % no
\else
Revised: \reviseddate                          % yes
\fi}
\eject
}       %  end of definition of \maketitlep

\let\maketitlepage\maketitlep
\let\maketitle\maketitlepage

\font\phead=cmsl9 scaled 950
\font=cmsl9 scaled 1050
\font\pnum=cmbx10 scaled 913
\font\lnum=cmbx10 
\font\pfoot=cmsl9 scaled 950
\font=cmsl9 scaled 1050
\ifplaintex
\headline{\vbox to 0pt{\vskip -4.5mm\line{\small\phead\ifnum
\count0=\startpage ISSN 1364-0380 (on line)
1465-3060 (printed) \hfill {\pnum\folio}\else\ifodd\count0\def\\{ }% 
\ifx\theshorttitle\relax\thetitle\else\theshorttitle\fi\hfill{\pnum\folio}
\else\def\\{ and }{\pnum\folio}\hfill\ifx\theshortauthors\relax\theauthors
\else\theshortauthors\fi\fi\fi}\vss}}
\footline{\vbox to 0pt{\vglue 0mm\line{\small\pfoot\ifnum\count0=\startpage
\copyright\ \gtp\hfill\else
\gt, Volume \thevolumenumber\ (\thevolumeyear)\hfill\fi}\vss
}}
\else
\makeatletter
\def\@oddhead{{\small\lhead\ifnum\count0=\startpage ISSN 1364-0380 (on line)
1465-3060 (printed) \hfill {\lnum\number\count0}\else\ifodd\count0
\def\\{ }\ifx\theshorttitle\relax \thetitle \else\theshorttitle\fi\hfill
{\lnum\number\count0}\else\def\\{ and }{\lnum\number\count0}
\hfill\ifx\theshortauthors\relax 
\theauthors\else\theshortauthors\fi\fi\fi}}\def\@evenhead{\@oddhead}
\def\@oddfoot{\small\lfoot\ifnum\count0=\startpage\copyright\ \gtp\hfill\else
\gt, Volume \thevolumenumber\ (\thevolumeyear)\hfill\fi}
\def\@evenfoot{\@oddfoot}
\makeatother
\fi

\newwrite\gtoutfile
\long\gdef\makeheadfile{  %%% start of definition of \makeheadfile
{\def\\{, }\def\s{ }
\immediate\openout\gtoutfile head.xxx
\immediate\write\gtoutfile{To: math@arxiv.org}
\immediate\write\gtoutfile{Subject: put or rep NNNNN:pppp}
\immediate\write\gtoutfile{--text follows this line--}
\immediate\write\gtoutfile{Proxy-for: \ifx\theasciiauthors\relax
\theauthors\else\theasciiauthors\fi\s<\ifx\theasciiemail\relax\theemail\else\theasciiemail\fi>}
\immediate\write\gtoutfile{\noexpand\\}
\immediate\write\gtoutfile{Authors: \ifx\theasciiauthors\relax
\theauthors\else\theasciiauthors\fi}
{\def\\{ }\immediate\write\gtoutfile{Title: \ifx\theasciititle\relax
\thetitle\else\theasciititle\fi}}
\immediate\write\gtoutfile{Subj-class: GT or SG or MG etc}
\immediate\write\gtoutfile{MSC-class: \theprimaryclass\ifx\thesecondaryclass\relax\else, \thesecondaryclass\fi}
\immediate\write\gtoutfile{Journal-ref: Geom. Topol. \thevolumenumber
(\thevolumeyear) \startpage-\finishpage}
\immediate\write\gtoutfile{Comments: Published by Geometry and Topology at}
\immediate\write\gtoutfile{\s\s http://www.maths.warwick.ac.uk/gt/GTVol\thevolumenumber/paper\thepapernumber.abs.html}
\immediate\write\gtoutfile{\noexpand\\}
\immediate\write\gtoutfile{}
\ifx\theasciiabstract\relax
\immediate\write\gtoutfile{\theabstract}\else
\immediate\write\gtoutfile{\theasciiabstract}\fi
\immediate\write\gtoutfile{}
\immediate\write\gtoutfile{\noexpand\\}
\immediate\write\gtoutfile{}
\immediate\closeout\gtoutfile}}  %%% end of definition of \makeheadfile

\def\maketitlepage{\maketitlep\makeheadfile}
\let\maketitle\maketitlepage

\lognumber{224}

\volumenumber{6}\papernumber{9}\volumeyear{2002}
\pagenumbers{269}{328}
\received{11 January 2002}
\revised{25 April 2002}
\published{20 May 2002}
\accepted{17 May 2002}

\proposed{Walter Neumann}
\seconded{Robion Kirby, Yasha Eliashberg}

\usepackage{amsmath,amssymb}
\usepackage{euscript,graphicx}
\input labelfig

\makeatletter
\renewcommand{\subsection}{\@startsection{subsection}{2}{\z@}{-3.5ex
plus-1ex minus-.2ex}{-\baselineskip}{\raggedright\reset@font\normalsize\bf}}

\makeatother
\def\gtsubsection#1{\subsection{}{\bf\hspace{-1.25em}#1}\qua}

\DeclareMathOperator{\Hom}{{\rm Hom}}
\DeclareMathOperator{\rank}{{\rm rank}}

\newcommand{\bP}{{{\mathbb P}}}
\newcommand{\bL}{{{\mathbb L}}}
\newcommand{\bC}{{\mathbb C}}
\newcommand{\bQ}{{{\mathbb Q}}}
\newcommand{\bR}{{\mathbb R}}
\newcommand{\bS}{{\mathbb S}}
\newcommand{\bT}{{\mathbb T}}

\newcommand{\bZ}{{\mathbb Z}}

\newcommand{\gq}{{\mathfrak q}}
\newcommand\aug{\mathfrak{aug}}
\newcommand{\dir}{{\mathfrak D}}

\newcommand{\LL}{\left(\left( }
\newcommand{\RR}{\right)\right)}
\newcommand{\bms}{\mbox{\boldmath$s$}}

\newcommand{\bSW}{\mbox{\boldmath$SW$}}
\newcommand{\cC}{\mathcal C}
\newcommand{\p}{\EuScript{P}}
\newcommand{\cv}{\EuScript{V}}
\newcommand{\gauge}{\EuScript{G}}
\newcommand{\f}{\EuScript{F}}
\newcommand{\modu}{{\mathfrak M}}
\newcommand{\R}{\mathbb R}
\newcommand{\coker}{\mbox{coker}}
\newcommand{\ra}{\rightarrow}
\newcommand{\im}{{\rm Im}}
\newcommand{\bif}{(\, , \,)}
\newcommand{\mmod}{\mbox{mod}}
\newcommand{\Z}{\mathbb Z}
\newcommand{\Q}{\mathbb Q}
\newcommand{\si}{{\sigma}}
\newcommand{\ssw}{{\bf sw}}
\newcommand{\C}{\mathbb C}
\newcommand{\eps}{\epsilon}
\newcommand{\B}{\EuScript{B}}
\newcommand{\ii}{{\bf i}}
\newcommand{\bmm}{\mbox{\boldmath$m$}}
\newcommand{\lan}{\langle}
\newcommand{\ran}{\rangle}
\newcommand{\pic}{{\rm Pic}^V_{top}(\Sigma)}
\newcommand{\picm}{{\rm Pic}_{top}(M)}

\newcommand{\tT}{\EuScript{T}}

\begin{document}

\title{Seiberg--Witten invariants and surface singularities}

\authors{Andr\'as N\'emethi\\Liviu I Nicolaescu}
\coverauthors{Andr\noexpand\'as N\noexpand\'emethi\\Liviu I Nicolaescu}
\asciiauthors{Andras Nemethi\\Liviu I Nicolaescu}

\address{Department of Mathematics, Ohio State University\\Columbus, OH 43210, USA\\\smallskip\\Department of Mathematics, University of Notre Dame\\Notre Dame, IN 46556, USA}
\asciiaddress{Department of Mathematics, Ohio State University\\Columbus, OH 43210, USA\\Department of Mathematics, University of Notre Dame\\Notre Dame, IN 46556, USA}
\emails{nemethi@math.ohio-state.edu\\nicolaescu.1@nd.edu}
\url{www.math.ohio-state.edu/\char'176nemethi/\\www.nd.edu/\char'176lnicolae/}



\asciikeywords{(Links of) surface singularities, (Q)-Gorenstein
singularities, rational singularities, Brieskorn-Hamm complete
intersections, geometric genus, Seiberg-Witten invariants of
Q-homology spheres, Reidemeister-Turaev torsion, Casson-Walker
invariant}

\begin{abstract}
{\small
We formulate a very general conjecture relating the analytical
invariants of a normal surface singularity to the Seiberg--Witten
invariants of its link provided that the link is a rational homology
sphere.  As supporting evidence, we establish its validity for a large
class of singularities: some rational and minimally elliptic
(including the cyclic quotient and ``polygonal'') singularities, and
Brieskorn--Hamm complete intersections.  Some of the verifications are
based on a result which describes (in terms of the plumbing graph) the
Reidemeister--Turaev sign refined torsion (or, equivalently, the
Seiberg--Witten invariant) of a rational homology 3--manifold $M$,
provided that $M$ is given by a negative definite plumbing.

These results extend previous work of Artin, Laufer and
S\thinspace S-T Yau, respectively of  Fintushel--Stern and Neumann--Wahl.}
\end{abstract}

\asciiabstract{
We formulate a very general conjecture relating the analytical
invariants of a normal surface singularity to the Seiberg-Witten
invariants of its link provided that the link is a rational homology
sphere.  As supporting evidence, we establish its validity for a large
class of singularities: some rational and minimally elliptic
(including the cyclic quotient and `polygonal') singularities, and
Brieskorn-Hamm complete intersections.  Some of the verifications are
based on a result which describes (in terms of the plumbing graph) the
Reidemeister-Turaev sign refined torsion (or, equivalently, the
Seiberg-Witten invariant) of a rational homology 3-manifold M,
provided that M is given by a negative definite plumbing.
These results extend previous work of Artin, Laufer and
S S-T Yau, respectively of Fintushel-Stern and Neumann-Wahl.}

\maketitlepage

\section{Introduction}

The main goal of the present paper is to formulate  a very
general conjecture  which relates the topological and the
analytical invariants  of a complex normal surface singularity
whose link is  a rational homology sphere.

The motivation for such a result comes from several directions.
Before we present  some of them, we fix some notations.

Let $(X,0)$ be a normal two-dimensional analytic  singularity. It
is well-known that from a topological point of view, it   is
completely characterized by its link $M$, which is an oriented
3--manifold. Moreover, by a result of Neumann \cite{NP}, any
decorated resolution graph of $(X,0)$ carries the same information
as $M$. A property of $(X,0)$ will be called \emph{topological} if
it can be determined from  $M$, or equivalently, from any
resolution graph of $(X,0)$.

It is interesting  to investigate, in which cases (some
of) the  analytical invariants (determined, say, from the local
algebra of $(X,0)$) are topological. In this
article we are mainly interested in the geometric genus $p_g$ of
$(X,0)$ (for details, see section 4).

Moreover, if $(X,0)$ has a smoothing with Milnor fiber $F$, then
one can ask the same question about the signature $\si(F)$ and the
topological Euler characteristic $\chi_{top}(F)$ of $F$ as well.
It is known (via some results of Laufer, Durfee, Wahl  and
Steenbrink) that for Gorenstein singularities, any of $p_g$,
$\si(F)$ and $\chi_{top}(F)$ determines  the remaining two  modulo
a  certain invariant $K^2+\#{\mathcal V}$ of the link $M$. Here
$K$ is the canonical divisor, and $\# {\mathcal V}$ is the number
of irreducible components of the exceptional divisor of the resolution. We want to point out that this invariant   coincides  with an invariant introduced by Gompf in \cite{Gompf} (see Remark \ref{ss: 48}).

The above program has a long history.  M Artin  proved in \cite{Artin,Artin2} that the rational singularities (ie $p_g=0$) can be
characterized  completely from the graph (and he computed even the multiplicity and embedding dimension of these singularities from the graph). In \cite{Lauferme}, H Laufer   extended  these results to minimally elliptic
singularities. Additionally, he noticed that the program breaks for
more complicated singularities (for details, see also section 4).
On the other hand, the first author noticed in \cite{Neminv} that
Laufer's counterexamples do not signal  the end of the program. He
conjectured  that if we restrict ourselves to the case of those
Gorenstein singularities whose links are rational homology spheres
then some numerical analytical invariants (including $p_g$) are
topological. This was carried out explicitly for elliptic
singularities in \cite{Neminv}.

On the other hand, in the literature there is no ``good''
topological candidate for $p_g$ in the very  general case. In fact, we
are searching  for a  {\em ``good'' topological upper bound}
in the following sense.
We want a topological upper bound for $p_g$ for {\em any} normal
surface singularity, which, additionally,  is optimal in the sense
that for Gorenstein singularities it yields exactly $p_g$.
Eg, such a ``good'' topological upper bound
for {\em elliptic} singularities  is
the length of the elliptic sequence, introduced and studied by  S\thinspace S-T Yau
(see, eg \cite{Yau}) and Laufer.

In fact, there are some other particular cases too,
when a possible candidate is present
in the literature. Fintushel and Stern proved in \cite{FS}  that
for a hypersurface Brieskorn  singularity whose link is an {\em
integral} homology sphere, the Casson invariant $\lambda(M)$  of
the link $M$ equals $\si(F)/8$ (hence, by the mentioned
correspondence,  it determines $p_g$ as well). This fact was generalized
by Neumann  and Wahl in \cite{NW}. They proved the same statement
for all Brieskorn--Hamm complete intersections and suspensions of plane
curve singularities (with the same assumption, that the link is an
integral homology sphere). Moreover, they conjectured the validity
of the formula for any isolated complete intersection  singularity
(with the same restriction about the link). For some other relevant 
conjectures the reader can also consult \cite{NW2}.

The result of Neumann--Wahl \cite{NW} was reproved and reinterpreted
by Collin and Saveliev (see \cite{CoS} and \cite{Co}) using equivariant
Casson invariant and cyclic covering techniques. But still, a possible
generalization for rational homology sphere links remained open.
It is important to notice that the ``obvious'' generalization of
the above identity for rational homology spheres, namely to expect
that $\si(F)/8$ equals the Casson--Walker invariant of the link,
completely fails.

In fact, our next conjecture states that one has to replace the
Casson invariant $\lambda(M)$ by a certain Seiberg--Witten
invariant of the link, ie, by the difference of
a certain Reidemeister--Turaev sign-refined torsion
invariant and the Casson--Walker invariant (the sign-change
is motivated by some sign-conventions already used in the literature).

We recall (for details, see section 2 and 3) that the
Seiberg--Witten invariants associates to any $spin^c$ structure
$\si$ of $M$ a rational number $\ssw^0_M(\si)$. In order to
formulate our conjecture, we need to fix a ``canonical'' $spin^c$
structure $\si_{can}$  of  $M$. This can be done as follows. The
(almost) complex structure on $X\setminus\{0\}$  induces a natural
$spin^c$ structure on $X\setminus \{0\}$. Its restriction to $M$
is,  by definition,  $\si_{can}$. The point is that this structure
depends only on the topology of   $M$ alone.

In fact, the $spin^c$  structures correspond in a natural way to
quadratic functions associated with the linking form of $M$; by
this correspondence $\si_{can}$ corresponds to the quadratic
function $-q_{LW}$ constructed by Looijenga and Wahl in \cite{LW}.

We  are now ready to state our conjecture.

{\bf Main Conjecture}\qua {\sl Assume that $(X,0)$ is a
normal surface singularity whose link $M$ is a rational homology
sphere. Let $\si_{can}$ be the canonical $spin^c$  structure on
$M$. Then, conjecturally, the following facts hold.

\vspace{2mm}

\noindent {\rm(1)}\qua For any $(X,0)$, there is a topological upper bound
for $p_g$ given by:
$$\ssw^0_M(\si_{can})-\frac{K^2+\#{\mathcal V}}{8}\geq p_g.$$
{\rm(2)}\qua If $(X,0)$ is $\Q$--Gorenstein, then in (1) one has equality.

{\rm(3)}\qua In particular, if $({\mathcal X},0)$ is a smoothing of a
Gorenstein singularity $(X,0)$ with Milnor fiber $F$, then
$$-\ssw^0_M(\si_{can})=\frac{\si(F)}{8}.$$}

If $(X,0)$ is numerically Gorenstein and $M$ is a $\Z_2$--homology
sphere then $\si_{can}$ is the unique spin structure of
$M$; if $M$ is an integral  homology sphere then in the above
formulae $- \ssw^0_M(\si_{can})=\lambda(M)$, the Casson invariant
of $M$.

In the above Conjecture, we have automatically built in the
following statements as well.

(a)\qua For any normal singularity $(X,0)$ the topological invariant
$$\ssw^0_M(\si_{can})-\frac{K^2+\#{\mathcal V}}{8}$$
is non-negative. Moreover, this topological invariant is zero if
and only if $(X,0)$ is rational. This provides a new topological
characterization of the rational singularities.

(b)\qua Assume that $(X,0)$ (equivalently, the link) is numerically
Gorenstein. Then the above topological invariant is 1 if and only
if $(X,0)$ is minimally elliptic (in the sense of Laufer). Again,
this is a new topological characterization of minimally  elliptic
singularities.

In this paper we will present evidence in support of the conjecture in the form of
explicit verifications. The computations are rather arithmetical,
involving non-trivial identities about generalized Fourier--Dedekind sums.
For the  reader's convenience, we have included a list of basic properties of the Dedekind sums
in Appendix B.

In general it is not easy to compute the Seiberg--Witten invariant.
In our examples we use two different approaches. First,
 the (modified) Seiberg--Witten invariant
 is the sum of the Kreck--Stolz invariant and the number of certain
monopoles \cite{Chen1,Lim,MW}. On the other hand, by a result of
the second author, it can also  be computed as the difference
of the Reidemeister--Turaev torsion and the Casson--Walker invariant
\cite{Nico5} (for more  details, see section 3). Both methods have
their advantages and difficulties.  The first method is rather
explicit when $M$ is a Seifert manifold (thanks to the results of
the second author in \cite{Nico2}, cf also with \cite{MOY}),
but frequently the
corresponding Morse function will be degenerate. Using the second
method, the computation of the Reidemeister--Turaev torsion leads
very often to complicated Fourier--Dedekind sums.

In section 5 we present some formulae for the needed invariant in
terms of the plumbing graph. The formula for the Casson--Walker
invariant was proved by Ratiu in his thesis \cite{Ratiu}, and can
be deduced from Lescop's surgery formulae as well \cite{Lescop}.
Moreover, we also provide a similar formula for the invariant
 $K^2+\#{\mathcal V}$  (which generalizes the corresponding
formula already known for cyclic quotient singularities by
Hirzebruch, see also \cite{Ishii,LW}). The most important result of this
section describes the Reidemeister--Turaev torsion (associated with any
$spin^c$ structure) in terms of the
plumbing graph. The proof is partially based on Turaev's surgery
formulae \cite{Tu} and the structure result  \cite[Theorem 4.2.1]{Tu5}. We have deferred it to  Appendix A.

In our examples we did not try to force the verification of the conjecture in
the largest generality possible, but we tried to supply  a rich and convincing
variety of examples which cover different aspects and cases.

In order to eliminate any confusion about different notations and
conventions in the literature, in most of the cases we provide our
working definitions.

\rk{Acknowledgements}The first author is partially supported by NSF
grant DMS-0088950; the second author is partially supported by NSF
grant DMS-0071820.

\section{The link and its canonical $spin^c$ structure}

\gtsubsection{Definitions}
\label{ss: 2.1}
Let $(X,0)$ be a normal surface singularity
embedded in  $(\C^N,0)$. Then for $\epsilon$ sufficiently small
the intersection $M:=X\cap S_{\epsilon}^{2N-1}$ of a
representative $X$ of the germ with the sphere
$S_{\epsilon}^{2N-1}$ (of radius $\eps$) is a compact oriented
3--manifold, whose oriented $C^{\infty}$ type does not depend on
the choice of the embedding and $\epsilon$. It is called the link
of $(X,0)$.

In this article we will assume that $M$ is a rational homology
sphere, and we write $H:=H_1(M,\Z)$. By Poincar\'e duality $H$ can
be identified with $H^2(M,\Z)$.

It is well-known that $M$ carries a symmetric non-singular
bilinear  form
$$b_M\co H\times H\to \Q/\Z$$
called the linking form of $M$. If $[v_1]$ and $[v_2]\in H$ are
represented by the 1--cycles $v_1$ and $v_2$, and for some integer
$n$ one has $nv_1= \partial w$, then $b_M([v_1],[v_2])=(w\cdot v_2)/n\
(\mmod \ \Z)$.

\gtsubsection{The linking form as discriminant form}
\label{ss: 2.2}
We briefly recall
 the definition of the discriminant form. Assume that $L$
is a finitely generated free Abelian  group with a symmetric
bilinear  form $\bif\co L\times L\to \Z$. Set $L' :=\Hom_{\Z}(L,\Z)$.
 Then there is a natural homomorphism $i_L\co L\to L' $ given by $x\mapsto
(x,\cdot)$ and a natural extension of the form $\bif$ to a
rational bilinear form $\bif_\Q\co L' \times L' \to \Q$. If
$d_1,d_2\in L' $ and $nd_j=i_L(e_j)$ ($j=1,2$)
 for some integer
$n$, then  $(d_1,d_2)_\Q=d_2(e_1)/n=(e_1,e_2)/n^2$.

If $L$ is non-degenerate (ie, $i_L$ is a monomorphism) then one
defines the discriminant space $D(L)$ by  coker$(i_L)$. In this
case there is a discriminant bilinear form
$$b_{D(L)}\co D(L)\times D(L)\to \Q/\Z$$
defined by $b_{D(L)}([d_1],[d_2])=(d_1,d_2)_\Q \ (\mmod\ \Z)$.

Assume that $M$ is the boundary of a oriented 4--manifold $N$ with
$H_1(N,\Z)=0$ and $H_2(N,\Z)$ torsion-free. Let $L$ be the
intersection lattice $\bigl(H_2(N,\Z),\bif\bigr)$. Then $L' $ can be
identified with $H_2(N,\partial N,\Z)$ and one has the exact sequence
$L\to L' \to H\to 0$.
 The fact that $M$ is a rational homology sphere implies that $L$ is
non-degenerate. Moreover $(H,b_M)=(D(L),-b_{D(L)})$. Sometimes it
is also convenient to regard $L' $ as $H^2(N,\Z)$ (identification
by Poincar\'e duality).

\gtsubsection{Quadratic functions and forms associated with $b_M$}
\label{ss: 2.3} A
map $q\colon H\to \Q/\Z$ is called quadratic {\em function} if
$b(x,y)=q(x+y)-q(x)-q(y)$ is a bilinear form on $H\times H$. If in
addition $q(nx)=n^2q(x)$ for any $x\in H$ and $n\in \Z$ then $q$
is called quadratic {\em form}. In this case we say that the
quadratic function, respectively form, is associated with $b$.
Quadratic {\em forms}
 are also called {\em quadratic refinements}  of the bilinear form $b$.

In the case of the link $M$, we denote by $Q^c(M)$ (respectively
by $Q(M)$) the set of quadratic  functions (resp. forms)
associated with $b_M$. Obviously, there is a natural inclusion
$Q(M)\subset Q^c(M)$.

  The set $Q(M)$ is non-empty.  It is a
$G:=H^1(M,\Z_2)$ torsor,  ie, $G$ acts freely and transitively on
$Q(M)$. The action can be easily described if we identify $G$ with
$\Hom(H,\Z_2)$ and we regard $\Z_2$ as $(\frac{1}{2}\Z)/\Z\subset
\Q/\Z$. Then, the difference of any two quadratic refinements of
$b_M$ is an element of $G$, which provides a natural action
$G\times Q(M)\to Q(M)$ given by $(\chi,q)\mapsto \chi+q$.

Similarly, the set $Q^c(M)$ is non-empty and it is a
$\hat{H}=\Hom(H,\Q/\Z)$ torsor. The free and transitive action
$\hat{H}\times Q^c(M)\to Q^c(M)$ is given by the same formula
$(\chi,q)\mapsto \chi+q$. In particular,  the inclusion
$Q(M)\subset Q^c(M)$ is $G$--equivariant  via the natural
monomorphism $G\hookrightarrow \hat{H}$. We prefer to replace the
$\hat{H}$ action on $Q^c(M)$ by an action of $H$. This action
$H\times Q^c(M)\to Q^c(M)$ is defined by $(h,q)\mapsto
q+b_M(h,\cdot)$. Then the natural monomorphism $G\hookrightarrow
\hat{H}$ is replaced by the Bockstein-homomorphism
$G=H^1(M,\Z_2)\to H^2(M,\Z)=H$. In the sequel we consider $Q^c(M)$
with this $H$--action.

Quadratic functions appear in a natural way. In order to see this,
 let $N$ be as in \ref{ss: 2.2}.
Pick a characteristic element, that is an element  $k\in L' $,  so
that $ (x,x)+k(x)\in2\Z$ for  any $x\in L$. Then for any  $d\in L'
$ with class $[d]\in H$, define
$$
q_{D(L),k}([d]):= \frac{1}{2}\, (d+k,d)_\Q\ (\mmod \ \Z).
$$
 Then $-q_{D(L),k}$  is a quadratic function associated with
$b_{M}=-b_{D(L)}$. If in addition $k\in \im(i_L)$, then
$-q_{D(L),k}$  is a {\em quadratic  refinement}  of  $b_{M}$.

There are two important examples to consider.

First assume that $N$ is an {\em almost-complex } manifold, ie,
its tangent bundle $TN$ carries an almost  complex structure. By
Wu formula, $k=-c_1(TN)\in L'$ is a characteristic element. Hence
$-q_{D(L),k}$ is a quadratic function associated with $b_M$. If
$c_1(TN)\in \im(i_L)$ then we obtain a quadratic refinement.

Next, assume that $N$ carries a {\em spin structure}. Then
$w_2(N)$ vanishes, hence by Wu formula $\bif$ is an even form.
Then one can take $k=0$, and $-q_{D(L),0}$ is a quadratic
refinement of $b_M$.

\gtsubsection{The $spin$ structures of $M$}
\label{ss: 2.4}  The 3--manifold $M$ is always spinnable. The set $Spin(M)$ of the possible spin structures of
$M$ is a $G$--torsor. In fact, there is a natural (equivariant)
identification of $\gq:Spin(M)\to Q(M)$.

In order to see this, fix a spin structure $\epsilon\in Spin(M)$.
Then there exists a simple connected  oriented spin 4--manifold $N$
with $\partial N=M$ whose induced spin structure on $M$ is exactly
$\epsilon$ (see, eg \cite[5.7.14]{GS}). Set
$L=(H_2(N,\Z),\bif)$. Then the quadratic refinement $-q_{D(L),0}$
(cf \ref{ss: 2.3}) of $b_M$ depends only on the spin structure
$\epsilon$ and not on the particular choice of $N$. The
correspondence $\epsilon\mapsto -q_{D(L),0}$ determines the
identification $\gq$ mentioned above.

\gtsubsection{The $spin^c$ structures on $M$}
\label{ss: 2.4b} We denote by $Spin^c(M)$ the space of isomorphism classes of $spin^c$
structures on $M$. $Spin^c(M)$ is in a natural way a
$H=H^2(M,\Z)$--torsor.
 We denote this action $H\times Spin^c(M)\ra Spin^c(M)$
 by $(h,\si)\mapsto h\cdot\si$.
For every $\si \in Spin^c(M)$ we denote by ${\bS}_\si$ the
associated  bundle of complex spinors, and by  $\det\si$ the
associated line bundle, $\det\si:=\det {\bS}_\si$.  We set
$c(\si):=c_1({\bS}_\si)\in H$. Note that $c({h\cdot\si})=
2h+c({\si})$.

$Spin^c(M)$ is equipped with a natural involution
$\si\longleftrightarrow\bar{\si}$ such that
\[
c(\bar{\si})=-c(\si)\ \ \mbox{and} \ \ \ \overline{h\cdot \sigma}=
(-h)\cdot \bar{\si}.
\]
There is a natural injection $\epsilon \mapsto \si(\epsilon)$  of
$Spin(M)$ into   $Spin^c(M)$. The image of   $Spin(M)$ in
$Spin^c(M)$  is
\[
\{ \si \in Spin^c(M);\;\;  c(\si)=0\}=\{\si \in
Spin^c(M);\;\;\si=\bar{\si}\}.
\]
Consider now a 4--manifold $N$ with lattices  $L$ and $L'$ as in
\ref{ss: 2.2}. We prefer to write $L' =H^2(N,\Z)$, and denote by $d\mapsto
[d]$  the restriction  map $L'=H^2(N,\Z)\to H^2(M,\Z)=H$.

Then $N$ is automatically a $spin^c$ manifold. In fact, the set of
$spin^c$ structures on $N$ is  parametrized by the set of
characteristic elements
$${\mathcal C}_{N}:=\{k\in L' :\, k(x)+(x,x)\in 2\Z\ \mbox{for
all } \ x\in L\}$$ via $\tilde{\sigma}\mapsto c(\tilde{\si})\in
{\mathcal C}_N$ (see. eg \cite[2.4.16]{GS}). The set
$Spin^c(N)$ is an $L' $ torsor with action
$(d,\tilde{\si})\mapsto d\cdot \tilde{\si}$. Let $r\co Spin^c(N)\to
Spin^c(M)$ be the restriction. Then $r(d\cdot
\tilde{\si})=[d]\cdot r(\tilde{\si})$ and
$c(r(\tilde{\si}))=[c(\tilde{\si})]$.

Moreover, notice that $r(\tilde{\si})=r(d\cdot \tilde{\si})$ if
and only if $[d]=0$, ie $d\in L$. If this is happening then
$c(\tilde{\si})- c(d\cdot \tilde{\si})\in 2L$.

\gtsubsection{Lemma}
\label{ss: 2.4c} {\sl There is a canonical $H$--equivariant identification
\[
\gq^c\co  Spin^c(M)\to Q^c(M).
\]
 Moreover, this identification is compatible with the $G$--equivariant identification $\gq\co Spin(M)\to Q(M)$ via the inclusions $Spin(M)\subset Spin^c(M)$ and $Q(M)\subset Q^c(M)$.}

\begin{proof} Let $N$ be as above. We first show that $r$
is onto. Indeed, take any $\tilde{\si}\in Spin^c(N)$ with
restriction $\si\in Spin^c(M)$. Then all the elements in the
$H$--orbit of $\si$ are induced structures. But this orbit is the
whole set. Next, define for any $\tilde{\si}$ corresponding to
$k=c(\tilde{\si})$ the quadratic function $q_{D(L),k}$. Then
$r(\tilde{\si})=r(d\cdot \tilde{\si})$ if and only if $d\in L$.
This means that $c(\tilde{\si})-c(d\cdot \tilde{\si})\in 2L$ hence
$c(\tilde{\si})$ and $c(d\cdot \tilde{\si})$ induce the  same
quadratic function.  Hence
$\gq^c(r(\tilde{\si})):=-q_{D(L),c(\sigma)}$ is well-defined.
Finally, notice that $\gq^c$ does not depend on the choice of $N$,
fact which shows its compatibility with $\gq$ as well (by taking
convenient spaces $N$). \end{proof}

\gtsubsection{$M$ as a plumbing manifold}
\label{ss: 2.5}
\setcounter{equation}{0}
Fix a sufficiently small (Stein) representative  $X$ of $(X,0)$ and let $\pi\co \tilde{X}\to
X$ be a resolution of the singular point $0\in X$. In particular,
$\tilde{X}$ is smooth, and $\pi$ is a biholomorphic  isomorphism
above $X\setminus \{0\}$. We will assume that the exceptional
divisor $E:=\pi^{-1}(0)$ is a normal crossing divisor with
irreducible components $\{E_v\}_{v\in {\mathcal V}}$. Let
$\Gamma(\pi)$ be the  dual resolution graph associated with $\pi$
decorated with the self intersection numbers $\{(E_v,E_v)\}_v$.
Since $M$ is a rational homology sphere, all the irreducible
components $E_v$ of $E$ are rational, and $\Gamma(\pi)$ is a tree.

It is clear that $H_1(\tilde{X},\Z)=0$ and $H_2(\tilde{X},\Z)$ is
freely generated by the fundamental classes $\{[E_v]\}_v$. Let $I$
be the intersection matrix $\{(E_v,E_w)\}_{v,w}$. Since $\pi $
identifies $\partial \tilde{X}$ with $M$, the results from \ref{ss: 2.2}
can be applied. In particular, $H=\coker(I)$ and $b_M=-b_{D(I)}$.
The matrix $I$ is negative definite.

The graph $\Gamma(\pi)$ can be identified with a plumbing graph,
and $M$ can be considered as an $S^1$--plumbing manifold whose
plumbing graph is $\Gamma(\pi)$. In particular, any resolution
graph $\Gamma(\pi)$ determines the oriented 3--manifold $M$
completely.

We say that two plumbing graphs (with negative definite
intersection forms) are equivalent  if one of them can be obtained
from the other by a finite sequence of blowups and/or blowdowns along rational $(-1)$--curves. Obviously, for a given
$(X,0)$, the resolution $\pi$, hence the graph $\Gamma(\pi)$ too,
is not unique. But different resolutions provide equivalent
graphs. By a result of W. Neumann \cite{NP}, the oriented diffeomorphism
type of $M$ determines completely the equivalence class of
$\Gamma(\pi)$. In particular, any invariant defined from the
resolution graph $\Gamma(\pi)$ (which is constant in its
equivalence class) is, in fact, an invariant of the oriented
$C^{\infty}$ 3--manifold $M$. This fact will be crucial in the next
discussions.

Now, we fix a resolution $\pi$ as above and identify $M=\partial
\tilde{X}$.
 Let $K$ be the canonical class (in $Pic(\tilde{X})$) of  $\tilde{X}$.
By the adjunction formula,
\[
-K\cdot E_v=E_v\cdot E_v+ 2
\]
 for any $v\in {\mathcal V}$. In fact, $K$ at homological level provides an
element $k_{\tilde{X}}\in L' $ which has the obvious property
\[
-k_{\tilde{X}}([E_v])=([E_v],[E_v])+2
\ \ \ \mbox{for any $v\in {\mathcal V}$}.
\]
Since the matrix $I$ is non-degenerate, this defines $k_{\tilde{X}}$ uniquely.

 $-k_{\tilde{X}}$ is known in the literature as the \emph{canonical
(rational) cycle} of $(X,0)$ associated with the resolution $\pi$.
More precisely, let $Z_K=\sum_{v\in {\mathcal V}}r_vE_v$, $r_v\in\Q$, be a
 rational cycle supported  by the exceptional divisor $E$, defined by
\begin{equation*}
Z_K\cdot E_v=-K\cdot E_v=E_v\cdot E_v+2 \ \ \
\mbox{for any $v\in {\mathcal V}$}.
\tag{$*$}
\label{2.5*}
\end{equation*}
Then the above linear system  has a unique solution, and
$\sum_vr_v[E_v] \in L\otimes \Q$ can be identified with
$(i_L\otimes \Q)^{-1}(-k_{\tilde{X}})$.

It is clear that $-k_{\tilde{X}}\in \im(i_L)$ if and only if all
the coefficients $\{r_v\}_v$ of $Z_K$  are integers. In this case
the singularity $(X,0)$ is called {\em numerically Gorenstein}
(and we will also say that ``$M$ is numerically Gorenstein'').

In particular, for any normal singularity $(X,0)$, the resolution
$\pi$  provides a quadratic function  $-q_{D(I),k_{\tilde{X}}}$
associated with  $b_M$, which is a quadratic form  if and only if
$(X,0)$ is numerically Gorenstein.

\gtsubsection{The universal property of $q_{D(I),k_{\tilde{X}}}$}
\label{ss: 2.6}
\setcounter{equation}{0} In
\cite{LW}, Looijenga and Wahl define a quadratic function
 $q_{LW}$ (denoted by $q$ in  \cite{LW}) associated  with $b_M$ from
the almost complex structure of the bundle $TM\oplus \R_M$ (where
$TM$ is the tangent bundle and $\R_M$ is the trivial bundle of
$M$). By the main universal property of $q_{LW}$ (see [loc. cit.],
Theorem 3.7) (and from the fact that any resolution $\pi$ induces
the same almost complex structure on $TM\oplus \R_M$) one gets
that for any $\pi$ as in \ref{ss: 2.5},  the identity
$q_{LW}=-q_{D(I),k_{\tilde{X}}}$ is valid.  This shows that
$q_{D(I),k_{\tilde{X}}}$ does not depend on the choice of the
resolution $\pi$.

This fact can be verified by elementary computation as well: one
can prove that $q_{D(I),k_{\tilde{X}}}$ is stable  with respect to
a blow up (of points of $E$).

\gtsubsection{The ``canonical'' $spin^c$ structure of a singularity
link}
\label{ss: 2.7}
\setcounter{equation}{0}
 Assume that $M$ is the link of $(X,0)$. Fix a
resolution $\pi\co \tilde{X}\to X$ as in \ref{ss: 2.5}. Then $\pi$
determines a ``canonical'' quadratic function
$q_{can}:=-q_{D(I),k_{\tilde{X}}}$ associated with  $b_M$ which
does not depend on the choice of $\pi$ (cf \ref{ss: 2.6}). Then the
natural  identification $\gq^c\co  Spin^c(M)\to Q^c(M)$ (cf
\ref{ss: 2.3}) provides  a well-defined $spin^c$ structure
$(\gq^c)^{-1}(q_{can})$.  Then the ``canonical'' $spin^c$
structure $\si_{can}$
 on $M$ is $(\gq^c)^{-1}(q_{can})$ modified by the
natural involution of $Spin^c(M)$. In particular, $c(\si_{can})=
-[k_{\tilde{X}}]\in H$.
(Equivalently, $\si_{can}$ is the restriction to $M$ of the $spin^c$ structure
given by the characteristic element
 $-k_{\tilde{X}}\in {\mathcal C}_{\tilde{X}}$.)
If $(X,0)$ is numerically Gorenstein then $\si_{can}$  is a spin
structure. In this case we will use the notation
$\eps_{can}=\si_{can}$ as well.

 We want to emphasize (again) that $\si_{can}$ depends only on the oriented $C^{\infty}$ type of
$M$ (cf also with \ref{ss: 2.8}). Indeed, one can construct
$q_{can}$ as follows. Fix an arbitrary plumbing graph
$\Gamma$ of $M$  with negative definite intersection form
(lattice) $L$. Then determine $Z_K$ by \ref{ss: 2.5}(\ref{2.5*}), and take
\[
q_{can}([d]):=-\frac{1}{2} (d-[Z_K],d)_\Q \;(\mmod \ \Z).
\]
Then $q_{can}$ does not depend on the choice of $\Gamma$.

It is remarkable that this construction provides an ``origin'' of the torsor space $Spin^c(M)$.

\gtsubsection{Compatibility with the (almost) complex structure}
\label{ss: acs}
\setcounter{equation}{0}
As we have already mentioned in \ref{ss: 2.6}, the result of Looijenga and
Wahl \cite{LW} implies the following:
the almost  complex structure on $X\setminus\{0\}$
determines a $spin^c$ structure, whose restriction to $M$ is $\si_{can}$.
Similarly, if $\pi\co \tilde{X}\to X$ is a resolution, then
the almost complex structure on $\tilde{X}$  gives a $spin^c$ structure
$\si_{\tilde{X}}$ on $\tilde{X}$,  whose restriction to $M$ is $\si_{can}$.
Here we would like to add the following discussion. Assume that
the intersection form $\bif_{\tilde{X}}$ is even, hence
$\tilde{X}$ has a unique spin structure $\eps_{\tilde{X}}$. The
point is that, in general, $\si_{\tilde{X}}\not=\eps_{\tilde{X}}$,
and their restrictions can be different as well, even if the
restriction of $\si_{\tilde{X}}$ is spin.

More precisely: $\bif_{\tilde{X}}$ is even if and only if
$k_{\tilde{X}}\in 2L'$; $r(\si_{\tilde{X}})\in Spin(M)$ if and
only if $k_{\tilde{X}}\in L$; and finally, $r(\si_{\tilde{X}})
=r(\eps_{\tilde{X}}) $ if and only if $k_{\tilde{X}}\in 2L$.

\gtsubsection{Remarks}
\label{ss: 2.8}
\setcounter{equation}{0}

(1)\qua  In fact, by the classification theorem of plumbing graphs given by
Neumann \cite{NP}, if $M$ is a rational homology sphere which is
not a lens space, then already $\pi_1(M)$ (ie, the homotopy type
of $M$) determines its orientation class and its canonical
$spin^c$ structure. Indeed, if one wants to recover the oriented
$C^{\infty}$ type of $M$ from its fundamental group, then by
Neumann's result the only ambiguity  appears for cusp
singularities (which are not rational homology spheres) and for
cyclic quotient singularities. The links of cyclic quotient
singularities are exactly the lens spaces. In fact, if we assume
the numerically Gorenstein assumption, even the lens spaces are
classified by their fundamental groups (since they are exactly the
du Val $A_p$--singularities).

(2)\qua  If $M$ is a numerically Gorenstein $\Z_2$--homology  sphere,
the definition of $\epsilon_{can}$ is obviously simpler: it is the
unique spin structure of $M$. If $M$ is an {\em integral }
homology sphere then it is automatically numerically Gorenstein,
hence the above statement applies.

(3)\qua  Assume that $(X,0)$ has a smoothing  with Milnor fiber $F$
whose homology group $H_1(F,\Z)$ has no torsion. Then the (almost)
complex structure of $F$ provides a $spin^c$ structure on $F$
whose restriction to $M$ is exactly $\si_{can}$. This follows
(again) by the universal property of $q_{LM}$ (\cite[Theorem 3.7]{LW}, ; cf also with \ref{ss: 2.6}).

Moreover, $F$ has a spin structure if and only if
 the intersection form $\bif$ of $F$ is even
 (see eg \cite[5.7.6]{GS}); and in this case,
the spin structure is unique. If $F$ is  spin, then its spin
structure $\eps_F$
 coincides with the $spin^c$ structure induced by the complex
structure (since the canonical bundle of $F$ is trivial). In
particular, if $F$ is spin, $\si_{can}$ is the restriction of
$\eps_F$, hence it is spin. This also proves that if $(X,0)$ has a
smoothing with even intersection form and without torsion in
$H_1(F,\Z)$, then it is necessarily numerically Gorenstein.

Here is worth  noticing   that the Milnor fiber of a smoothing of
a Gorenstein singularity has  even intersection form \cite{Seade}.

(4)\qua Clearly, $q_{can}$  depends only on $Z_K$ ($\mmod \ 2\Z$).

\gtsubsection{The invariant $K^2+\#{\mathcal V}$}
\label{ss: 2.11}
\setcounter{equation}{0}
Fix a resolution
$\pi\co \tilde{X} \to X$ of $(X,0)$ as in \ref{ss: 2.5}, and consider
$Z_K$ or $k_{\tilde{X}}$. The rational number $Z_K\cdot
Z_K=(k_{\tilde{X}},k_{\tilde{X}})_\Q$ will be denoted by $K^2$.
Let $\#{\mathcal V}$ denote the number of irreducible components
of $E=\pi^{-1}(0)$. Then $K^2+\#{\mathcal V}$ does not depend on
the choice  of the resolution $\pi$. In fact, the discussion in
\ref{ss: 2.5} and \ref{ss: 2.7} shows that it is an invariant of $M$.
Obviously, if $(X,0)$ is numerically Gorenstein, then
$K^2+\#{\mathcal V}\in \Z$.

\gtsubsection{Notation}
\label{ss: no}
\setcounter{equation}{0}Let $\tilde{X}$ as above. Let $\{E_v\}_{v\in {\mathcal V}}$ be the set of irreducible
exceptional divisors and $D_v$ a  small transversal disc to $E_v$.
Then $\{[E_v]\}_v$ (resp. $\{[D_v]\}_v$) are the free generators
of $L=H_2(\tilde{X},\Z)$ (resp. $L' =H_2(\tilde{X},M,\Z)$) with
$[D_v]\cdot [E_w]=1$ if $v=w$ and $=0$ otherwise. Moreover,
$g_v:=[\partial D_v]$ ($v\in {\mathcal V}$)
 is a generator set of $L' /L=H$. In fact $\partial D_v$ is a generic
fiber of the $S^1$--bundle over $E_v$ used in the plumbing
construction of $M$. If $I$ is the intersection matrix defined by
the resolution (plumbing) graph, then $i_L$ written in the bases
$\{[E_v]\}_v$ and $\{[D_v]\}_v$ is exactly $I$.

Using this notation, $k_{\tilde{X}}\in L'$ can be expressed as
$\sum_v(-e_v-2)[D_v]$, where $e_v=E_v\cdot E_v$. For the degree of $v$  (ie, for $\#\{w: E_w\cdot E_v=1\}$) we
will use the notation $\delta_v$. Obviously
\[
\sum_v\delta_v=-2\times \mbox{Euler characteristic of the plumbing graph}+2\#\cv=2\#\cv -2.
\]
Most of the examples considered later are star-shaped graphs. In
these cases it is convenient to express the corresponding
invariants of the Seifert 3--manifold $M$ in terms of their
Seifert invariants. In order to eliminate any confusion about the
different notations and conventions in the literature, we list
briefly the definitions and some of the needed properties.

\gtsubsection{The unnormalized Seifert invariants}
\label{ss: s1}
\setcounter{equation}{0} Consider a Seifert fibration $\pi\co M\to \Sigma$. In
our situation   $M$ is a rational homology sphere and  the base
space $\Sigma$ is an $S^2$ with genus 0 (and we will not emphasize this
fact anymore).

Consider a set of points $\{x_i\}_{i=1}^{\nu}$ in such a way that
the set of fibers $\{\pi^{-1}(x_i)\}_i$ contains the set of
singular fibers. Set $O_i:= \pi^{-1}(x_i)$. Let $D_i$ be a
small disc  in $X$ containing  $x_i$,
$\Sigma':=\Sigma\setminus \cup_iD_i$ and $M':=\pi^{-1}(\Sigma')$.  Now,
$\pi\co M'\to \Sigma'$ admits sections, let $s\co \Sigma'\to M'$ be one of them.
Let $Q_i:=s(\partial D_i)$ and let $H_i$  be  a circle fiber in
$\pi^{-1}(\partial D_i)$. Then in $H_1(\pi^{-1}(D_i),\Z)$ one has
$H_i\sim \alpha_iO_i$ and $Q_i\sim -\beta_iO_i$ for some integers
$\alpha_i>0 $ and $
 \beta_i$ with $(\alpha_i,\beta_i)=1
$. The set $((\alpha_i,\beta_i)_{i=1}^\nu)$ constitute the set of
\emph{(unnormalized) Seifert invariants}. The number
$$
e:=- \sum_i(\beta_i/\alpha_i)
$$
 is called the \emph{(orbifold) Euler number} of $M$. If $M$ is a link of singularity then $e<0$.

Replacing the section by another one, a different choice  changes
each $\beta_i$ within its residue class modulo $\alpha_i$ in such
a way that the sum $e=-\sum_i(\beta_i/\alpha_i)$ is constant.

The elements $q_i=[Q_i] \ (1\leq i\leq \nu)$ and the class $h$ of
the generic fiber $H$ generate the group $H=H_1(M,\Z)$. By the
above construction is clear that:
\[
H=\mbox{ab}\langle q_1,\ldots q_\nu,h\ |\ q_1\ldots q_\nu=1,
\ q_i^{\alpha_i}h^{\beta_i}=1,\ \mbox{for all $i$} \rangle.
\]
 Let
$\alpha:=\mbox{lcm}(\alpha_1,\ldots,\alpha_\nu)$. The order of the
group $H$ and of the subgroup $\langle h \rangle$ can be
determined by (cf \cite{Neu}):
\[
|H|=\alpha_1\cdots \alpha_\nu |e|, \ \ \ \ |\langle h\rangle|=\alpha|e|.
\]

\gtsubsection{The normalized Seifert invariants and plumbing graph}
\label{ss: s2}
\setcounter{equation}{0}
We write
$$
e=b+\sum \omega_i/\alpha_i
$$
for some integer $b$, and $0\leq \omega_i<\alpha_i$ with
$\omega_i\equiv  -\beta_i \ (\mmod \ \alpha_i)$.  Clearly, these
properties define $\{\omega_i\}_i$ uniquely. Notice that $b\leq
e<0$. For the uniformity of the notations, in the sequel we assume
$\nu\geq 3$.

For each $i$, consider the continued fraction $\alpha_i/ \omega_i=
b_{i1}-1/(b_{i2}-1/(\cdots
-1/b_{i\nu_i})\cdots)$. Then (a
possible) plumbing graph  of $M$ is a star-shaped graph  with
$\nu$ arms. The central vertex has decoration $b$ and the arm
corresponding to
 the index $i$ has $\nu_i$ vertices, and they are
decorated by $b_{i1},\ldots, b_{i\nu_i}$ (the vertex decorated by
$b_{i1}$ is connected by the central  vertex).

We will distinguish those vertices $v\in {\mathcal V}$  of the
graph which have $\delta_v\not= 2$. We will denote by $\bar{v}_0$
the central vertex (with $\delta=\nu$), and by $\bar{v}_i$ the
end-vertex of the $i^{th}$ arm (with $\delta=1$) for all $1\leq
i\leq \nu$.
In this notation, $g_{\bar{v}_0}=h$, the
class of the generic fiber. Moreover, using the plumbing
representation of the group $H$, we have another presentation for
$H$, namely:
\[
H=\mbox{ab}\langle g_{\bar{v}_1},\ldots g_{\bar{v}_\nu},h\ |\
h^{-b}=\prod_{i=1}^\nu g_{\bar{v}_i}^{\omega_i},
 \ \ h=g_{\bar{v}_i}^{\alpha_i}
\  \mbox{for all $i$} \rangle.
\]

%\end{document}

\section{Seiberg--Witten invariants of ${\bQ}$--homology spheres}

\noindent In this section we consider an oriented   rational
homology $3$--sphere
 $M$.  We set $H:=H^2(M,\Z)$. When working with the group algebra
$\Q[H]$ of $H$ it is more convenient to use the
multiplicative notation for the group operation of $H$.

\gtsubsection{The Seiberg--Witten invariants of $M$}
\setcounter{equation}{0}
\label{ss: 3.2} To describe the
Seiberg--Witten invariants we need to fix some additional geometric
data belonging to  the space of  parameters
\[
\p=\{ (g,\eta);\quad g=\mbox{Riemann metric},\;\;\eta=\mbox{closed
two-form} \}.
\]
For each $spin^c$ structure $\si$ on $M$  (cf  \ref{ss: 2.4b}), we
have the space of configurations   $\cC_\si$ (associated with
$\si$) consisting  of pairs ${\sf C}=(\psi, A)$, where $\psi$ is a
section of ${\bS}_\si$ and $A$ is a Hermitian connection on $\det
\si$.  The gauge   group $\gauge:= {\rm Map}\,(M, S^1)$ acts on
$\cC_\si$.    Moreover, it acts freely on the irreducible part
\[
\cC_\si^{irr}=\{ (\psi, A)\in \cC_\si;\quad \psi\not\equiv 0\},
\]
and the quotient   $\B_\si^{irr}:=\cC_\si^{irr}/\gauge$   can be
equipped with a structure of Hilbert manifold. Every parameter
$u=(g,\eta)\in \p$ defines a $\gauge$--invariant function \linebreak $\f_{\si, u}\co  \cC_\si\ra {\bR} $ whose  critical points  are called the
$(\si,g,\eta)$--{\em Seiberg--Witten monopoles}.
  In particular, $\f_{\si, u}$ descends  to  a smooth function
$[\f_{\si, u}]\co \B_\si^{irr}\ra {\bR} $. We denote by
$\modu_{\si, u}^{irr}$ its critical set.

The  first Chern class $c(\si)$ of ${\bS}_\si$ is a torsion
element of $H^2(M,\Z)$,  and thus the curvature of any connection
on $\det\si$ is an {\em exact} $2$--form.  In particular we can
find an unique $\gauge$--equivalence class of  connections  $A$ on
$\det \si$ with the property that
$$
F_A=\ii\eta. \eqno{(\dagger)}
$$
Using  the metric $g$ on $M$ (which is part of our parameter $u$)
and  a connection $A_u$  satisfying $(\dagger)$,  we obtain a
$spin^c$--Dirac operator $\dir_{A_u}$. To define the Seiberg--Witten
invariants we need to work with {\em good} parameters, ie,
parameters $u$ such that the following two things happen.

\vspace{1mm}

\noindent $\bullet$\qua  The Dirac operator $\dir_{A_u}$ is invertible.

$\bullet$\qua The function $[\f_{\si, u}]$ is Morse, and $\modu_{\si,u}$ consists of finitely many points.

\vspace{1mm}

The space of good parameters is {\em generic}. Fix such a good
parameter $u$. Then   each critical point has a well defined
$\Z_2$--valued Morse index
\[
\bmm\co  \modu^{irr}_{\si, u}\ra \{\pm 1\}
\]
and we set
$$
\ssw_M(\si,u)=\sum_{x\in \modu_{\si,u}^{irr}}\bmm(x)\in {\bZ}.
$$
This integer depends on the choice of the parameter $u$ and thus
it is not a topological  invariant.   To obtain an invariant we
need to  alter this monopole count.

The  eta invariant of $\dir_{A_u}$  depends only on  the gauge
equivalence class of  $A_u$, and we will denote it by
$\eta_{dir}(\si,u)$. The metric $g$ defines an odd signature
operator  on $M$ whose eta invariant we denote by
$\eta_{sign}(u)$.
  Now define the  {\em Kreck--Stolz} invariant associated with the data
$(\si,u)$ by
\[
KS_M(\si, u):= 4\eta_{dir}(\si,u) +\eta_{sign}(u) \in {\bQ}.
\]
We have the following result.
\gtsubsection{Theorem}
\label{ss: 3.3}{\cite{Chen1,Lim,MW}}\qua {\sl The rational number
\[
\frac{1}{8}KS_M(\si, u)+\ssw_M(\si, u)
\]
is independent   of $u$ and thus  it is a topological invariant of
the pair $(M,\si)$. We denote this number by \ $\ssw_M^0(\si)$.
Moreover
\begin{equation*}
\ssw_M^0(\si)=\ssw^0_M(\bar{\si}). \tag{$*$}
\label{3.3*}
\end{equation*} }

\noindent It is convenient to rewrite the collection
$\{\ssw^0_M(\si)\}_\si$ as a function $H\to \Q$ (see eg
the Fourier calculus below). For every $spin^c$
structure $\si$  on $M$  we consider
\[
\bSW^0_{M,\si}:=\sum_{h\in H} \ssw^0_M(h^{-1}\cdot \si)h\in
{\bQ}[H].
\]
Equivalently, $\bSW^0_{M,\si}$, as a function $H\ra {\bQ}$, is
defined by $\bSW^0_{M,\si}(h)= \ssw^0_M(h^{-1}\cdot \si)$. The
symmetry condition \ref{ss: 3.3}(\ref{3.3*})  implies
\[
\bSW^0_{M,\si}(h)=\bSW^0_{M,\bar{\si}}(h^{-1})\ \ \ \mbox{for all}
\ h\in H.
\]
This description is very difficult to use in concrete computations
unless we have very specific information about the geometry of
$M$. This is the case of the Seifert 3--manifolds,  see
\cite{Nico0, Nico2} for the complete presentation. In the next subsection
we recall some facts needed in our computations. The
interested reader is invited to consult [loc. cit.] for more
details.

\gtsubsection{The Seiberg--Witten invariants of Seifert manifolds}
\label{ss: sws}  We
will use the notations of \ref{ss: s1} and \ref{ss: s2}; nevertheless,
in \cite{FuS, MOY, Nico2} (and in general, in the gauge theoretic literature) some other notations became
generally accepted too. They will be mentioned  accordingly.

In \cite{MOY, Nico2} a Seifert manifold is regarded as   the unit
circle sub-bundle of an (orbifold) $V$--line bundle over a
$2$--dimensional $V$--manifold (orbifold) $\Sigma$.  The
$2$--dimensional orbifold in our case is  ${\bP}^1$ (with $\nu$
conical singularities each with  angle $\frac{\pi}{\alpha_i}$,
$i=1,\ldots,\nu$).

The  space of isomorphisms classes  of topological $V$--line
bundles over $\Sigma$ is an Abelian group ${\rm
Pic}^V_{top}(\Sigma)$. Its is a subgroup of ${\bQ}\times
\prod_{i=1}^\nu\, \Z_{\alpha_i}$, and correspondingly we denote
its elements by $(\nu+1)$--uples
\[
L(c; \frac{\tau_1}{\alpha_1},\cdots ,\frac{\tau_\nu}{\alpha_\nu}),
\]
where $0\leq \tau_i <\alpha_i$, $i=1,\cdots,\nu$. The number $c$
is called the \emph{rational  degree},   while the   fractions
$\frac{\tau_i}{\alpha_i}$ are called the  \emph{singularity data}.
They are subject to a single  compatibility condition
\[
c-\sum_{i=1}^\nu\frac{\tau_i}{\alpha_i}\in {\bZ}.
\]
To  any $V$--line bundle $L(c; \frac{\tau_i}{\alpha_i},\;1\leq
i\leq \nu)$   we canonically associate  a \emph{smooth} line
bundle $|L|\ra \Sigma={\bP}^1$ uniquely  determined by the
condition
\[
\deg|L|=c-\sum_{i=1}^\nu\frac{\tau_i}{\alpha_i}.
\]
The {\em canonical} $V$--line bundle $K_\Sigma$   has  singularity
data
$ (\alpha_i-1)/\alpha_i$ for $1\leq i\leq \nu$, and
$\deg|K_\Sigma|=-2$, hence  rational degree
\[
\kappa:= \deg^V
K_\Sigma=-2+\sum_{i=1}^\nu\Bigl(1-\frac{1}{\alpha_i}\Bigr).
\]
The Seifert manifold  $M$  with non-normalized Seifert invariants
$((\alpha_i,\beta_i)_{i=1}^\nu)$ (or, equivalently, with
normalized Seifert invariants $(b;(\alpha_i,\omega_i)_{i=1}^\nu)$,
cf \ref{ss: s2}), is the unit circle bundle of the $V$--line bundle
${\bL}_0$ with rational degree $\ell=e$, and singularity data
\[
\frac{\omega_i}{\alpha_i},\;\;1\leq i\leq \nu, \ \
(\omega_i\equiv -\beta_i\ (\mmod \ \alpha_i)).
\]
Denote by $\lan {\bL}_0\ran \subset \pic$ the cyclic group
generated by ${\bL}_0$. Then one has the following exact sequence:
\[
0\to \lan {\bL}_0\ran\to \pic\stackrel{\pi^*}{\to } \picm\to 0,\]
where $\pi^*$ is  the  pullback map induced by the natural
projection $\pi\co M\to\Sigma$. Therefore, the above exact sequence
identifies for every $L\in \pic$ the pullback $\pi^*(L)$ with the
class $[L]\in  \pic/\lan {\bL}_0\ran$.

For every $L\in \pic$, $c:=\deg^VL$ we set
\[
\rho(L):=\frac{\deg^VK_\Sigma-2c}{2\ell}= \frac{\kappa}{2\ell}-
\frac{c}{\ell}\in {\bQ}.
\]
For every class $u\in \picm$ we can find an unique $E_u\in \pic$
such that $u=[E_u]\ \ \mbox{and } \ \ \rho(E_u)\in[0,1)$.
We say that $E_u$ is the \emph{canonical representative} of $u$.
As explained in \cite{MOY,Nico2} there is a natural bijection
\[
\picm\ni u\mapsto \si(u)\in Spin^c(M)
\]
with the property that $\det \si(u)= 2u-[K_\Sigma]\in \picm$.
The canonical $spin^c$ structure $\si_{can}\in Spin^c(M)$
 corresponds to $u=0$.
In fact, $\picm$ can be identified in a natural way to $H$ via the
Chern class. Then $\si(u)$, in terms of the $H$--action described in
\ref{ss: 2.4b},
is given by $u\cdot \si_{can}$. In this case, if one writes
$\rho_0:=\rho(E_0)$ one has
\[
\rho_0=\Bigl\{\frac{\kappa}{2\ell}\Bigr\}\ \ \mbox{and}\ \ E_0:=
n_0{\bL}_0, \ \ \mbox{with}\ \
n_0:=\Bigl\lfloor \frac{\kappa}{2\ell}\Bigr\rfloor.
\]
We denote the orbifold invariants of $E_0$ by
$\frac{\gamma_i}{\alpha_i}$. Observe that
\[
\frac{\gamma_i}{\alpha_i}=\Bigl\{
\frac{n_0\omega_i}{\alpha_i}\Bigr\}.
\]
The Seifert manifold $M$  admits a natural metric, the so called
\emph{Thurston  metric} which we denote by $g_0$. The
$(\si_{can},g_0,0)$--monopoles were   explicitly described  in
\cite{MOY, Nico2}.

The space $\modu_0^*$ of irreducible $(\si_{can},g_0,0)$
monopoles  on $M$ consists  of several  components  parametrized
by  a subset of
\[
S_0=\bigl\{ E=E_0+n{\bL}_0\in \pi;\;\; 0< |\nu(E)|\leq \frac{1}{2}
\deg^VK_\Sigma\bigr\},
\]
where
\[
 \nu(E):=\deg^V (E_0+n {\bL}_0) -\frac{1}{2}\deg^V K_\Sigma.
\]
More precisely, consider the sets
\[
S_0^+:=\Bigl\{ E\in S_0;\;\; \nu(E)<0,\;\;\deg|E|\geq
0\Bigr\},
\]
\[
S_0^-= \Bigl\{ E\in S_0;\;\;\nu(E)>0,\;\;\deg|K_\Sigma-E|\geq 0\}.
\]
To every $E\in S_0^+$ there corresponds a component $\modu^+_E$ of
$\modu^*_0$ of dimension $2\deg|E|$,  and to every $E\in S_0^-$
there corresponds a component $\modu_E^-$  of $\modu^*_0$ of
dimension $2\deg|K_\Sigma-E|$.

The Kreck--Stolz invariant $KS(\si_{can},g_0,0)$ is given  by (see \cite{Nico2})
\[
KS_M(\si_{can},g_0,0)=\ell+1-4\ell\rho_0(1-\rho_0)
\]
\[
+4\nu\rho_0 -4\sum_{i=1}^\nu\bms(\omega_i,\alpha_i)-8\sum_{i=1}^\nu\bms(\omega_i,\alpha_i; \frac{\gamma_i+\rho_0\omega_i}{\alpha_i}, -\rho_0)
\]
\[
+4\left\{
\begin{array}{ccc}
-\sum_{i=1}^n\LL \frac{r_i\gamma_i}{\alpha_i}\RR & {\rm if} & \rho_0=0\\
& &\\
\frac{2+\kappa}{2}(1-2\rho_0)-\sum_{i=1}^\nu
\Bigl\{\frac{r_i\gamma_i+\rho_0} {\alpha_i}\Bigr\} & {\rm if} &
\rho_0\neq 0,
\end{array}
\right.
\]
where
\[
r_i\omega_i\equiv 1\ (\mmod \ \alpha_i), \ i=1,\ldots,\nu.
\]
Above, $\bms(h,k;x,y)$ is the Dedekind--Rademacher sum defined in
Appendix B, where we list some of its basic   properties as
well.

The following result is a consequence of the analysis carried out in \cite{Nico0,Nico2}.

\gtsubsection {Proposition}
\label{ss: sws2} {\sl
{\rm(a)}\qua If $\rho_0\neq 0$ and $\modu^*_0$ has only zero dimensional
components then $(g_0,0)$ is a good parameter and
\[
\ssw^0_M(\si_{can})=\frac{1}{8}KS_M(\si_{can},g_0,0)+|S_0^+|+|S_0^-|.
\]
{\rm(b)}\qua If $g_0$ has positive scalar curvature then $(g_0,0)$ is a
good parameter, $S_0^+=S_0^-=\emptyset$ and
\[
\ssw^0_M(\si_{can})=\frac{1}{8}KS(\si_{can},g_0,0).
\]
} Notice that part (b) can be applied for the links of quotient
singularities.

One of the main obstructions is, that in many cases, the above
theorem cannot be applied (ie, the natural parameter provided by
the natural  Seifert metric is not ``good'', cf \ref{ss: 3.2}).

Fortunately, the Seiberg--Witten invariant has an alternate
combinatorial description as well. To formulate it we need to
review a few basic topological  facts.

\gtsubsection{The Reidemeister--Turaev torsion}
\label{ss: 3.4}
\setcounter{equation}{0}
According to Turaev \cite{Tu5}   a  choice of a $spin^c$ structure on $M$ is
equivalent to a choice of an Euler structure. For  every $spin^c$
structure $\si$ on $M$, we denote by
$$\tT_{M,\si}=\sum_{h\in H}\tT_{M,\si}(h)\, h\in {\bQ}[H],$$
 the sign refined
{\em Reidemeister--Turaev torsion}  determined by the Euler
structure associated to $\si$. (For its detailed description, see
\cite{Tu5}.)
 Again,  it is convenient to think of $\tT_{M,\si}$ as a function
$H\ra {\bQ}$ given by $h\mapsto \tT_{M,\si}(h)$.  The Poincar\'{e}
duality implies that $\tT_{M,\si}$ satisfies the symmetry condition
\begin{equation*}
\tT_{M,\si}(h)=\tT_{M,\bar{\si}}(h^{-1})\ \ \ \mbox{for all} \ h\in H.
\tag{$*$}
\label{3.4.*}
\end{equation*}
Recall that the  {\em augmentation map} $\aug\co  {\bQ}[H]\ra {\bQ}$ is defined by
\[
\sum a_h\, h\mapsto \sum a_h.
\]
 It is known that $\aug(\tT_{M, \si})=0$.

\gtsubsection{The Casson--Walker invariant and the modified
Reidemeister--Turaev torsion}\label{ss: 3.5}
\setcounter{equation}{0}
Denote by $\lambda(M)$ the Casson--Walker invariant of $M$
normalized as in Lescop's book (cf \cite[Section 4.7]{Lescop}),  and denote by
$H\ni h\mapsto \tT_{M,\si}^0(h)$ the
{\em modified Reidemeister--Turaev torsion}
\[
\tT^0_{M,\si}(h):=\tT_{M,\si}(h)- \lambda(M)/|H|.
\]
We have the following result.
\gtsubsection{Theorem} {\rm\cite{Nico5}}\qua
\label{ss: 3.6}
 {\sl $\bSW^0_{M,\si}(h)=\tT^0_{M,\si}(h)$ for all $\si\in Spin^c(M)$
and $h\in H$. }
\gtsubsection{The Fourier transform}
\label{ss: 3.7}
\setcounter{equation}{0} Later we will need a dual
description of these invariants in terms of Fourier transform.
Denote by $\hat{H}$  the Pontryagin dual of $H$, namely $\hat{H}:=
\Hom(H, U(1))$.
 The Fourier transform of any function $f\co  H \ra {\bC}$ is  the function
 \[
\hat{f}\co \hat{H}\ra {\bC},\quad \hat{f}(\chi)=\sum_{h\in H}
f(h)\bar{\chi}(h).
 \]
 The function $f$ can be recovered from its Fourier transform via the \emph{ Fourier inversion formula}
 \[
 f(h)=\frac{1}{|H|}\sum_{\chi\in \hat{H}}\hat{f}(\chi)\chi(h).
 \]
 Notice  that $\aug(f)=\hat{f}(1)$, in particular $\hat{\tT}_{M,\si}(1)=\aug(\tT_{M, \si})=0$.
By the above identity,
\begin{equation*}
\begin{aligned}[b]
 \ssw^0_M(\si)=\bSW^0_{M,\si}(1)=-\frac{1}{|H|}\lambda(M)+\frac{1}{|H|}
\sum_{\chi\in \hat{H}}\hat{\tT}_{M,\si}(\chi)\\
=-\frac{1}{|H|}\lambda(M)+\tT_{M,\si}(1).
\end{aligned}
\tag{1}
\label{3.7.2}
\end{equation*}
The symmetry condition \ref{ss: 3.4}(\ref{3.4.*}) transforms into
\begin{equation*}
\hat{\tT}_{M,\si}(\chi)=
\hat{\tT}_{M,\bar{\si}}(\bar{\chi}).
\tag{2}
\label{3.7.1}
\end{equation*}
It is convenient to use the notation
$\sum_{\chi}'$ for a summation where $\chi$ runs over all the {\em
non-trivial} characters of  $\hat{H}$.

\gtsubsection{The identification $Spin^c(M)\to Q^c(M)$ via the Seiberg--Witten invariant}
\label{ss: swiden} Sometimes it is important to have an efficient way
to recover the $spin^c$ structure $\si$ (or, equivalently, the
quadratic function $\gq^c(\si)$, cf \ref{ss: 2.4c}) from
the Seiberg--Witten invariant $\bSW^0_{M,\si}$, or from
$\tT_{M,\si}$. In order to do this,
we first recall that Turaev in
\cite[Theorem 4.3.1]{Tu5} proves the following identity for {\em any} $\si$
and $g,h\in H$.
\[
\tT_{M,\si}(1)-\tT_{M,\si}(h)-\tT_{M,\si}(g)+\tT_{M,\si}(gh)=
-b_M(g,h)\ (\mmod\ \Z).
\]
Clearly there is a similar identity for
$\tT^0_{M,\si}$ instead of  $\tT_{M,\si}$. By Fourier inversion, this reads
\begin{equation*}
\frac{1}{|H|}{\sum_{\chi}}'\hat{\tT}_{M,\si}(\chi)(\chi(h)-1)(\chi(g)-1)=
-b_M(g,h)\ (\mmod\ \Z).
\tag{1}
\end{equation*}
This identity has a ``refinement'' in the following sense (see
\cite[3.3]{N3}): for any $spin^c$ structure $\si$, the
map   $H\ni h\mapsto q_{\si}'(h)$, defined by,
\begin{equation*}
\begin{aligned}
q_{\si}'(h):=\tT_{M,\si}(1)-\tT_{M,\si}(h)=\tT^0_{M,\si}(1)-\tT^0_{M,\si}(h)\\
=\bSW^0_{M,\si}(1)-\bSW^0_{M,\si}(h)\;\;({\rm mod}\; {\bZ}),
\end{aligned}
\end{equation*}
is a quadratic function associated with $b_M$. Moreover, the
correspondence $ {\gq}_{sw}^c :Spin^c(M)\to  Q^c(M)$ given by
$\si\mapsto q_{\si}'$ is a bijection.

\gtsubsection{Proposition}
\label{ss: tur} {\sl $ {\gq}_{sw}^c={\gq}^c$, ie, the
above bijection is exactly the canonical identification $\gq^c$
considered in \ref{ss: 2.4c}.}

\begin{proof} Let us denote $\gq^c(\si)$ by $q_{\si}$.
 Since both  maps ${\gq}^c$ and ${\gq}_{sw}^c$ are $H$--equivariant,
 it suffices to show that $q_\eps=q'_\eps$ for some  $spin$ structure
$\eps$ on $M$.  Fix  $\eps \in Spin(M)$, and as in
\ref{ss: 2.4}, pick  a simple connected oriented 4--manifold $N$ with
$\partial N=M$ together with a  $\tilde{\eps}\in Spin(N)$
such that the restriction of $\tilde{\eps}$ to $M$ is $\eps$. In
particular, $q_\eps=-q_{D(L),0}$, cf \ref{ss: 2.4}.
 For every  $h\in H$  set $\si_h:=
h\cdot \si(\eps)$. Pick $\tilde{h}\in H^2(N,\Z)=L'$ such that
$[\tilde{h}] =h$. For $h=0$ we choose $\tilde{h}=0$. Set
$\tilde{\si}_h:=\tilde{h}\cdot \si(\tilde{\eps})\in Spin^c(N)$.
Then  $c_1(\det {\tilde{\si}}_h)=2\tilde{h}$.   Observe that
\[
\bSW^0_{M,\eps}(h)=\frac{1}{8} KS(\si_h)\ (\mmod \ \Z).
\]
But the Atiyah--Patodi--Singer index theorem implies
(see eg \cite[page 197]{Lim}):
\[
\frac{1}{8}
KS(\si_h)=\frac{1}{8}(c_1(\det{\tilde{\si}}_h),c_1(\det{\tilde{\si}}_h))_\Q-\frac{1}{8}\mbox{signature}(N) \ (\mmod\ \Z),
\]
for any $h\in H$. Thus $q'_\eps(h)=-\frac{1}{2}\, (\tilde{h},\tilde{h})_\Q=q_\eps(h) \
(\mmod\ \Z)$. \end{proof}

Via the Fourier transform, the above identity is equivalent to the following
one, valid for any $h\in H$:
\begin{equation*}
\frac{1}{|H|}{\sum_{\chi}}'\hat{\tT}_{M,\si}(\chi)(\chi(h)-1)=
-\gq^c(\si)(h)\ (\mmod\ \Z).
\tag{2}
\end{equation*}

\gtsubsection{Remark}\label{ss: blij} The above discussion can be compared
 with the  following identity. Let us keep the
notations of \ref{ss: 2.3}. Let $\sigma(L)$
denote the signature of $L$, and $k\in L'$ a characteristic element.
Then the ($\mmod\ 8$)--residue class
of $\sigma(L)-(k,k)_\Q\in \Q/8\Z$ depends only on the  quadratic
function  $q=q_{D(L),k}$.
 In fact one has the
following formula of van der Blij \cite{Blij} for the Gauss sum:
$$\gamma(q):=|H|^{-1/2}\, \sum_{x\in H}\,
e^{2\pi i q(x)}=e^{\frac{\pi i}{4}(\sigma(L)-(k,k)_\Q)}.$$ If
$k\in \im(i_L)$ then $\sigma(L)-(k,k)_\Q=\sigma(L)-(k,k)\in
\Z/8\Z$.

\section{Analytic invariants and the main conjecture}

\gtsubsection{Definitions}
\label{ss: 4.1} Let $(X,0)$ be a normal surface
singularity. Consider the  holomorphic line bundle
$\Omega^2_{X\setminus\{0\}}$ of holomorphic 2--forms on $X\setminus\{0\}$.
If this line bundle is holomorphically trivial then we say
that $(X,0)$ is \emph{Gorenstein}.
If some power of this line bundle is holomorphically trivial then we say
that $(X,0)$ is $\Q$--\emph{Gorenstein}.
If $\Omega^2_{X\setminus\{0\}}$ is {\em topologically}  trivial we say
that $(X,0)$ is \emph{numerically Gorenstein}. The first two conditions
are  analytic, the third depends only on the link $M$ (cf \ref{ss: 2.5}).

\gtsubsection{The geometric genus}
\label{ss: 4.2}
Fix a resolution $\pi\co \tilde{X}\to
X$ over a sufficiently small Stein representative  $X$ of the germ
$(X,0)$. Then $p_g:=\dim \, H^1(\tilde{X}, {\mathcal
O}_{\tilde{X}})$ is finite and independent of the choice of $\pi$.
It is called the {\em geometric genus } of $(X,0)$. If
$p_g(X,0)=0$ then the singularity $(X,0)$ is called \emph{rational}.

\gtsubsection{Smoothing invariants}
\label{ss: 4.3}
Let $(X,0)$ be as above. By a
\emph{smoothing} of $(X,0)$ we mean a proper flat analytic germ
$f\co ({\mathcal X},0)\to (\C,0)$ with an isomorphism
$(f^{-1}(0),0)\to (X,0)$. Moreover, we assume that $0$ is an
isolated singular point of the germ $({\mathcal X},0)$.

If ${\mathcal X}$ is a sufficiently small contractible Stein
representative of $({\mathcal X},0)$, then for sufficiently small
$\eta$ ($0<|\eta|\ll 1$)  the fiber $F:= f^{-1}(\eta)\cap
{\mathcal X}$ is smooth, and its diffeomorphism  type is
independent of the choices. It is a connected
 oriented real 4--manifold with boundary
$\partial F$ which can be identified with the link $M$ of $(X,0)$.

We will use the following notations: $\mu(F)=\rank H_2(F,\Z)$
(called the \emph{Milnor number}); $\bif_F$= the intersection form  of
$F$ on $H_2(F,\Z)$; $(\mu_0,\mu_+,\mu_-)$ the Sylvester invariant
of $\bif_F$; $\si(F):=\mu_+-\mu_- $ the signature of $F$. Notice
that the Milnor fiber $F$, hence its invariants too,  in general
depend on the choice of the (irreducible component) of the
smoothing.

If $M=\partial F$ is a rational homology sphere then  $\mu_0=0$, hence
$\mu(F)=\mu_++\mu_-$. It is known that for a smoothing of a
Gorenstein singularity $\rank H_1(F,\Z)=0$ \cite{GSt}. Therefore,
in this case $\mu(F)+1$ is the topological Euler characteristic
$\chi_{top}(F)$ of $F$.

The following relations connect the invariants $p_g,\ \mu(F)$ and
$\sigma(F)$. The next statement is formulated for rational
homology sphere links, for the general statements the reader can
consult the original sources \cite{Durfee,Laufermu,Steenbrink}
(cf also with \cite{LW}).

\gtsubsection{Theorem}
\label{ss: 4.4} {\sl  Assume that the link $M$ is a rational
homology sphere. Then the following identities hold.

{\rm(1)\qua [Wahl, Durfee, Steenbrink]}\qua $4p_g=\mu(F)+\si(F)$.

In addition, if $(X,0)$ is Gorenstein, then

{\rm(2)\qua [Laufer, Steenbrink]}\qua $\mu(F)=12 p_g+K^2+\#{\mathcal V}$,
where $K^2+\#{\mathcal V}$ is the topological invariant of $M$
introduced in \ref{ss: 2.11}.

In particular, for Gorenstein singularities, (1) and (2)
give} $\si(F)+8p_g+K^2+\#{\mathcal V}=0$.

This shows that modulo the link-invariant
$K^2+\#{\mathcal V}$ there are two (independent) relations
connecting $p_g,\ \mu(F)$ and $\si(F)$, provided that $(X,0)$ is
Gorenstein.  So, if by some other argument one can recover one of
them from the topology of $M$, then all of them can be determined
from $M$.

In general, these invariants cannot be computed from $M$. Here one
has to emphasize two facts. First, if $M$ is not a rational
homology sphere, then one can construct easily (even hypersurface)
singularities with the same link but different $(\mu,\sigma,p_g)$.
On the other hand, even if we restrict ourselves to rational
homology links, if we consider all the possible analytic
structures of $(X,0)$, then again $p_g$ can vary. For example, in
the case of ``weakly'' elliptic singularities, there is a
topological upper bound of $p_g$ (namely, the length of the
elliptic sequence, found by Laufer and S\thinspace S-T Yau) which equals
$p_g$ for Gorenstein singularities; but $p_g$ drops to 1  for a
generic analytic structure (fact proved by Laufer). For more
details and examples, see the series of articles of S\thinspace S-T Yau
(eg \cite{Yau}), or   \cite{Neminv}.
On the other hand, the first author in \cite{Neminv} conjectured
that for Gorenstein singularities with rational homology sphere
links the invariants $(\mu,\sigma,p_g)$ can be determined from the
topology of $(X,0)$ (ie, from the link $M$); (cf also with the
list of conjectures in \cite{NW2}). The conjecture  is true for
rational singularities \cite{Artin,Artin2},  minimally
elliptic singularities \cite{Lauferme}, ``weakly'' elliptic
singularities \cite{Neminv}, and some special hypersurface singularities
\cite{FS,NW}, and special
complete intersections \cite{NW};  in all cases
 with explicit formulae for $p_g$. But in general,  even a conjectural
topological candidate (computed from $M$) for $p_g$ was
completely open.
The next conjecture provides exactly this topological candidate
(which is also a ``good'' topological upper bound, cf introduction).

\gtsubsection{The Main Conjecture}
\label{ss: 4.5} {\sl Assume that $(X,0)$ is a normal
surface singularity whose link $M$ is a rational homology sphere.
Let $\si_{can}$ be the canonical $spin^c$  structure on $M$. Then,
conjecturally, the following facts hold:

\vspace{2mm}

{\rm(1)}\qua For any $(X,0)$, there is a topological upper bound
for $p_g$ given by:
$$\ssw^0_M(\si_{can})-\frac{K^2+\#{\mathcal V}}{8}\geq p_g.$$
{\rm(2)}\qua If $(X,0)$ is $\Q$--Gorenstein, then in (1) one has equality.

{\rm(3)}\qua In particular, if $({\mathcal X},0)$ is a smoothing of a
Gorenstein singularity $(X,0)$ with Milnor fiber $F$, then
$$-\ssw^0_M(\eps_{can})=\frac{\si(F)}{8}.$$}
If $(X,0)$ is numerically Gorenstein and $M$ is a $\Z_2$--homology
sphere then $\eps_{can}=\si_{can}$ is the unique spin structure of
$M$; if $M$ is an integral  homology sphere then in the above
formulae $-\ssw^0_M(\eps_{can})=\lambda(M)$, the Casson invariant
of $M$.

\gtsubsection{Remarks}
\label{ss: 4.6} 

(1)\qua Assume that $(X,0)$ is a hypersurface
Brieskorn singularity whose link is an {\em integral} homology
sphere. Then $\lambda(M)=\si(F)/8$ by
 a result of Fintushel and Stern \cite{FS}. This fact was generalized  for
Brieskorn--Hamm  complete intersections  and for suspension hypersurface
singularities ($(X,0)=\{g(x,y)+z^n=0\}$) with $H_1(M,\Z)=0$ by
Neumann and Wahl \cite{NW}. In fact, for a normal complete
intersection surface singularity with $H=H_1(M,\Z)=0$,   Neumann
and Wahl conjectured $\lambda(M)=\si(F)/8$. This conjecture was
one of the starting points of our investigation.

The result of Neumann--Wahl \cite{NW} was re-proved and reinterpreted
by Collin and Saveliev (see \cite{CoS} and \cite{Co}) using equivariant
Casson invariant and cyclic covering techniques.

(2)\qua The family of $\Q$--Gorenstein singularities is rather large:
it contains eg the rational singularities \cite{Artin,Artin2},
the singularities with good $\C^*$--actions and with rational homology sphere
links \cite{Neu}, the minimally elliptic singularities \cite{Lauferme},
and all the isolated complete intersection  singularities. Neumann--Wahl have conjectured in \cite{NW2} that all the singularities in \ref{ss: 4.5} (2) are finite abelian quotients of complete intersection singularities.

(3)\qua If one wants to test the Conjecture for rational or elliptic
singularities (or in any example where $p_g$ is known), one should
compute the  corresponding Seiberg--Witten invariant. But, in some
cases, even if   all the terms in the main conjecture can (in
principle) be  computed, the identification of these contributions
in the main formula can create difficulties (eg involving
complicated identities  of Dedekind sums and lattice point
counts).

\gtsubsection{Remark}
\label{ss: conjadd}Notice, that in the above Conjecture, we have 
automatically built in the following statements as well.

\smallskip

(1)\qua For any normal singularity $(X,0)$ the topological
invariant
\[
\ssw^0_M(\si_{can})-\frac{K^2+\#{\mathcal V}}{8}
\]
is non-negative. Moreover, this topological invariant is zero if
and only if $(X,0)$ is rational. This provides a new topological
characterization of the rational singularities.

\smallskip

(2)\qua Assume that $(X,0)$ (equivalently, the link) is numerically
Gorenstein. Then the above topological invariant is 1 if and only
if $(X,0)$ is minimally elliptic (in the sense of Laufer). Again,
this is a new topological characterization of minimally  elliptic
singularities.

\gtsubsection{Remark}\label{ss: 48} The invariant $K^2+\#{\mathcal V}$ appears not only in the type of
results listed in \ref{ss: 4.4}, but also in  other topological contexts. For example, it can be identified with the Gompf  invariant $\theta(\xi)$ defined in \cite[4.2]{Gompf} (see also \cite[11.3.3]{GS}).
This appears as an ``index defect'' (similarly to the signature defect of Hirzebruch) (cf also with \cite{Durfee} and \cite{LW}).

 More precisely, the almost complex structure on $TM\oplus \R_M$ (cf \ref{ss: 2.6})
determines a {\em contact structure} $\xi_{can}$ on $M$ (see eg
\cite[page 420]{GS}), with $c_1(\xi_{can})$ torsion element. Then
the Gompf invariant $\theta(\xi_{can})$, computed via $\tilde{X}$, is
 $K^2-2\chi_{top}(\tilde{X})-3\si(\tilde{X})=K^2+\#{\mathcal V}-2$.

In fact, in our situation, by \ref{ss: 2.11},  $\theta(\xi_{can})$ can be recovered
from the oriented $C^{\infty}$ type of $M$ completely. In the Gorenstein case, in the presence of a smoothing, $\theta(\xi_{can})$
computed from the Milnor fiber $F$,
 equals $-2-2\mu(F)-3\si(F)$. The identity $K^2+\#{\mathcal V}+2\mu(F)+3\si(F)=0$ can be deduced
from \ref{ss: 4.4} as well.

The goal of the remaining part of the present paper is to describe the needed
topological invariants in terms of the plumbing graphs of the link, and
finally,  to provide
a list of examples supporting the Main Conjecture.

\section{Invariants computed from the plumbing graph}

\gtsubsection{Notation}
\label{ss: 5.1} The goal of this section is to list some
formulae for the invariants $K^2+\#{\mathcal V}$, $\lambda(M)$ and
$\tT_{M,\sigma}$ from the resolution graph of $M$ (or,
equivalently, from any negative definite plumbing). The formulae
are made explicit for star-shaped graphs in terms of their
Seifert invariants. For notations, see  \ref{ss: no}, \ref{ss: s1} and
\ref{ss: s2}.

Let $I^{-1}$ be the inverse of the intersection matrix $I$. For
any $v,w\in {\mathcal V}$, $I^{-1}_{vw}$   denotes the
$(v,w)$--entry of $I^{-1}$.    Since $I$
is negative definite, and the graph is connected,  $I_{vw}^{-1}<0$
for each entry $v,w$. Since $I$ is described  by a tree, these
entries have the following interpretation as well. For any two
vertex $v,w\in {\mathcal V}$, let $p_{vw}$ be the unique minimal
path in the graph connecting $v$ and $w$, and let $I_{(vw)}$ be
the matrix obtained from $I$ by deleting all the lines and columns
corresponding to the vertices on the path $p_{vw}$ (ie,
$I_{(vw)}$ is the intersection matrix of the complement graph of
the path). Then $I_{vw}^{-1}=-|\det(I_{(vw)})/\det(I)|$.

For simplicity we will write $E_v,\, D_v$, \ldots instead of
$[E_v],\, [D_v] \ldots$.

\gtsubsection{The invariant $K^2+\#{\mathcal V}$ from the plumbing
graph}
\label{ss: 5.2}Let $Z_K=\sum_v\, r_v E_v$. Since
$Z_K=(i_L\otimes \Q)^{-1}(\sum (e_v+2)D_v)$ (cf \ref{ss: no}),
one has clearly  $r_v=\sum_w(e_w+2)I^{-1}_{vw}$.

Then a ``naive formula'' of $K^2=Z_K^2$ is
$$Z_K^2=\sum_v\, r_v Z_K\cdot E_v= \sum_v\, r_v(e_v+2)=
\sum_{v,w} (e_v+2)(e_w+2)I_{vw}^{-1}.$$
However, we prefer a different
form for $r_v$ and $Z_K^2$ which involves only a small part of the
entries of $I^{-1}$. Indeed, let us consider the class
\[
D=\sum_vi_L(E_v)+\sum_v(2-\delta_v)D_v\in L'.
\]
Then clearly $D\cdot E_v=e_v+\delta_v+2-\delta_v=e_v+2$, hence $D=(i_L\otimes\Q)Z_K$. Therefore,
\[
Z_K=\sum_vE_v+\sum_v (2-\delta_v)(i_L\otimes\Q)^{-1} (D_v),
\]
 hence
\[
r_v=1+\sum_w(2-\delta_w)I^{-1}_{vw}.
\]
Moreover,

$Z_K^2=Z_K\cdot D=$\vspace{-2mm}
$$\sum_vZ_K\cdot E_v+\sum_v(2-\delta_v)Z_K\cdot D_v= \sum_v(e_v+2)+\sum_v(2-\delta_v)r_v,$$
hence by the second formula for $r_v$ we deduce
\[
K^2+\#{\mathcal V}=\sum_ve_v+3\#{\mathcal V}+2 +\sum_{v,w}(2-\delta_v)
(2-\delta_w)I^{-1}_{vw}.
\]
 In particular, this number depends only
of those entries of $I^{-1}$
 whose index set
runs over the rupture points ($\delta_v\geq 3$) and the
end-vertices ($\delta_v=1$) of the graph.

For cyclic quotient singularities,  the above formula for $K^2$
goes back to the work of Hirzebruch.  In fact,
the right hand side can also be expressed
 in terms of Dedekind sums,  see eg \cite[5.7]{LW}
 and \cite{Ishii} (or \ref{ss: 4.8} here).

\gtsubsection{The Casson--Walker invariant from plumbing}
\label{ss: 5.3}
We recall a formula for the Casson--Walker invariant for  plumbing 3--manifolds
proved  by A Ratiu in his  dissertation \cite{Ratiu}. In fact, the
formula can also be recovered from the surgery formulae of Lescop
\cite{Lescop} (since any plumbing graph can be transformed into a
precise surgery data, see eg A1). The first author thanks Christine  Lescop
for providing him  all the details and information about it. We have
\[
-\frac{24}{|H|}\,\lambda(M)=\sum_ve_v+3\#{\mathcal V}+\sum_v(2-\delta_v)\,
I_{vv}^{-1}.
\]
If $M=L(p,q)$, then this can be transformed into $\lambda(L(p,q))=
p\cdot \bms(q,p)/2$. (Here we emphasize that by our notations,
$L(p,q)$ is obtained by $-p/q$--surgery  on the unknot in $S^3$, as
in \cite[page 158]{GS}, and {\em not} by $p/q$--surgery as in
\cite[page 108]{W}.)

\gtsubsection{The Casson--Walker invariant for Seifert manifolds}
\label{ss: s3}
Assume that $M$ is a Seifert manifold as in \ref{ss: s1} and \ref{ss: s2}.
Using \cite{Lescop}, Proposition 6.1.1, one has the following
expression:
\[
-\frac{24}{|H|}\lambda(M)=
\frac{1}{e}\Big(2-\nu +\sum_{i=1}^\nu\frac{1}{\alpha_i^2}\Big)+
e+3+12\sum_{i=1}^\nu \,  \bms(\beta_i,\alpha_i).
\]
(Warning: our notations for the Seifert invariants  differ  slightly  from
those used in \cite{Lescop}; and also, our $e$ and $b$ have opposite signs.)

\gtsubsection{$K^2+\#{\mathcal V}$ for Seifert manifolds}
\label{ss: s4}
Using \ref{ss: 5.3}  we deduce
\[
-\frac{24}{|H|}\lambda(M)=\sum_ve_v+3\#{\mathcal V}+\sum_{i=1}^\nu
I_{\bar{v}_i\bar{v}_i}^{-1}+
(2-\nu)I_{\bar{v}_0\bar{v}_0}^{-1}.
\]
For Seifert manifolds,  \ref{ss: 5.2} can be rewritten as

$
K^2+\#{\mathcal V}=$\vspace{-2mm}
$$
2+ \sum_ve_v+3\#{\mathcal V}+\sum_{i=1}^\nu
I_{\bar{v}_i\bar{v}_i}^{-1}+ (2-\nu)^2I_{\bar{v}_0\bar{v}_0}^{-1}+
2\sum_{i=1}^\nu (2-\nu)I_{\bar{v}_0\bar{v}_i}^{-1}+
\sum_{i\not=j}I_{\bar{v}_i\bar{v}_j}^{-1}.
$$
 Using the interpretation of the entries of $I^{-1}$ given in \ref{ss: 5.1}, one gets easily that
\begin{equation*}
I_{\bar{v}_0\bar{v}_0}^{-1}=\frac{1}{e};\ \
I_{\bar{v}_i\bar{v}_0}^{-1}=\frac{1}{e\alpha_i};\ \
I_{\bar{v}_i\bar{v}_j}^{-1}=\frac{1}{e\alpha_i\alpha_j} \
(i\not=j;\ 1\leq i,j\leq \nu).
\tag{1}
\label{s41}
\end{equation*}
Therefore, these identities and \ref{ss: s3} give:
\[
K^2+\#{\mathcal V}=
\frac{1}{e}\Big(2-\nu +\sum_{i=1}^\nu \frac{1}{\alpha_i}\Big)^2 +e +
5+12\sum_{i=1}^\nu \, \bms(\beta_i,\alpha_i).
\]
It is instructive to compare this expression with \ref{ss: s3} and also
with the
 coefficient $r_0$ of $Z_K$, namely with
$
r_0=1+ (2-\nu +\sum_i1/\alpha_i)/e.
$
\gtsubsection{The Reidemeister--Turaev torsion}
\label{ss: rt1}
In the remaining part of this section we provide a formula for the
torsion $\tT_M$ of $M$ using the plumbing representation of $M$.
We decided not to distract the reader's attention
from the main message of the paper  and we  deferred its proof to Appendix A.

\gtsubsection{Theorem}
\label{ss: th1} {\sl Let $M$ be an oriented  rational homology 3--manifold
represented by a negative definite plumbing graph $\Gamma$. (Eg,
let $M$ be the link of a normal surface
singularity $(X,0)$, and $\Gamma=\Gamma(\pi)$ be
one of its  resolution graphs.) In the sequel, we keep the
notations used above (cf \ref{ss: no} and  \ref{ss: 5.1}). For any
$spin^c$ structure $\si \in Spin^c(M)$, consider the unique
element $h_\si \in H$ such that $h_\si\cdot \si_{can}=\si$. Then
for any $\chi\in \hat{H}$, $\chi\not=1$, the following  identity
holds:
\[
\hat{\tT}_{M,\si}(\bar{\chi})
=\bar{\chi}(h_{\si})\cdot \prod_{v\in {\mathcal V}}\
\Bigl(\chi(g_v)-1\Bigl)^{\delta_v-2}.
\]
 The right hand side should
be understood as follows. If $\chi(g_v)\not=1$ for all $v$ (with
$\delta_v\not=2$), then
the expression is well-defined. Otherwise, the right hand
 side is computed via a limit (regularization procedure). More
precisely, fix a vertex  $v_*\in {\mathcal V}$ so that
$\chi(g_{v_*})\not=1$, and let $\vec{b}^*$ be the column
vector with entries =1 on the place $v_*$ and zero otherwise. Then
find  $\vec{w}^*\in \Z^n$ with entries $\{w_v^*\}_v$
 in such a way that $I\vec{w}^*=-m^*\cdot \vec{b}^*$  for some integer
$m^*>0$. Then
\[
\hat{\tT}_{M,\si}(\bar{\chi})=\bar{\chi}(h_{\si})\cdot \lim_{t\to 1}\,
\prod_{v\in {\mathcal V}}\ \Bigl(t^{w_v^*}\,
\chi(g_v)-1\Bigl)^{\delta_v-2}.
\] }
\noindent
The above limit always exists.
Moreover, once  $v_*$ is fixed, the vector $\vec{w}^*$  is unique
modulo a positive multiplicative factor (which does nor alter the limit).
In fact, the above limit is independent even of the choice of
$v_*$ (as long as $\chi(g_{v_*})\not=1$). This follows also from the general
theory (cf also \ref{ss: rt4} and \ref{ss: limit}),
but it also has  an elementary
combinatorial proof given in Lemma \ref{le}.
(This can be read independently
from the other parts of the proof.)
In fact, by Lemma \ref{le},  the set of vertices $v_*$,
providing   a suitable $\vec{w}^*$ in the limit expression,  is even
larger than the set identified in the theorem:
 one can take any $v_*$ which satisfy
{\em either} $\chi(g_{v_*})\not=1$ {\em or}  it has an adjacent vertex $u$
with $\chi(g_{u})\not=1$.

\gtsubsection{The torsion of Seifert manifolds}
\label{ss: sm}
In this paragraph we use the notations of \ref{ss:  s1} and \ref{ss: s2}.
Recall that we
introduced $\nu+1$ distinguished vertices $\bar{v}_i, \
0\leq i\leq \nu$, whose degree is $\not=2$
(and $g_{\bar{v_0}}$ is the central vertex).
Fix $\chi$ and first assume that $\chi(g_{\bar{v}_i})\not=1$ for
all $0\leq i\leq \nu$. Then, the above theorem reads as
\[
\hat{\tT}_{M,\si_{can}}(\bar{\chi})=
\frac{\bigl(\chi(g_{\bar{v}_0})-1\bigr)^{\nu-2}}
{\prod_{i=1}^\nu\ \bigl(\chi(g_{\bar{v}_i})-1\bigr)}.
\]
 If there is one index $1\leq i\leq \nu$ with $\chi(g_{\bar{v}_i})=1$, then
necessarily $\chi(g_{\bar{v}_0})=1$ as well. If
$\chi(g_{\bar{v}_0})=1$, then either $\tT_{M,\si_{can}}(\bar{\chi})=0$,
or for exactly $\nu-2$ indices $i$ ($1\leq i\leq \nu$)  one has
$\chi(g_{\bar{v}_i})=1$.  In this later case the limit is
non-zero. Let us analyze this case more closely.
Assume that $\chi(g_{\bar{v}_i})\not=1$ for $i=1,2$.
Using the last statement of the previous subsection,
it is not difficult to verify that $\vec{w}^*$ computed from any
vertex on these two arms  provide the same limit.
In fact, by the same argument (cf \ref{le}), one gets that even the central
vertex $\bar{v}_0 $ provides a suitable set of weight
$\vec{w}^*$ (for any $\chi$).
The relevant weights  can be computed via \ref{ss: s4}(1),
and  with the notation $\alpha:=\mbox{lcm}(\alpha_1,\ldots,\alpha_\nu)$
 one has:
\[
\hat{\tT}_{M,\si_{can}}(\bar{\chi})=
\lim_{t\to 1}\,
\frac{\bigl(t^{\alpha}\chi(g_{\bar{v}_0})-1\bigr)^{\nu-2}}
{\prod_{i=1}^\nu\ \bigl(t^{\alpha/\alpha_i}\,
\chi(g_{\bar{v}_i})-1\bigr)} \ \ \ \mbox{for any } \ \chi\in \hat{H}\setminus
\{1\}.\]
Notice the mysterious similarity of this expression
 with the Poincar\'e series of the
graded affine ring associated  with the universal abelian cover of
$(X,0)$, provided that $(X,0)$ admits a good $\C^*$--action,
cf \cite{Neu}.

\section{Brieskorn--Hamm  rational homology spheres}

\gtsubsection{Notation}
\label{ss: b1} Fix $n\geq 3$
 positive integers $a_i\geq 2$ ($i=1,\ldots,n$).
For any set of complex numbers $C=\{c_{j,i}\}, i=1,\ldots, n$;
$j=1,\ldots, n-2$, one can consider the affine variety
\[
X_C(a_1,\ldots,a_n):=\left\{z\in \C^n\, :\, c_{j,1}z_1^{a_1}+
\cdots +c_{j,n}z_n^{a_n} =0 \mbox{ for all } 1\leq j\leq
 n-2\right\}.
\]
It is well-known (see \cite{Hamm})
that for generic $C$, the intersection of
$X_{C}(a_1,\ldots,a_n)$  with the unit sphere $S^{2n-1}\subset
\C^n$ is an oriented  smooth 3--manifold whose diffeomorphism type
is independent of the choice of the coefficients $C$. It is
denoted by $M=\Sigma(a_1,\ldots,a_n)$.

\subsection{}\hspace{-0.25em}
\label{ss: b2}
In fact (cf \cite{JN,NeR}), $M$ is an oriented  Seifert 3--manifold with Seifert  invariants
\[
\Bigl(g; \underbrace{(\alpha_1,\beta_1),\cdots
,(\alpha_1,\beta_1)}_{s_1}; \cdots;
\underbrace{(\alpha_n,\beta_n),\cdots
,(\alpha_n,\beta_n)}_{s_n}\Bigr)
\]
where $g$ denotes the  genus of the base of the Seifert fibration,
and the pairs of coprime positive integers $(\alpha_i,\beta_i)$
(each considered $s_i$ times) are the orbit invariants (cf
\ref{ss: s1} and \ref{ss: s2} for notations). Recall that the rational degree
 of this  Seifert fibration is
\begin{equation*}
e=-\sum_{i=1}^n\, s_i\cdot \frac{\beta_i}{\alpha_i}<0.
\tag{e}
\label{b2e}
\end{equation*}
Set
\[
 a:={\rm lcm}\,(a_i;\;\;1\leq i\leq n),\;\;
q_i:=\frac{a}{a_i},\;\; 1\leq i\leq n; \ \;\;A:=\prod_{i=1}^na_i.
\]
\noindent The  Seifert invariants   are as follows
(see \cite[Section 7]{JN} or \cite{NeR}):
\begin{equation*}
\alpha_i:=\frac{a}{{\rm lcm}\,(a_j;\;\;j\neq i)},
\ \ \ \
s_i:= \frac{\prod_{j\neq i}a_j}{{\rm lcm}\,(a_j;\;\;j\neq
i)}=\frac{A\alpha_i}{a a_i},\ \ \ \
(1\leq i\leq n);
\end{equation*}
\begin{equation*}
g:=\frac{1}{2}\Bigl(2 +(n-2)\frac{A}{a}-\sum_i s_i\Bigr). \tag{g}
\end{equation*}
By these notations $-e=A/(a^2)$. Notice that the integers
$\{\alpha_j\}_{j=1} ^n$ are pairwise coprime. Therefore the
integers $\beta_j$ are  determined from (e). In fact $\sum _j\
q_j\beta_j =1$, hence $\beta_iq_j\equiv 1 \ (\mmod \ \alpha_j)$
for any $j$. Similarly as above, we set $\alpha:={\rm
lcm}(\alpha_1,\ldots,\alpha_n)$.

It is clear that $M$ is a ${\Q}$--homology sphere if and only if
$g=0$. In order to be able to compute the Reidemeister--Turaev torsion,
we need a good characterization of $g=0$ in terms of the integers $\{a_i\}_i$.
For hypersurface singularities this is given in \cite{Bries}. This
characterization was partially extended for complete intersections in
\cite{Hamm}. The next proposition provides a complete characterization (for the case when the link is 3--dimensional).

\gtsubsection{Proposition}\label{ss: bicis} {\sl Assume
that $X_C(a_1,\ldots,a_n)$
is a Brieskorn--Hamm  isolated complete intersection singularity as in
\ref{ss: b1} such that its link $M$ is a 3--dimensional  rational homology
sphere. Then $(a_1,\ldots, a_n)$ (after a possible permutation) has (exactly)
one of the following forms:

{\rm(i)}\qua  $(a_1,\ldots,a_n)=(db_1,db_2,b_3,\ldots,b_n)$,
   where the integers $\{b_j\}_{j=1}^n$ are
pairwise coprime, and $\gcd(d,b_j)=1$ for any $j\geq 3$;

{\rm(ii)}\qua  $(a_1,\ldots,a_n)=(2^cb_1,2b_2,2b_3,b_4\ldots,b_n)$,
   where the integers $\{b_j\}_{j=1}^n$ are odd and pairwise coprime,
and $c\geq 1$. }
\begin{proof}  The proof  will be carried out in several steps.

\smallskip

\noindent {\bf Step 1}\qua Fix any four distinct indices
 $i,j,k,l$. Then $d_{ijkl}:=\gcd(a_i,a_j,a_k,a_l)=1$. Indeed,
if a prime $p$ divides $d_{ijkl}$, then $p^2|(A/a)$ and $p^2|s_i$ for all $i$.
Hence by \ref{ss: b2}(g) one has $p^2|2$.

\smallskip

\noindent {\bf Step 2}\qua  Fix any three distinct indices $i,j,k$
and set $ d_{ijk}:=\gcd(a_i,a_j,a_k)$. If a prime $p$ divides
$d_{ijk}$ then $p=2$.  Indeed, if  $p|d_{ijk}$ then $p| (A/a)$ and
$p| s_j$ for all $j$, hence $p|2$ by \ref{ss: b2}(g).

\smallskip

\noindent {\bf Step 3}\qua There is at most one triple
 $i< j<k$ with $d_{ijk}\not=1$.
This follows from steps 1 and 2.

\smallskip

\noindent {\bf Step 4}\qua Assume that $d_{ijk}=1$ for all triples $i,j,k$.
For any $i\not=k$ set $d_{ik}:=\gcd(a_i,a_k)$. Then $A/a=\prod_{i<k}
d_{ik}$, and there are similar identities for each $s_j$. Then
\ref{ss: b2}(g) reads as
\begin{equation*}
2+(n-2)\prod_{i<k}d_{ik}-\sum_{j=1}^n\prod_{i<k;i,k\neq j}d_{ik}=0.
\tag{eq1}
\end{equation*}

\noindent {\bf Step 5}\qua Assume that $d_{123}\not=1$,
hence $$(a_1,\ldots,a_n)=(2^cb_1,2^ub_2,2^vb_3,b_4\ldots,b_n),$$ with
$c\geq  u\geq v$, and all $b_i$ odd numbers  (after a permutation
of the indices). Then $u=v=1$.
For this use similar argument as above with
 $4|(A/a)$, $4|s_i$ for $i\geq 4$,
$s_1$ and $s_2$ (resp. $s_3$) is divisible exactly by the $v^{th}$
(resp. $u^{th}$) power of 2.

Using this fact, write for each par $i\not=k$, $d_{ik}:=\gcd(b_i,b_k)$.
We deduce as above that  $A/a=4\prod_{i<k}d_{ik}$, and there are
similar identities  for each $s_j$. Then \ref{ss: b2}(b) transforms into
\begin{equation*}
\frac{1}{2} +(n-2)\prod_{i<k}d_{ik}-\sum_{j=1}^n\varepsilon_j\prod_{
i<k;i,k\not=j} d_{ik}=0,
\tag{eq2}
\end{equation*}
where $\varepsilon_j=1/2$ for $j\leq 3$ and $=1$ for $j\geq 4$.

\smallskip

\noindent {\bf Step 6}\qua The equation (eq2) has only one solution with
all $d_{ik}$ strict positive integer, namely $d_{ik}=1$ for all $i,k$.
Similarly, any  set of solutions of (eq1) has at most one $d_{ik}$
strict greater than 1, all the others being equal to 1.

This can be proved eg  by induction. For example, in the case of
(eq2), if one replaces the set of integers $d_{ik}$ in the left hand
side of the equation with the same set  but in which one of them is increased by one
unit, then the new expression is strictly greater than the old one. A similar
argument works for (eq1) as well. The details are left to the reader.
\end{proof}

\gtsubsection{Verification of the conjecture in the case \ref{ss: bicis}(i)}
\label{(i)} We start to list the properties of Brieskorn--Hamm complete
intersections of the form (i).

\noindent  $\bullet$\qua $\alpha_j=b_j$ for $j=1,\ldots, n$;

\noindent  $\bullet$\qua $s_1=1, \ s_2=1$, and  $s_j=d$ for $j\geq 3$;
in particular, the number of ``arms'' is $\nu=2+(n-2)d$;

\noindent  $\bullet$\qua $\alpha=B:=\prod_{j=1}^nb_j$, and $-e \cdot B=1$,
hence  by \ref{ss: s1} the generic orbit  $h$  is homologically
trivial.

\noindent $\bullet$\qua Using the group representation \ref{ss: s2}, and the fact
that $h$ is trivial, one has
\[
H=\mbox{ab}\langle \ g_{ij},\;\;1\leq j\leq n,\;\;1\leq i\leq s_j\
| \ g_{ij}^{\alpha_j}=1 \ \mbox{for all } \ i,j;
\;\;\prod_{i,j}g_{ij}^{\omega_j}=1 \ \rangle.
\]
Since the integers $\alpha_j$ are pairwise coprime, taking the
$\alpha/\alpha_j$ power of the last relation, and using that
$\gcd(\alpha_j, \omega_j)=1$, one obtains that
\[
H=\bigoplus_{j\geq 3} \mbox{ab}\langle \ g_{ij},\;\;1\leq i\leq d
\ | \ g_{ij}^{b_j}=1 \ \mbox{for all } \ i; \;\;\prod_{i}g_{ij}=1
\ \rangle\ \cong\bigoplus_{j\geq 3}(\Z_{b_j})^{d-1}.
\]
In particular, $|H|= \prod_{j\geq 3} \, b_j^{d-1}$.

\vspace{2mm}

$\bullet$\qua {\bf The Reidemeister--Turaev torsion of $M$}\qua
By \ref{ss: sm}
\begin{equation*}
\hat{\tT}_{M,\si_{can}}(\chi)= \lim_{t\to 1}\
\frac{(t^{\alpha}-1)^{d(n-2)}}{ (t^{\alpha/\alpha_1}-1)
(t^{\alpha/\alpha_2}-1) \prod_{j\geq 3} \prod_{i=1}^d
(t^{\alpha/\alpha_j}\bar{\chi}(g_{ij})-1)}.
\label{b4t}
\end{equation*}
As  explained in \ref{ss: sm}, for a fixed $\chi$, the
expression $\hat{\tT}_{M,\si_{can}}(\chi)$ is nonzero if and only
if for exactly two pairs $(i,j)$ (where $j\geq 3$ and $1\leq
i\leq b$) $\chi(g_{ij})\not=1$. Analyzing the group structure, one
gets easily that these two pairs must have the same $j$. For a
fixed $j$, there are $d(d-1)/2$ choices for the  set of indices
$\{i_1,i_2\}$. For fixed $j$, the set of nontrivial characters of
the group
\[
\mbox{ab}\langle \
g_{ij},\;\;1\leq i\leq d \ | \ g_{ij}^{b_j}=1 \ \mbox{for all } \
i; \;\;\prod_{i}g_{ij}=1 \ \rangle,
\]
satisfying $\chi(g_{ij})=1$ for all $i\not=i_1,i_2$ is clearly
$\hat{\Z}_{b_j}\setminus \{1\}$, and in this case
 $\chi(g_{i_1j})\chi(g_{i_2j})=1$ as well.
Therefore
\[
\tT_{M,\si_{can}}(1)=\frac{1}{|H|}\cdot
\sum_{j\geq 3}\ \frac{\alpha^{d(n-2)}}{ \frac{\alpha}{\alpha_1}
\frac{\alpha}{\alpha_2} (\frac{\alpha}{\alpha_3})^d\cdots
(\frac{\alpha}{\alpha_j})^{d-2}\cdots
(\frac{\alpha}{\alpha_n})^d}
\]
\[
\times  \frac{d(d-1)}{2}\cdot {\sum_{\Z_{b_j}}}'\frac{1}{(\zeta-1)(\bar{\zeta}-1)}.
\]
 Recall that $|H|=\prod_{j\geq 3}b_j^{d-1}$. Hence, by (\ref{eq: trig3}) of
Appendix B and an easy  computation:
\[
\tT_{M,\si_{can}}(1)=\frac{B\cdot d(d-1)}{24}\sum_{j\geq 3}\bigl(1-\frac{1}{
b_j^2}\bigr).
\]
\noindent $\bullet$\qua {\bf The Casson--Walker invariant}\qua
From \ref{ss: s3} we get
\begin{equation*}
-\frac{\lambda(M)}{|H|}= -\frac{B}{24}\Bigl(-d(n-2)+\sum_{j=1}^n
\frac{s_j}{b_j^2}\Bigr)-\frac{1}{24B}+\frac{1}{8}+\frac{1}{2}
\sum_{j=1}^ns_j\bms(\beta_j, b_j).
\end{equation*}
\noindent $\bullet$\qua {\bf The signature of the Milnor fiber}\qua Since $X_C$
is an isolated complete intersection singularity, its singular point is
Gorenstein. Hence, by \ref{ss: 4.4}, it is enough to verify only part
(3) of the main conjecture, part (2) will follow automatically.

The signature  $\si(F)=\si(a_1,\ldots, a_n)$ of
the Milnor fiber $F$ of a Brieskorn--Hamm
singularity is computed by Hirzebruch \cite{Hi} in terms of
cotangent sums. Nevertheless, we will use the version proved in
\cite[1.12]{NW}.
This, in the case (i),  (via \ref{eq: trig4})  reads as
\begin{equation*}
\si(F)=-1
+\frac{1}{3B}\Bigl(1-(n-2)d^2B^2+B^2\sum_{j=1}^n\frac{s_j^2}
{b_j^2}\Biggr)-4\sum_{j=1}^ns_j\bms(q_j,b_j),
\end{equation*}
where $q_j=s_jB/b_j$ for all $1\leq j\leq n$. Since
$\beta_jq_j\equiv 1 \ (\mmod \ b_j)$ (cf \ref{ss: b2}), one has
$\bms(q_j,b_j)=\bms(\beta_j,b_j)$ for all $j$. Now, by a simple
computation  one can verify the conjecture.

\gtsubsection{Verification of the conjecture in the case \ref{ss: bicis}(ii)}
\label{(ii)}  The discussion is rather similar to the previous case, the only
difference (which is not absolutely negligible) is that now $h$ is not
trivial. This creates some extra work in the torsion computation.
In the sequel we write $B:=\prod_jb_j$.

\noindent $\bullet$\qua $a=2^cB$ and $A=2^{c+2}B$. Moreover,
 $\alpha_1=2^{c-1}b_1$ and $\alpha_j=b_j$ for $j\geq 2$.

\noindent  $\bullet$\qua $s_j=2$ for $j\leq 3$ and $s_j=4$ for $j\geq 4$.
The number of ``arms'' is $\nu=4n-6$.

\noindent  $\bullet$\qua $\alpha=2^{c-1}B$, hence  $-e^{-1}=2^{c-2} B$.
Therefore $|\langle h\rangle|=2$ and $|H|=2^cB^3/(b_1b_2b_3)^2$.

\noindent $\bullet$\qua The self intersection number $b$ of the central
exceptional divisor is even. Indeed, equation (e) implies that
$1\equiv \beta_1b_2\cdots b_n\ (\mmod \ \alpha_1)$. Since $\omega_1\equiv
-\beta_1\ (\mmod\ \alpha_1)$, one has $2^{c-1}|1+\omega_1b_2\cdots b_n$.
This, and the first formula of \ref{ss: s2} implies that $b$ is even.

\noindent $\bullet$\qua Using \ref{ss: s2}, and the fact
that $h^{-b}=1$ is automatically satisfied,  we obtain the following presentation for $H$:
\[
\mbox{ab}\Bigl\langle g_{ij},\;\;1\leq j\leq n,\;\;1\leq i\leq s_j;\  h\Bigl|
\; g_{ij}^{\alpha_j}=h \ \mbox{for all } \; i,j;
\;\;\prod_{i,j}g_{ij}^{\omega_j}=1, \ h^2=1 \Bigr\rangle.
\]
Clearly $\langle h\rangle\approx \Z_2$, and there is
an exact sequence $0\to \langle h\rangle \to H\to Q\to 0$ with
\[
Q=\bigoplus_{j\geq 1} \mbox{ab}\langle \ g_{ij},\;\;1\leq i\leq s_j
\ | \ g_{ij}^{\alpha_j}=1 \ \mbox{for all } \ i; \;\;\prod_{i}g_{ij}=1
\ \rangle\ \cong\bigoplus_{j\geq 1}(\Z_{\alpha_j})^{s_j-1}.
\]

\vspace{2mm}

\noindent  $\bullet$\qua {\bf The Reidemeister--Turaev torsion of $M$}\qua
We have to distinguish two types of characters $\chi\in \hat{H}\setminus
\{1\}$ since $\chi(h)$ is either $+1$ or $-1$.
The sum over characters $\chi$ with $\chi(h)=1$
(ie, over $\hat{Q}\setminus \{1\}$) can be computed similarly as in case
(i), namely it is
\[
\frac{1}{|H|}\cdot\!
\sum_{j\geq 1}\ \frac{\alpha^{\nu-2}}{ (\frac{\alpha}{\alpha_1})^{s_1}
\cdots (\frac{\alpha}{\alpha_{j-1}})^{s_{j-1}}
(\frac{\alpha}{\alpha_j})^{s_j-2}
(\frac{\alpha}{\alpha_{j+1}})^{s_{j+1}}\cdots
(\frac{\alpha}{\alpha_n})^{s_n}}\cdot\!
{\sum_{\Z_{\alpha_j}}}'\frac{s_j(s_j-1)}{2(\zeta-1)(\bar{\zeta}-1)}
\]
\[
=\frac{2^{c-1}B}{24}\sum_{j\geq 1}\frac{s_j(s_j-1)}{2}\bigl(1-\frac{1}{
\alpha_j^2}\bigr).
\]
The sum over characters $\chi$ with $\chi(h)=-1$ requires no ``limit
regularization'', hence it is $(-2)^{\nu-2}/(\prod_jP_j)$, where for any
fixed  $j$ the expression
$P_j$ has the form $\prod_*(\zeta_i-1)$,
where the product runs over $1\leq i\leq s_j$,
$\zeta_i^{\alpha_j}=-1$ with restriction $\prod_i\zeta_i=1$.

Using the identity $-2=\zeta_i^{\alpha_j}-1$ one gets
$$\frac{1}{P_j}=\frac{1}{(-2)^{s_j}}{\prod}_*(1+\zeta_i+\cdots +
\zeta_i^{\alpha_j-1}).$$
By an elementary argument, this is exactly $(\alpha_j/2)^{s_j}$.
Therefore, this second contribution  is $2^{c-1}B/8$, hence
$$\tT_{M,\si_{can}}(1)=\frac{2^{c-1}B}{8}+
\frac{2^{c-1}B}{24}\sum_{j\geq 1}\frac{s_j(s_j-1)}{2}\bigl(1-\frac{1}{
\alpha_j^2}\bigr).
$$
\noindent $\bullet$\qua {\bf The Casson--Walker invariant}\qua
From \ref{ss: s3} and $\nu=4n-6$ one gets
\begin{equation*}
-\frac{\lambda(M)}{|H|}= -\frac{2^{c-2}B}{24}\Bigl(-4(n-2)+\sum_{j=1}^n
\frac{s_j}{\alpha_j^2}\Bigr)-\frac{1}{3\cdot 2^{c+1}B}+\frac{1}{8}+\frac{1}{2}
\sum_{j=1}^ns_j\bms(\beta_j, \alpha_j).
\end{equation*}
\noindent $\bullet$\qua {\bf The signature of the Milnor fiber}\qua
Using the above identities about the Seifert invariants
(and  $a_j=4\alpha_j/s_j$ too),
 \cite[1.12]{NW} reads as

$\si(F)=$\vspace{-4mm}
$$-1
+\frac{1}{3\cdot 2^{c-2}B}\Bigl(1-(n-2)2^{2c}B^2+2^{2c-4}B^2\sum_{j=1}^n
\frac{s_j^2}{\alpha_j^2}\Biggr)-4\sum_{j=1}^ns_j\bms(\beta_j,\alpha_j).
$$
Now, the verification of the statement of the conjecture is elementary.

\section{Some rational singularities}

\gtsubsection{Cyclic quotient singularities}
\label{ss: 4.8} The link of a cyclic quotient singularity
$(X_{p,q},0)$ ($0<q<p,\; (p,q)=1$) is the lens space $L(p,q)$.
$X_{p,q}$ is numerically Gorenstein if and only if $q=p-1$, case which will be considered in \ref{ss: 4.7}. In all other cases $\si_{can}$
is not spin. In all the cases $p_g=0$. Moreover (see eg
\cite[5.9]{LW},   or  \cite{Ishii}, or \ref{ss: 5.2}):
\[
K^2+\#{\mathcal V}= \frac{2(p-1)}{p}-12\cdot \bms(q,p).
\]
On the other hand, the Seiberg--Witten invariants of $L(p,q)$  are
computed in \cite{Nico3} (where a careful reading will identify
$\ssw^0_M(\si_{can})$ as well). In fact, cf (\cite[3.16]{Nico3}):
\[
\hat{\tT}_{M,\si_{can}}(\chi)= \frac{1}{(\bar{\chi}-1)(\bar{\chi}^q-1)},
\]
 fact
which follows also from \ref{ss: th1}. Therefore, using
\cite[18a]{RG} (or \ref{eq: trig2bis}), one gets that
\[
\tT_{M,\si_{can}}(1)=\frac{p-1}{4p}-\bms(q,p).
\]
The Casson--Walker contribution is $\lambda(L(p,q))/p=\bms(q,p)/2$
(cf \ref{ss: 5.3}).
Hence one has equality in part (1) of the Conjecture.

\gtsubsection{Particular case: the  $A_{p-1}$--singularities}
\label{ss: 4.7}
Assume that $(X,0)\!=\!(\{x^2+y^2+z^p= 0\},0)$. Then by \ref{ss: 4.8} (or 
\cite[Section 2.2.A]{Nico3}) one has $\ssw^0_M(\eps_{can})=(p-1)/8$ (in \cite{Nico3} this spin structure is denoted by $\si_{spin}$). On the other hand, $\bif_F$ is negative definite of rank $p-1$, hence $\si(F)=-(p-1)$. (Obviously, the $A_{p-1}$--case
can also  be deduced from section 6.)

\gtsubsection{The $D_n$--singularities}
\label{ss: d}
For each $n\geq 4$, one denotes
by $D_n$ the singularity at the origin of the weighted homogeneous
complex hypersurface $x^2y+ y^{n-1} +z^2=0$. It is convenient to
write $p:=n-2$. We invite the reader to recall the notations of
\ref{ss: sws} about orbifold invariants.

The normalized Seifert invariants are
\[
(b,(\omega_1,\alpha_1), (\omega_2,\alpha_2),(\omega_3,\alpha_3))=(-2, (p,p-1),
(2,1),(2,1)).
\]
 Its rational degree is $\ell=-1/p$.  Observe that $ \kappa=\ell$.

The link $M$ is the unit circle bundle of the $V$--line bundle
${\bL}_0$ with rational degree $\ell$ and singularity data
$((p-1)/p,1/2,1/2)$. Therefore, ${\bL}_0=K_\Sigma$. Hence,
$\rho(\si_{can})=\{\kappa/2\ell\}=0$. The canonical representative
of $\si_{can}$ is then the trivial line bundle $E_0$. It has
rational degree $0$ and singularity data $\vec{\gamma}=(0,0,0)$.
The Kreck--Stolz invariant is then
$$KS_{M}(\si_{can},g_0,0) =7-4\sum_{i=1}^3
\bms(\omega_i,\alpha_i)-8\sum_{i=1}^3\bms(\omega_i,\alpha_i;
\frac{\omega_i} {2\alpha_i}, -1/2)-4\sum_{i=1}^3
\Bigl\{\frac{1}{2\alpha_i}\Bigr\}.$$
Using the fact that
$\bms(1,2; 1/4,-1/2)=0$, this expression equals:
\[
6+\frac{p}{3}+\frac{2}{3p}+ 8\bms(1,p,\frac{1}{2p},1/2).
\]
Now using the reciprocity formula for the generalized Dedekind sum, one has
\[
8\bms(1,p,\frac{1}{2p},1/2)=\frac{4/3-2p+2p^2/3}{p}=\frac{4}{3p}-2+
\frac{2p}{3}.
\]
Thus
\[
8\ssw^0_{M}(\si_{can})=KS_{M}(\si_{can},g_0)=
4+\frac{p}{3}-\frac{4}{3p}+\frac{4}{3p}-2+\frac{2p}{3}=2+p.
\]
On the other hand, the signature of the Milnor fiber is
$-n=-(p+2)$, confirming again the Main Conjecture.

\gtsubsection{The $E_6$ and $E_8$ singularities}
\label{ss: ee}
Both $E_6$ (ie, $x^4+y^3+z^2=0$) and $E_8$ (ie, $x^5+y^3+z^2=0$) are
Brieskorn (hypersurface) singularities, hence the result of section 6 can be
applied.
The link of $E_8$ is an
\emph{integral} homology sphere,  hence the validity of
the conjecture in this case was proved in \cite{FS}.
The interested reader can verify the conjecture  using the machinery
of  \ref{ss: sws} and \ref{ss: sws2} as well.

\gtsubsection{The $E_7$ singularity}
\label{ss: e}
It is given by the complex hypersurface $ x^3+xy^3+z^2=0$. The group $H$ is $\Z_2$. The
normalized Seifert invariants are $(-2, (2,1), (3,2), (4,3))$, the
rational degree is $-1/12$. We deduce as above  that${\bL}_0=K_\Sigma$,
with $\rho_0=\rho(\si_{can})=1/2$. The canonical  representative
is again the trivial line bundle $E_0$. Its singularity data are
trivial. The Seiberg--Witten invariant of $\si_{can}$ is determined
by the Kreck--Stolz invariant alone. A direct computation shows
that $KS_{M}(\si_{can},g_0)=7$.
But the signature of the Milnor fiber is $\si(E_7)=-7$ as well,
hence the statement of the conjecture is true.

\gtsubsection{Another family  of  rational singularities}
\label{ss: rat}
Consider a singularity\break $(X,0)$ whose link $M$ is described by the
negative definite plumbing given in Figure 1. (It is clear that in this
case  $M$  it is {\em not}  numerically Gorenstein.)
\begin{figure}[ht!]
\centerline{\small
%\ShowGrid 
\SetLabels 
\E(0*0.1){$-2$}\\
\E(0.35*0.1){$-2$}\\
\E(0.65*0.1){$-2$}\\
\E(.98*0.1){$-2$}\\
\E(.55*0.35){$-2$}\\
\E(.55*.98){$-2$}\\
\E(0.47*0.1){$-3$}\\
\E(0.2*0.1){$m-1$}\\
\E(0.8*0.1){$m-1$}\\
\E(0.35*-0.05){$(m-1)x$}\\
\E(0.65*-0.05){$(m-1)x$}\\
\E(0.5*-0.05){$mx$}\\
\E(0.4*0.35){$(m-1)x$}\\
\E(0.55*0.7){$m-1$}\\
\E(0.45*.98){$x$}\\
\endSetLabels 
\AffixLabels{{\includegraphics[width=4in]{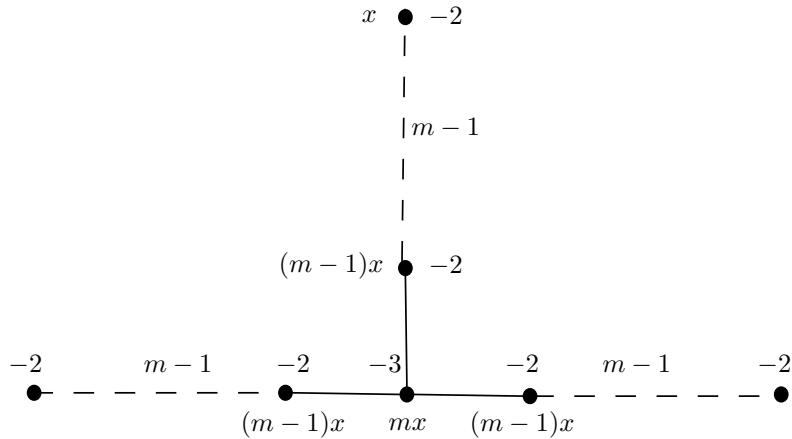}}}}
\vspace{4mm}
\caption{The resolution graph of the rational singularity $(X,0)$}
\end{figure}

\noindent The number of $-2$ spheres on any branch is  $m-1$,
where  $m\geq 2$. It is easy to verify that the $(X,0)$ is a
rational (with Artin cycle $\sum_vE_v$) (see, eg \cite{Nem1}).
$M$ is   Seifert  manifold with normalized Seifert invariants
$(-3, (m,m-1),(m,m-1),(m,m-1))$ and rational degree $l=-3/m$.

To compute  the Seiberg--Witten  invariant of $M$ associated with
$\si_{can}$ we use again \ref{ss: sws2}.

The canonical $V$--line bundle of $\Sigma$   has singularity data
$(m-1)/m$ (three times) and  rational degree $\kappa=1+\ell$. The
link $M$ is the unit circle bundle of the $V$--line bundle
${\bL}_0$ with rational degree $\ell$ and singularity data
$(m-1)/m$ (three times). Therefore
$$\rho_0=\Bigl\{-\frac{m-3}{6}\Bigr\}.$$
To apply \ref{ss: sws2} we need $\rho_0\neq 0$, ie,
\[
m\not\equiv 3\ (\mmod \ 6).
\]
The canonical representative   of  $\si_{can}$ is the $V$--line
bundle $E_0$ with
\[E_0= n_0{\bL}_0,\;\; n_0=\Bigl\lfloor\frac{3-m}{6}\Bigr\rfloor.
\]
The reader is invited to recall the definition of $S_0^\pm$. We
start with the computation of $S_0^+$. Notice that
\[
-\frac{1}{2}\deg^VK_\Sigma \leq \nu(E)<0\Longleftrightarrow 0\leq
n\ell<\frac{1}{2}\deg^VK_\Sigma.
\]
Hence
\begin{equation*}
\frac{1}{2\ell}\deg^VK_\Sigma < n \leq  0, \ \ \mbox{ie,} \ \ \
-\frac{m-3}{6}< n \leq 0.
\tag{1}
\label{rat1}
\end{equation*}
The singularity data of $n{\bL}_0$ are all equal to $\{-n/m\}$
(three times). We deduce
\[
\deg| n{\bL}_0|=\deg^V( n{\bL}_0)- 3 \Bigl\{ \frac{-n}{m}\Bigr\}=
3\Bigl\lfloor\frac{-n}{m}\Bigr\rfloor.
\]
Now observe that  (1) implies
\[
0\leq \frac{-n}{m} <\frac{m-3}{6m}, \ \mbox{hence} \
\Bigl\lfloor\frac{-n}{m}\Bigr\rfloor=0
\]
for every $n$ subject to the condition (1). Since $m\not\equiv 3\
(\mmod \ 6)$,  we deduce
\begin{equation*}
|S_0^+|=\Bigl\lfloor\frac{m-3}{6}\Bigr\rfloor+1=-n_0.
\tag{2}
\label{rat2}
\end{equation*}
Moreover, all the connected components corresponding to the
elements in $S_0^+$ are points. \noindent Similarly, the condition
$ 0<\nu(E)\leq\frac{1}{2}\deg^VK_\Sigma$ implies
\[
\frac{1}{2}\deg^V K_\Sigma< n\ell \leq
\deg^VK_\Sigma\Longrightarrow \frac{1}{\ell}\deg^VK_\Sigma \leq n
< \frac{1}{2\ell}\deg^VK_\Sigma.
\]
Hence
\begin{equation*}
 -\frac{m-3}{3}\leq n <-\frac{m-3}{6}  .
\tag{3}
\label{rat3}
\end{equation*}
The singularity data of the $V$--line bundle  $K_\Sigma -n{\bL}_0$
are all equal to $ \{ (n-1)/m\}$. We deduce
\[
\deg|K_\Sigma -n{\bL}_0|= 1+3\frac{n-1}{m} -3\{(n-1)/m\}=
1+3\Bigl\lfloor\frac{n-1}{m}\Bigr\rfloor.
\]
But this number is negative (because of (3)), hence
$S_0^-=\emptyset$.  These considerations show that Proposition
\ref{ss: sws2} is applicable. Set
\[
-\frac{m-3}{6}=-k+\rho_0, \;\;0<\rho_0<1, \ \ \mbox{$k$ non-negative integer}.
\]
Then  $ n_0=-k$. The canonical representative  is
$E_0= -k {\bL}_0$. It has  degree  $-k\ell$. Its singularity data
are all equal to $\gamma_i/\alpha_i=k/m$. Then in the formula of
$KS_M$ one has $\omega_i=m-1$, $r_i=-1$, $\gamma_i=k$ for all $i$.
Hence
\[
 KS_M(\si_{can})
=\ell+1-4\ell\rho_0(1-\rho_0)+12\rho_0+2(3+\ell)(1-2\rho_0)
-12\Bigl\{\frac{\rho_0-k}{m}\Bigr\}
\]\vspace{-4mm}
\[
-12\bms(m-1,m)-24\bms(m-1, m;\frac{k+\rho_0(m-1)}{m},-\rho_0) .
\]
Observe now that
\[
-\bms(m-1,m)=\bms(1,m)= \frac{m}{12}+\frac{1}{6m}-\frac{1}{4},
\]\vspace{-4mm}
\[
-\bms(m-1,m;\frac{k+\rho_0(m-1)}{m},-\rho_0))=-\bms(-1,m,\frac{k+\rho_0(m-1)}{m}-\rho_0,-\rho_0)
\]\vspace{-4mm}
\[
=\bms(1,m,\frac{\rho_0-k}{m},-\rho_0).
\]
Moreover, from the definition of  Dedekind sum we obtain
\[
\bms(1,m,\frac{\rho_0-k}{m},-\rho_0)=\bms(1,m,\frac{k-\rho_0}{m},0)=
\bms(1,m)+\frac{k(k-1)}{2m}-\frac{k-1}{2}.
\]
Finally, by an elementary but tedious computation we get
\[
KS_M(\si_{can})=3m-\frac{m}{3}-2 -8k.
\]
The Seiberg--Witten invariant  of the canonical $spin^c$ structure
is then
\[
8\ssw_M^0(\si_{can})=KS_M(\si_{can})+8|S_0^+|=KS_M(\si_{can})+8k=3m-\frac{m}{3}-2.
\]
The coefficients of $Z_K$ are labelled on the graph,
 where  the unknown $x$ is determined from the adjunction formula applied
to the central $-3$--sphere; namely $ -3mx +3(m-1)x=-1$, hence
$x=1/3$. Then $Z_K^2=\sum r_v(e_v+2)=-r_0=-m/3$. The number of
vertices of this graph is $3m-2$ so
\begin{equation*}
8p_g+K^2+\#{\mathcal V}=3m-2-\frac{m}{3}.
\tag{4}
\end{equation*}
This confirms once again the Main Conjecture.

\gtsubsection{The case $m=3$}
\label{ss: mmm} In the previous example we verified the
conjecture for all $m\not\equiv 3 \ (\mmod \ 6)$. For the other
values the method given by \ref{ss: sws2} is not working. But this
fact does not contradict the conjecture. In order to show this, we
indicate briefly how one can verify the conjecture in the  case
$m=3$ by the torsion computation.

In this case $|H|=27$ and $h$ has order 3. First consider the set
of characters $\chi$ with $\chi(g_{\bar{v}_0})=1$ (there are 9
altogether). They satisfy $\prod_i\chi(g_{\bar{v}_i})=1$. If
$\chi(g_{\bar{v}_i})\not=1$ for all $i$ (2 cases), or if $\chi=1$
(1 case) then $\hat{\tT}(\chi)=0$. If $\chi(g_{\bar{v}_i})=1$ for
exactly one index $i$, then the contribution in
$\sum_{\chi}\hat{\tT}(\chi)$ is 2 for each choice of the index,
hence altogether 6.

Then, we consider those characters $\chi$ for which
$\chi(g_{\bar{v}_0})\not=1$ (18 cases). Then one has to compute
the sum
$$\sum \frac{1-\zeta_1^3}{(1-\zeta_1)(1-\zeta_2)(1-\zeta_1^{-1}\zeta_2^{-1})},
$$
where the sum runs over $\zeta_1,\zeta_2\in \Z_9,\
\zeta_1^3=\zeta_2^3\not=1$. A computation shows that  this is 9.
Therefore, $\tT_M(1)/|H|=(6+9)/27=5/9$.

The Casson--Walker invariant can be computed easily from the
Seifert invariants, the result is $\lambda(M)/|H|=-7/36$.
Therefore, the Seiberg--Witten invariant is 5/9+7/36=3/4. But this
number equals $(K^2+\#{\mathcal V})/8$ (cf \ref{ss: rat}(4) for $m=3$
and $p_g=0$).

\section{Some minimally elliptic  singularities}

\gtsubsection{``Polygonal'' singularities}
\label{ss: me1}   Let $(X,0)$ be a normal surface singularity with
resolution graph given by Figure 2.
\begin{figure}[ht!]
\centerline{\small
%\ShowGrid 
\SetLabels 
\E(0.4*0.55){$E_0$}\\
\E(0.38*0.45){$2-\nu$}\\
\E(0.9*0.9){$E_2$}\\
\E(0.95*0.8){$-a_2$}\\
\E(1.05*0.55){$E_1$}\\
\E(1.05*0.45){$-a_1$}\\
\E(.95*0.2){$E_\nu$}\\
\E(.9*0.1){$-a_\nu$}\\
\endSetLabels 
\AffixLabels{{\includegraphics[width=1.8in]{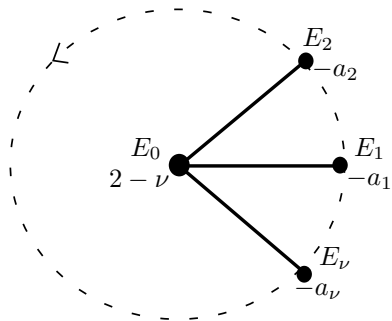}}}}
\vspace{1mm}
\caption{The resolution graph of a ``polygonal'' singularity}
\end{figure}

Here  we assume that $\nu>2$. The negative definiteness of the
intersection form implies that the integers $(a_1,\cdots,a_\nu)\in
{\bZ}^\nu_{>1}$ satisfy
\[
\ell:= 2-\nu+ \sum_{i=1}^\nu\frac{1}{a_i}<0.
\]
An elementary computation shows that $Z_K=2 E_0+\sum_{i=1}^\nu
E_i$. If $\nu>3$ then this cycle is exactly the minimal cycle
$Z_{min}$ of Artin. If $\nu=3$ then the graph is not
 minimal, but after blowing down the central irreducible exceptional divisor
one gets the identity $Z_K=Z_{min}$. In particular, $(X,0)$ is
minimally elliptic (by Laufer's criterion, see \cite{Lauferme}).
Hence $p_g=1$. Moreover, by a calculation $K^2=8-\sum_{i=1}^\nu a_1$,  and thus
\[
8p_g + (K^2+\nu+1)= 17+\nu-\sum_{i=1}^\nu a_i.
\]
Now we will compute the Seiberg--Witten invariant via \ref{ss: sws2}.
The Seifert manifold  $M$ is the unit circle bundle of the
$V$--line bundle ${\bL}_0$ with rational degree $\ell$, and
singularity data $ 1/a_i,\;\;1\leq i\leq \nu$.

The canonical $V$--line bundle $K_\Sigma$   has  singularity data
$(a_i-1)/a_i,\;\;1\leq i\leq \nu,$ \ and rational degree $
\kappa:=-\ell>0$.
 Note that $K_0=-{\bL}_0\in \pic$.
We have $\rho_0=1/2$, $n_0=\lfloor-1/2\rfloor=-1$. The canonical
representative of $\si_{can}$ is the $V$--line bundle
$E_0=-{\bL}_0=K_\Sigma$.

Resolving the inequalities for $S_0^\pm$, one gets
\[
S_0^+ =\{ n{\bL}_0\in {\bZ};\;\; \ell/2 \leq (n+1/2)\ell <0\}
=\{0\cdot \bL_0\},
\]
\[
S_0^-=\{ n{\bL}_0\in {\bZ} ;\;\;0< (n + 1/2)\ell\leq
-\ell/2\}=\{-1\cdot \bL_0\}.
\]
Hence, we  have only two components
$\modu^+_{0},\;\;\modu^-_{K_\Sigma}$, both of dimension $0$.  Thus
$\modu_0^*$ consists of only two  monopoles. Thus $(g_0,0)$ is a
good parameter. The Kreck--Stolz invariant of $\si_{can}$ is
\[
KS_M(\si_{can})=1-2\nu-8\sum_{i=1}^\nu\bms( 1, a_i,
\frac{a_i-1/2}{a_i},-1/2)-4
\sum_{i=1}^\nu\bms(1,a_i)+2\sum_{i=1}^3\frac{1}{a_i}.
\]
The last identity can be further simplified using the identities
from the Appendix B, namely
\[
\bms( 1, a_i, \frac{a_i-1/2}{a_i},-1/2)=
\bms(1,a_i)=\frac{a_i}{12} +\frac{1}{6a_i}-\frac{1}{4}.
\]
Therefore
\[
8\ssw^0_M(\si_{can})=17-2\nu-12\Bigl(\sum_{i=1}^\nu(\frac{a_i}{12}
+\frac{1}{6a_i}
-\frac{1}{4})\Bigr)+2\sum_{i=1}^\nu\frac{1}{a_i}=17+\nu-\sum_{i=1}^\nu
a_i.
\]
We have thus verified the conjecture in this case too.

\gtsubsection{A singularity whose graph is not  star-shaped}
\label{ss: nss}
All the
examples we have analyzed so far had  star shaped resolution
graphs. In this section we consider a different situations which
will indicate  that   the validity of the Main Conjecture extends
beyond singularities whose link is a Seifert manifold. (In this
subsection we will use some standard result about hypersurface
singularities. For these result and the  terminology, the
interested reader can consult \cite{AGV}.)

Consider the isolated plane curve singularity  given by the local
equation \linebreak  $g(x,y):=(x^2+y^3)(x^3+y^2)=0$. We define the surface
singularity $(X,0)$  as the 3--fold cyclic cover of $f$, namely
$(X,0)$  is a hypersurface singularity in $(\C^3,0)$ given by
$f(x,y,z):=g(x,y)+z^3=0$.

The singularity (smoothing) invariants of $f$ can be computed in
many different ways. First notice that it is not difficult to draw
the embedded resolution graph of $g$, which gives all the
numerical smoothing invariants of  $g$.
For example, by A'Campo's formula \cite{A'C} one gets that
the Milnor number of $g$ is 11. Then by Thom--Sebastiani theorem
(see, eg \cite{AGV}, page 60)
$\mu(f)=11\cdot 2=22$. The signature $\si(F)$ of the Milnor fiber
of $F$ can be computed by the method described in \cite{nemsignat}
or \cite{Nem2}; and it is $-18$. Now, by the relations \ref{ss: 4.4}
one gets $p_g(X,0)=1$ and $K^2+\#{\mathcal V}= 10$.

In fact, by the algorithm given in \cite{nemsignat}, one can
compute easily the resolution graph of $(X,0)$ as well (see Figure 3).
\begin{figure}[ht!]
\vspace{3mm}
\centerline{\small
%\ShowGrid 
\SetLabels 
\E(0.01*1.04){$-2$}\\
\E(0.15*1.04){$-2$}\\
\E(0.29*1.04){$-2$}\\
\E(0.43*1.04){$-2$}\\
\E(0.5*1.04){$-3$}\\
\E(0.57*1.04){$-2$}\\
\E(0.71*1.04){$-2$}\\
\E(0.84*1.04){$-2$}\\
\E(0.99*1.04){$-2$}\\
\E(0.11*.5){$-2$}\\
\E(0.8*.5){$-2$}\\
\E(0.11*.03){$-2$}\\
\E(0.8*.03){$-2$}\\
\E(0.01*.85){$\zeta^0$}\\
\E(0.185*.85){$\zeta^0$}\\
\E(0.29*.85){$\zeta^1$}\\
\E(0.43*.85){$\zeta^2$}\\
\E(0.5*.85){$\zeta^0$}\\
\E(0.57*.85){$\zeta^1$}\\
\E(0.71*.85){$\zeta^2$}\\
\E(0.88*.85){$\zeta^0$}\\
\E(0.99*.85){$\zeta^0$}\\
\E(0.185*.5){$\zeta^2$}\\
\E(0.88*.5){$\zeta^1$}\\
\E(0.185*.03){$\zeta^1$}\\
\E(0.88*.03){$\zeta^2$}\\
\endSetLabels 
\AffixLabels{{\includegraphics[width=4in]{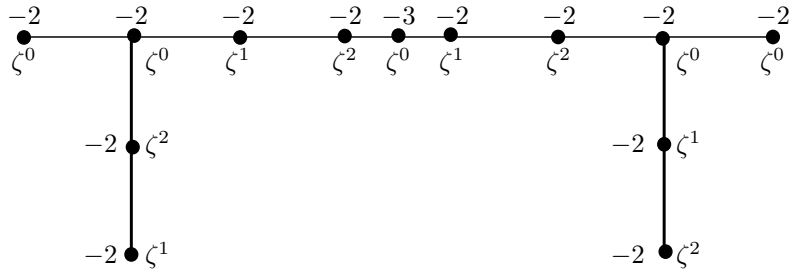}}}}
\vspace{2mm}
\caption{The resolution graph of $(x^2+y^3)(x^3+y^2)=0$}
\end{figure}

Then it is not difficult to verify that the graph satisfies
Laufer's criterion for a minimally elliptic singularity, in
particular this also gives that  $p_g=1$.

Using either way, finally one obtains $p_g+(K^2+\#{\mathcal V})/8=9/4$.
Using the correspondence between the characteristic polynomial $\Delta(t)$ of
the monodromy action  (which can be again easily computed from the
Thom--Sebastiani theorem) and the torsion of $H$ (namely that
$|\Delta(1)|=|H|$), one obtains $H=\Z_3$.

Using the formula for the Casson--Walker invariant from the
plumbing graph one gets $\lambda(M)/|H|=-49/36$.

Finally we have to compute the torsion. There are only two
non-trivial characters. One of them appears one the resolution
graph (ie,  $\chi(g_v)=\zeta^{n_v}$ with $\zeta^3=1$). The other is
its conjugate. Using the general formula for plumbing graphs, one
gets $\tT_{M,\si_{can}}(1)=8/9$.

Since $8/9+49/36=9/4$, the conjecture is true.

\bigskip\bigskip
\appendix\cl{\Large\bf Appendices}\vspace{-6mm}
\small
\section{The Reidemeister--Turaev torsion for plumbings}
\setcounter{equation}{1}
In  this section we prove Theorem \ref{ss: th1}, which  describes the torsion $\tT_M$ of $M$ in terms of plumbing data. The proof has two  parts.

In the first part  we use  Turaev's surgery results \cite{Tu}  formulated  in Fourier theoretic terms which will allow us to replace  formal objects (elements in group algebras) by  analytic ones (meromorphic functions of several variables). We  obtain a first  rough description of $\tT_M$ in terms of surgery data which   has a  $spin^c$ structure ambiguity.

In the second part, we eliminate the ambiguity about the $spin^c$ structure using Turaev's  structure theorem  \cite[Theorem 4.2.1]{Tu5},  and  the identities \ref{ss: swiden} (1) and (2) from our section 3, which completely   determine the $spin^c$ structure  from the Fourier transform of a sign-refined torsion.

\gtsubsection{The surgery data}\label{ss: sur}
We consider an {\em integral} surgery data:
$M$ is a rational homology $3$--sphere   described by the Dehn
surgery on the oriented link   ${\mathcal L}=L_1\cup\cdots \cup L_n
\subset S^3$ with integral surgery coefficients.
We will assume that $n>1$.
 We denote by $E$ the complement of this link.
The manifold $M$ is obtained from $E$ by
attaching $n$ solid tori $Z_1,\cdots, Z_n$.
 We denote by $\mu_i\in H_1(E,\Z)$ the meridian of $L_i$.
Similarly as in the case of plumbing, we can construct lattices
$G:=H_2(M,E,\Z)$ and $G':=\Hom_\Z(G,\Z)$ and a presentation $P\co G\to G'$ for $H:=H_1(M,\Z)$.

 Indeed, the exact sequence $0\to H_2(M,E,\Z)\to H_1(E,\Z)\to H_1(M,\Z)\to 0$ produces a short exact sequence
$0\to G\to G'\to H\to 0$. Moreover, we can identify $G=\Z^n$ via  the canonical basis consisting of classes $[D\times *]$, one
for each solid torus  $Z_i=D\times S^1$;
and $G'=\Z^n$ via the canonical basis determined by the oriented meridians $\{\mu_i\}_i$.  Sometimes we regard  $P$ as a matrix written in these  bases.

Recall that a plumbing graph provides  a canonical surgery presentation
 in such a way that the 3--manifolds obtained by plumbing, respectively by
the surgery, are the same. This presentation is the following:
the components of the link ${\mathcal L}\in S^3$ are in one-to-one
correspondence with the vertices of the graph (in particular, the index set
$I_n=\{1,\ldots,n\}$ is identified with
${\mathcal V}$); all these components are trivial knots in $S^3$;
their framings are the decorations of the
corresponding vertices; two knots corresponding to two vertices connected by
an edge form a Hopf link, otherwise the link is ``the simplest possible''.
In this way we obtain an integral  surgery data  in such a way
that  the matrix $P$ becomes exactly the intersection matrix $I$.

If $\pi\co  G'\ra H$ is the natural projection, then the cores $K_i$ of the attached solid tori $Z_i$ determine the homology classes $\pi (K_i) $ in $H$, and $\pi(K_i)=\pi(-\mu_i)$.

For every  fixed $i\in I_n$  we denote by $E_i$ the manifold
obtained by performing the surgery only along the knots $K_j$,
$j\in I_n\setminus \{i\}$. Equivalently, $E_i$ is the exterior of $K_i$
in $M$.
We set $G_i:=H_2(M,E_i,\Z)$ and $G'_i:=H_1(E_i,\Z)$.
$G_i'$ is generated by the set $\{\mu_j\}_{j=1}^n$ subject to the
relations provided by the $j$-th columns of  $P$ for each $j\not =i$.
There is a natural projection  $G'\ra G'_i$ denoted by $\pi_i$.
Sometimes, for simplicity, we write $K_j$ ($j\not=i$)
 for its projections as well.
We write also  $\tilde{G}:=\Hom(G', {\bC}^*)$ and
 $\tilde{G_i}:=\Hom(G'_i, {\bC}^*)$.
It is natural and convenient to introduce the following definition.

\gtsubsection{Definition}\label{ss: rt2}  A surgery presentation  of a rational homology sphere is called {\em non-degenerate}  \index{Dehn surgery! non-degenerate} if the homology class $\pi_i (K_j)$ has  infinite order in $G'_i$  for any $j\neq i$.

\medskip

\noindent The non-degenerate surgeries can be recognized as follows:
{\em  the surgery is non-degen\-er\-ate
if and only if every off diagonal element of  $P^{-1}$ is  nontrivial.}

Indeed, the fact that for some $j\neq i$ the class
$\pi_{i}(\mu_{j})$ has finite order in $G'_i$
is equivalent with the existence of
 $n\in {\bZ}^*$ and $\vec{v}\in G$ with
$i$-th component $ v_{i}=0$  such that $n\mu_{j}=P\cdot \vec{v}$. But
this says $v_i=nP^{-1}_{ij}$.
Notice that in our case, when  the matrix $P$  is exactly
the negative definite intersection matrix $I$ associated with a connected
(resolution)  graph,
by a well-known result, the surgery presentation is non-degenerate.

We can now begin  the presentation of the surgery formula for the Reidemeister--Turaev torsion.

\gtsubsection{Proposition}\label{ss: rt4} {\sl  Suppose
that  the rational homology $3$--sphere
$M$ is described by a non-degenerate  Dehn surgery.
Fix a   relative $spin^c$ structure
$\tilde{\si}$  on $E$. For any $j$, it induces a relative $spin^c$ structures
$\si_j$ on $E_j$ and a $spin^c$ structure $\si$ on $M$. Let
$\tT_{E_j,\si}$ be  the (sign-refined)  Reidemeister--Turaev
torsion of $E_j$  determined by $\si_j$.  Fix  $\chi\in \hat{H}\setminus \{1\}$
 and $i\in I_n$ such that $\chi(K_i)\neq 1$. Then the following hold.

\smallskip

\noindent {\rm(a)}\qua The  Fourier transform  $\hat{\tT}_{E_i,\si}$  of the torsion
of $E_i$ extends to a holomorphic function on  $\tilde{G}_i\setminus\{1\}$
uniquely determined by the equality
\[
\hat{\tT}_{E_i,\si}(\chi)\cdot \prod_{j\neq i}({\chi}^{-1}(K_j)-1) =
\hat{\tT}_{E,\tilde{\si}}(\chi),\;\;\mbox{for all} \  \chi\in \tilde{G}_i.
\]
Here $\hat{\tT}_{E,\tilde{\si}}$ is the holomorphic extension of Fourier transform of the Alexander--Conway polynomial $\tT_{E,\tilde{\si}}$ of the link
complement $E$, associated with the $spin^c$ structure $\si$ (normalized as in \cite[Section 8]{Tu}).

\smallskip

\noindent {\rm(b)}\qua
$$
\hat{\tT}_{M,\si}(\chi)=\frac{\hat{\tT}_{E_i,\si}(\chi)}
{\bar{\chi}(K_i)-1}.
$$}
\begin{proof}
$\tilde{G}$ is  complex $n$--dimensional torus, and  the  Fourier
transform of the torsion of $E$ extends to a holomorphic function $\chi\mapsto
\hat{\tT}_E(\chi)$ on $\tilde{G}$. The   elements ${K}_j$ also
define  holomorphic functions on $\tilde{G}$ by
$\chi\mapsto \chi({K}_j)^{-1}-1$. Moreover,
$\tilde{G}_i$ is an union of 1--dimensional complex tori
 and the Fourier transform of $\tT_{E_i}$ extends a holomorphic
function $\chi\mapsto \hat{\tT}_{E_i}(\chi)$ on $\tilde{G_i}$.
Since  the elements $K_j$ ($j\not=i$)  have infinite orders in
$G'_i$, we deduce from \cite[Lemma 17.1]{Tu1}, \cite[Section 2.5]{N2}, that   $\hat{\tT}_{E_i}$ is
\emph{the  unique holomorphic extension of the meromorphic function}
\[
\tilde{G}_i\setminus \{1\}\ni\chi\mapsto \frac{\hat{\tT}_{E,\tilde{\si}}(\chi)}{\prod_{j\not=i}({\chi}^{-1}(K_j)-1)}.
\]
Part (b) follows from the surgery formula \cite[Lemma 5.1]{Tu}.
\end{proof}

\gtsubsection{The ``limit'' expression}\label{ss: limit}Let us now explain how we will  use the above  theoretical results.
For each $\chi\in \hat{H}\setminus \{1\}$  pick an arbitrary $i$ with
$\chi(K_i)\neq 1$.    Then $\chi$ belongs to $\tilde{G}_i$ too.
The group $\tilde{G}_i$ is an union of  complex tori,  we denote by
${\bT}_{\chi,i}$ the irreducible  component   containing  $\chi$.
In fact,  there exists $w^i\in G$ such that
${\bT}_{\chi,i}=\{t\ast_{w^i}\chi; t\in {\bC}^*\}$, where
\[
t\ast_{w^i}\chi(v):=t^{v(w^i)}\chi(v)
\;\;\mbox{for all} \  t\in {\bC}^* \ \mbox{and}  \ \;\; v\in G'.
\]
A possible set of ``weights''  $w^i$ can be  determined easily.
$G_i$ is a \emph{free} Abelian group of rank
$1$   which injects into $G$. We can choose $w^i$ to be an arbitrary non-trivial
element of  $G_i$. Obviously, $w^i$ depends on the index $i$. In general, there is no universal choice of the index $i$ which is
suitable for any character $\chi$.

Using the matrix notation, $w^i$ can be regarded as a vector
$\vec{w}^i$ so that  $I\vec{w}^i$ is an (integer) multiple of
$\vec{b}^i$, where $\vec{b}^i$ is the vector whose $i$-th entry is 1,
all the other entries are zero.
Then the above proposition reads as follows:
\[
\hat{\tT}_{M,\si}(\chi)=\frac{1}{\bar{\chi}(K_i)-1}\cdot \lim_{t\ra 1}
\frac{\hat{\tT}_{E,\tilde{\si}}
(t\ast_{w^i}\chi)}{\prod_{j\neq i}\bigl(t\ast_{w^i}{\chi}
(K_j)^{-1} -1\bigr)}.
\]
In particular, by switching the index set $I_n$ to
${\mathcal V}$, if $\chi(\mu_v)\not=1$ for all $v$, one has:
\[
\hat{\tT}_{M,\si}(\chi)=
\frac{\hat{\tT}_{E,\tilde{\si}}(\chi)}
{\prod_v\bigl(\chi (\mu_v) -1\bigr)}.
\]
\gtsubsection{The first part of the proof}\label{ss: fpart}
According to Turaev \cite{Tu},
for any $\tilde{\si}\in$\break $Spin^c(E,\partial E)$
the Alexander--Conway polynomial $\tT_{E,\tilde{\si}}$ has the form
$$\tT_{E,\tilde{\si}}=g\,\prod_v\, (\mu_v-1)^{\delta_v-1},$$
where $g\in G'$ depends on $\tilde{\si}$ by a normalization rule
established by Turaev \cite[Section 8]{Tu} (described in terms of ``charges'').
The generator set  $\{g_v\}_{v\in {\mathcal V}}$ of $H$
(defined via the plumbing)
and  $\{[\mu_v]\}_{v\in {\mathcal V}}$ (defined via the surgery) can be
identified as follows. (Here we will identify their Poincar\'e  duals.)
Consider a resolution $\tilde{X}\to X$ of $(X,0)$ as above.
The lattice inclusion $I\co L\to L'$ (ie, $H^2(\tilde{X},\partial \tilde{X},\Z)\to
H^2(\tilde{X},\Z)$) can be identified with the lattice inclusion
$P=G\to G'$ (ie, $H^2(E,\Z)\to H^2(E,\partial E,\Z)$). Indeed, let $D^4$
be the 4--dimensional ball  with boundary $S^3$ with  ${\mathcal L}\in S^3$.
Then $\tilde{X}$ can be obtained from $D^4$ by attaching $n=\#
{\mathcal V}$ copies of
2--handles $D^2\times D^2$. Let the union of these handles be denoted by
$H$. Clearly, $S^3\setminus \mbox{int}(E)$ is a union $T$ of solid
tori. Then the isomorphism $L'\to G'$ is given by the following sequence of
isomorphism:
$$H^2(E,\partial E)
\stackrel{(1)}{\longleftarrow}
H^2(S^3,T)
\stackrel{(2)}{\longleftarrow}
H^2(D^4,T)
\stackrel{(3)}{\longleftarrow}
H^2(\tilde{X},H)
\stackrel{(4)}{\longrightarrow}
H^2(\tilde{X}).
$$
Above, (1) is an excision, (2) is given by the triple $(D^4,S^3,T)$,
(3) is excision, and (4) is a restriction isomorphism. Moreover,
under this isomorphism, the basis $\{\mu_v\}_v$ correspond exactly to
the basis $\{D_v\}_v$. This also shows that $g_v=[\mu_v]$ for all $v$.
Now, fix a character $\chi\in \hat{H}$, $\chi\not=1$. Set $v_*$ and
$\vec{w}^*\in \Z^n$ as in Theorem \ref{ss: th1}, ie, with
$I\vec{w}^*=-m^*\vec{b}^*$.
Then, using the notations of \ref{ss: limit}, clearly $\vec{w}^*\in \Z^n=G$,
 $\mu_v\in G'$, and $\mu_v(\vec{w}^*)$ is exactly the $v$--component
$w_v^*$ of $\vec{w}^*$. Hence
$t\ast_{\vec{w}^*}\chi(\mu_v)=t^{w_v^*}\chi(g_v)$.
Then by the above notations and \ref{ss: limit}  we conclude that
for any $\si\in Spin^c(M)$
\begin{equation*}
\hat{\tT}_{M,\si}(\chi)=
{\chi}(h)\cdot \lim_{t\to 1}\,
\prod_{v\in {\mathcal V}}\
\Bigl(t^{w_v^*}\, \chi(g_v)-1\Bigl)^{\delta_v-2},
\end{equation*}
for some $h=h(\si)\in H$ which depends (bijectively)  on $\si$.
(Clearly, the limit is not effected by the choice of  $m^*$.)
Now, notice that if we use the identity
$\hat{\tT}_{M,\si}(\bar{\chi})=\hat{\tT}_{M,\bar{\si}}(\chi)$ (cf
\ref{ss: 3.7}(3)), Theorem \ref{ss: th1} is equivalent with the following identity
\begin{equation*}
\hat{\tT}_{M,\bar{\si}}(\chi)=\bar{\chi}(h_\si)\cdot \lim_{t\to 1}\,\prod_{v\in {\mathcal V}}\Bigl(t^{w_v^*}\, \chi(g_v)-1\Bigl)^{\delta_v-2}.
\tag{$t$}
\label{tag: t}
\end{equation*}
The above discussion shows clearly that this is true,  modulo the ambiguity  about $h_\si$. This ambiguity
 (ie, the fact that in the above expression exactly $h_\si$ should be inserted)
is verified via \ref{ss: swiden}(2) (since there is exactly one $h$
which satisfies \ref{ss: swiden}(2) with a fixed $spin^c$ structure $\si$).

\gtsubsection{Additional discussion about the ``weights''}
\label{weights} Before we start the second part,  we clarify an important
fact about the behavior of the weights considered above.
Recall that above, for a fixed $\chi\not=1$, we chose  $v_*$ with
$\chi(g_{v_*})\not=1$. This can be rather unpleasant in any  Fourier
formula, since for different characters we have to take different vertices
$v_*$.
Therefore, we also wish to analyze the case of an arbitrary $v_0$ (disregarding
the fact that $\chi(g_{v_0})$ is 1 or not) instead of $v_*$.

\gtsubsection{Lemma}\label{le} {\sl Fix a character
 $\chi\in \hat{H}\setminus\{1\}$.

{\rm(a)}\qua For an arbitrary vertex $v_0$, consider a vector
$\vec{w}^0$, with components $\{w_v^0\}_v$,
satisfying $I\vec{w}^0=-m^0\vec{b}^0$ for some positive $m^0$.
Then the limit
$$ \lim_{t\to 1}\, \prod_{v\in {\mathcal V}}\
\Bigl(t^{w_v^0}\, \chi(g_v)-1\Bigl)^{\delta_v-2}$$
exists and it is finite.

\smallskip

{\rm(b)}\qua Let $I_{\chi}:=\{ v\, :\, \mbox{either $\chi(g_v)\not=1$ or $v$ has an
adjacent vertex $u$ with $\chi(g_u)\not=1$}\}$.
Then the above limit  is the same for any $v_0\in I_{\chi}$. }

\begin{proof} First we fix some notations.  We say that

$\bullet$\qua  a subgraph
$\Gamma'$ of the plumbing graph $\Gamma$ satisfies the property (P) if
$\sum_{\Gamma'}(\delta_v-2)\geq 0$, where the sum $\sum_{\Gamma'}$ runs over
the vertices of $\Gamma'$ (and $\delta_v$ is the degree of $v$ in $\Gamma$).

$\bullet$\qua $\Gamma'$ is a ``full'' subgraph of $\Gamma$ if any two vertices
of $\Gamma'$ adjacent in $\Gamma$ are adjacent in $\Gamma'$ as well.
For any subgraph $\Gamma'$, we denote by ${\mathcal V}(\Gamma')$ its  set
of vertices.

$\bullet$\qua  a full proper subgraph $\Gamma'$ of $\Gamma$ has property (C)
if it has a unique vertex (say $v_{end}$) which is connected by an edge of
$\Gamma$ with a vertex in ${\mathcal V}(\Gamma)\setminus
{\mathcal V}(\Gamma')$.
For any $\chi\in \hat{H}$,
let $\Gamma^1$ be the full subgraph of $\Gamma$  with set of vertices
$\{v\in {\mathcal V}\,|\,\chi(g_v)=1\}$.
Next, fix a character $\chi\in \hat{H}\setminus \{1\}$ and a vertex $v$
of $\Gamma$. Then
\begin{equation*}
e_vg_v+\sum g_u=0 \ \ \mbox{in $H$, hence }\ \
\chi(g_v)^{e_v}\cdot \prod \chi(g_u)=1,
\tag{1}
\end{equation*}
where the sum (resp. product)
runs over the adjacent vertices $u$ of $v$ in $\Gamma$.
Therefore, if $v$ is in $\Gamma^1$ then
\begin{equation*}
\#\{u\, :\, u\ \mbox{adjacent to $v$ and $u\not\in {\mathcal V}(\Gamma^1)$}\}
\not=1.
\tag{2}
\end{equation*}

{\bf The proof of (a)}\qua
 We have to show that $\Gamma^1$ satisfies (P).
Let $\Gamma^{1,c}$ be
one of its connected components, and denote by $\delta_{v}^{1,c}$
the degree of $v$ in $\Gamma^{1,c}$.
Since $\Gamma^{1,c}$ is a tree,  one has $\sum_{\Gamma^{1,c}}
(\delta_v^{1,c}-2)=-2$.
Since $\Gamma^{1,c}$ is a proper subgraph of the connected
graph $\Gamma$,  there exists at least one edge of $\Gamma$ which
is not an edge of $\Gamma^{1,c}$, but it has one of its end-vertices in
$\Gamma^{1,c}$. In fact, (2) shows that there are at least two such
edges. Therefore, $\Gamma^{1,c}$ satisfies (P).

\smallskip

{\bf The proof of (b)}\qua First we claim the following fact.

\smallskip

(F)\qua {\sl Let  $\Gamma'$ be a full proper subgraph of $\Gamma$  which satisfies (C).
Then for any
\[
v_0\in ({\mathcal V}(\Gamma)\setminus {\mathcal V}(\Gamma')) \cup \{v_{end}\}
\]
the solution $\{w_v^0\}_v$ of $I\vec{w}^0=-m^0\vec{b}^0$ has the following special property: the subset \linebreak $\{v_v^0\}_{v\in {\mathcal V}(\Gamma')}$,
modulo a multiplicative constant, is independent of the choice of $v_0$.}

\smallskip

Indeed, the subset $\{v_v^0\}_{v\in {\mathcal V}(\Gamma')}$,
modulo  a multiplicative constant, is completely determined by the set of
relations  of type (1) considered for vertices $v\in {\mathcal V}(\Gamma')\setminus \{v_{end}\}$. Since the intersection form associated with $\Gamma'$ is
non-degenerate, this system has a maximal rank.
Now, we make a partition of ${\mathcal V}(\Gamma^{1,c})$ (cf part (a)
for the notation).
Each set $S$ of the partition defines a full subgraph $\Gamma^{1,c,j}$
of $\Gamma^{1,c}$ with $S={\mathcal V}(\Gamma^{1,c,j}$).
The partition is defined in such a way that each
$\Gamma^{1,c,j}$ is a maximal  subgraph satisfying both properties (P)
and (C).
One way to construct such a partition is the following.

Let us start with $\Gamma^{1,c}$. By (a), it satisfies (P).
If it does not satisfies (C), then take two of its vertices, both having
adjacent vertices outside of $\Gamma^{1,c}$. Eliminate next all the edges
of $\Gamma^{1,c}$ situated on the path connecting these two vertices, and then, if necessary,  repeat the above procedure for
 the connected components of the remaining
graph. After a finite number of steps all the connected components will
satisfy  both properties (P) and (C).

Now, fact (F) can be applied for all these subgraphs $\Gamma^{1,c,j}$.
In the limit we regroup the product corresponding to the subsets
${\mathcal V}(\Gamma^{1,c,j})$, and the result follows. \end{proof}

\gtsubsection{The second part of the proof: preliminaries}\label{ss: rr}
Our next goal is to show that the right hand side of ($t$)
satisfies the formulae \ref{ss: swiden}(1) and (2) for the $spin^c$
structure $\bar{\si}$.  This clearly ends our proof.

For this, let us fix  a vertex $v_0\in {\mathcal V}$
and we plan to verify \ref{ss: swiden} (1) and (2) for $h=g_{v_0}$.
In the sequel we prefer to fix $m^0$ in the equation of $v_0$, namely we
let $m>0$ be the smallest positive integer  so that
\begin{equation*}
\mbox{$I\vec{w}^0=-m\vec{b}^0$  has a solution $\vec{w}^0=\{w_v^0\}_{v\in
{\mathcal V}}\in \Z^{\#{\mathcal V}}$.}
\tag{3}
\end{equation*}
Clearly $\gcd(\{w_v^0\}_v)=1$ (and each $w_v^0>0$, fact not really important
here).
For a non-trivial character $\chi\in \hat{H}$ with $\chi(g_{v_0})\not=1$,
the vertex $v_0$ is a good candidate for $v_*$
(or, at least, the weights $\vec{w}^*$ in ($t$) can be replaced by the
weights $\vec{w}^0$ since they provide the same limit, cf \ref{le}).
But for characters $\chi$  with $\chi(g_{v_0})=1$ the limit in \ref{le}
(consider for $v_0$) can be different from the limit needed in ($t$)
(where one has  $v_*$). Nevertheless, the products of these
(probably different) limits with $\chi(g_{v_0})-1$  are the same (namely zero)
(and in
 \ref{ss: swiden}(1) and (2) we need only these type of products !). More
precisely, for {\em any} $\chi\in \hat{H}\setminus\{1\}$:
\begin{equation}
(\chi(g_{v_0})-1)\cdot \lim_{t\to 1}\, \prod_{v\in {\mathcal V}}\
\Bigl(t^{w_v^*}\, \chi(g_v)-1\Bigl)^{\delta_v-2}\!\!
=(\chi(g_{v_0})-1)\cdot \lim_{t\to 1}\, \prod_{v\in {\mathcal V}}\
\Bigl(t^{w_v^0}\, \chi(g_v)-1\Bigl)^{\delta_v-2}
\tag*{}
\end{equation}
Therefore, in all our verifications, we can use only one set of weights,
namely $\vec{w}^0=\{w_v^0\}_v$,
 given exactly by the vertex $v_0$, and this is good for all $\chi\in\hat{H}
\setminus \{1\}$.
In the sequel we drop the upper index 0, and we simply write
$w_v$ instead  $w_v^0$.
 Let us introduce the notation
$$\hat{R}_{\chi}(t):= %\lim_{t\to 1}\,
\prod_{v\in {\mathcal V}}\
\Bigl(t^{w_v}\, \chi(g_v)-1\Bigl)^{\delta_v-2}.$$
We have to show that $\lim_{t\to 1} \bar{\chi}(h_{\si})\hat{R}_{\chi}(t)$
satisfies the formulae \ref{ss: swiden}(1) and (2) for the $spin^c$ structure
$\bar{\si}$. Since in these formulae  we need the product of
this limit with $\chi(g_{v_0})-1$, we   set $\bar{\delta}_v:=
\delta_v$ for any $v\not=v_0$, but $\bar{\delta}_{v_0}:=\delta_{v_0}+1$, and
define
$$\hat{P}_{\chi}(t):= % \lim_{t\to 1}\,
\hat{R}_{\chi}(t)\cdot
\Bigl(t^{w_{v_0}}\, \chi(g_{v_0})-1\Bigl)=
%\lim_{t\to 1}\,
\prod_{v\in {\mathcal V}}\
\Bigl(t^{w_v}\, \chi(g_v)-1\Bigl)^{\bar{\delta}_v-2}.$$
In the case of the
trivial character $\chi=1$, we define $\Delta(t)$ via the identity:
\begin{equation*}
\frac{\Delta(t)}{t-1}:=\hat{P}_{1}(t)=
\prod_{v\in {\mathcal V}}\
(t^{w_v}-1)^{\bar{\delta}_v-2}.
\tag{4}
\end{equation*}
Since $\sum_v(\bar{\delta}_v-2)=-1$, one gets that $\Delta(t)$ has no
pole or zero at $t=1$, in fact:
\begin{equation}
\Delta(t)=
\prod_{v\in {\mathcal V}}\
(t^{w_v-1}+\cdots+t +1)^{\bar{\delta}_v-2}.
\tag{5}
\end{equation}
Let $L_0$ be a fixed
generic fiber of the $S^1$--bundle over $E_{v_0}$ used in the plumbing
construction  of $M$ (cf \ref{ss: no}).
Set $G_0:=H_1(M\setminus L_0,\Z)$.

The reader familiar with the theorem of A'Campo about the zeta function
associated with the monodromy action of a Milnor fibration,
certainly realizes that $\hat{P}_1(t)$ is such a zeta function, and
$\Delta(t)$ is a characteristic polynomial of a monodromy operator.
The next proof will not use this possible  interpretation.
Nevertheless, in \ref{ss: fox} we will show that $\Delta(t)\in \Z[t]$, and
$\Delta(1)$ is the order of the torsion subgroup of $G_0$.

 Since $H_2(M,M\setminus L_0,\Z)=\Z$,
one has the exact sequence:
\begin{equation*}
0\to \Z\stackrel{i}{\to}G_0\stackrel{p}{\to}H\to 0,
\end{equation*}
where $i(1_\Z)=\tilde{g}_{\infty}$:=the homology class in $M\setminus L_0$
of the meridian of $L_0$ viewed as a knot in $M$.
Let  $\tilde{g}_v$ be the homology class in $G_0$ of
$\partial D_v$, defined similarly as $g_v\in H$, cf \ref{ss: no}.
Obviously,  $\{\tilde{g}_v\}_{v\in {\mathcal V}}$ is a generator set for
$G_0$. Define $\varphi\co G_0\to \Z$ by $\tilde{g}_v\mapsto w_v$.
The equations (3) guarantee that this is well-defined. Moreover, since
gcd$(\{w_v\}_v)=1$, $\varphi$ is onto. Then clearly, its kernel $T$ is
exactly the subgroup of torsion elements of $G_0$.
Let $j\co T\to G_0$ be the natural inclusion.
Again by (3), $\varphi(\tilde{g}_{\infty})=m$, hence
the composition $\varphi\circ i$
in multiplication by $m$. These facts can be summarized in the
following diagram (where $r$ is induced by $\varphi$):

$$\begin{array}{ccccc}
 & & T & \stackrel{1}{\longleftrightarrow} & T\\
 & & \Big\downarrow\vcenter{\rlap{$j$}}
 & &  \Big\downarrow\vcenter{\rlap{$j'$}} \\
\Z&\stackrel{i}{\hookrightarrow}& G_0& \stackrel{p}{\longrightarrow}&H\\
\Big\updownarrow\vcenter{\rlap{1}} & &
\Big\downarrow\vcenter{\rlap{$\varphi$}} & &
\Big\downarrow\vcenter{\rlap{$r$}} \\
\Z &\stackrel{\cdot m}{ \hookrightarrow}&\Z &
 \stackrel{q}{\longrightarrow}&\Z_m
\end{array}
$$

\noindent It is convenient to identify $\Z_m$ with a subgroup of $\Q/\Z$ via
$\Z_m\ni \hat{a}\mapsto \frac{a}{m}\in \Q/\Z$.

\gtsubsection{Lemma}\label{ss: rbm} {\sl For any $h\in H$, $r(h)=b_M(g_{v_0},h)$
via the above identification.}
\begin{proof} It is enough to verify the identity for each $g_v,\ v\in
{\mathcal V}$. In that case, $r(g_v)=q(\varphi(\tilde{g}_v))=
\hat{w}_v=w_v/m\in \Q/\Z$. But by \ref{ss: 2.2} and \ref{ss: rr}(3),
$b_M(g_{v_0},g_v)=-(\mu_{v_0},\mu_{v})_\Q=-\mu_v(I^{-1}\mu_{v_0})
=w_v/m$ as well.
\end{proof}
\noindent Fix $g\in G_0$ so that $\varphi(g)=1_\Z$. This provides
automatically a splitting of the exact sequence
$0\to T\to G_0\to \Z\to 0$, ie,
a morphism $s\co G_0\to T$ with $s\circ j=id_T$ and
$s(g)=1_T$. In the sequel, we extend any morphism and character
 to the corresponding group-algebras over $\Z$
(and we denote them  by the same symbol).
 For any character  $\xi$ of $G_0$ we define the  representation
$\xi_t\co G_0\to \C[[t,t^{-1}]]^*$ given by $\xi_t(x)=\xi(x)t^{\varphi(x)}$.
(This can be identified with a family of characters.
Indeed, for any fixed $t\in \C^*$ and character $\xi$,  one can define the
character $\xi_t$ given by $x\mapsto
\xi(x)t^{\varphi(x)}$. Eg, for $\xi=p\circ \chi$,
$\xi_t$  is just a more convenient
notation for the action $t\ast_{\vec{w}^0}\chi$,
cf \ref{ss: limit} and \ref{ss: fpart}.)

Now, the point is that the identity (4) has a generalization in the following sense.

\gtsubsection{Theorem}\label{ss: fox} {\sl
For any  character $\xi\in  \hat{G}_0$  define
$$\hat{P}_\xi(t):=
 \prod_{v\in {\mathcal V}}
\Bigl(\xi_t(\tilde{g}_v)-1\Bigl)^{\bar{\delta}_v-2}=
 \prod_{v\in {\mathcal V}}
\Bigl(t^{w_v}\, \xi(\tilde{g}_v)-1\Bigl)^{\bar{\delta}_v-2}.$$
Then there exist an element $\delta\in \Z[G_0]$ such  that the following hold.

\smallskip

{\rm(a)}\qua  For any  $\xi\in  \hat{G}_0$
\[
\hat{P}_{\xi}(t)=\frac{\xi_t(\delta)}{\xi_t(g)-1}.
\]
{\rm(b)}\qua $1_t(\delta)=\Delta(t)$; $1(\delta)=\aug(\delta)=\Delta(1)=|T|$.

\smallskip

{\rm(c)}\qua $s(\delta)=\Sigma_T$, where $\Sigma_T:=\sum_{x\in T}x\in \Z[T]$.}

\begin{proof} From the first part of the proof (cf \ref{ss: rt4}.(b)
and \ref{ss: fpart}(t))
follows that $\lim_{t\to 1}$\break $\hat{P}_{\xi}(t)$,
modulo a multiplicative factor of type
 $\mp \xi(x)$, for some $x\in G_0$, is the Fourier transform
of the Reidemeister--Turaev torsion $\tT$ on $M\setminus L_0$ associated with
some $spin^c$ structure (whose identification is not needed here).
By \cite[4.2.1]{Tu5}, $\tT-\Sigma_T/(1-g)\in \Z[G_0]$, identity valid in
$Q(G_0)$, the ring of quotients of $\Q[G_0]$.
By the first statement
$\hat{P}_\xi(t)=\mp \xi_t(x\tT)$ for any $\xi\in \hat{G}_0$.
Hence, for some $A\in \Z[G_0]$ one has:
\begin{equation*}
\hat{P}_{\xi}(t)=\xi_t(A)\pm \xi_t\Big(\frac{x\cdot\Sigma_T}{g-1}\Big).
\tag{$*$}
\end{equation*}
This identity multiplied by $\xi_t(g-1)$,
for $\xi=1$ and $t\to 1$, and via (4), provides
 $\Delta(1)=\pm |T|$. By (5),  $\Delta(1) $ is positive, hence
in ($*$)  $\pm 1=+1$. Moreover,  $\Delta(1)=|T|$.
Now, if one defines $\delta:=A(g-1)+x\Sigma_T$,
then (a) and (c) follow easily, and $1_t(\delta)=\Delta(t)$
is exactly (4). \end{proof}

\noindent  In order to verify \ref{ss: swiden}(1) and
(2), we will apply the above theorem to special characters of the type
$\xi=\chi\circ p$, where $\chi\in \hat{H}$.
It is clear that for any $y\in \Z[G_0]$ and $h\in H$, the sum
$\sum \chi(h)\cdot \chi\circ p(y)=\sum \chi(h\cdot p(y))$
over $\chi\in \hat{H}$ is an integer multiple of $|H|$. Hence:
\begin{equation*}
\frac{1}{|H|}{\sum_{\chi\in \hat{H}}}'\chi(h)\cdot \chi\circ p(y)=-
\frac{1}{|H|}1(h)\cdot 1(y)=-\frac{1}{|H|}\aug(y)\ (\mmod \ \Z).
\tag{6}
\end{equation*}
Using the splitting of $G_0$ into $T\times \Z$ given by $g$, one can
easily verify that in $Q(G_0)$
\begin{equation*}
\frac{y-s(y)}{g-1}\in \Z[G_0] \ \
\mbox{for any $y\in \Z[G_0]$}.
\tag{7}
\end{equation*}
In the sequel, we write simply $\chi_t$ for $(\chi\circ p)_t$.

\gtsubsection{Verification of \ref{ss: swiden}(1)}\label{(1)} Now
we will verify
that $\bar{\chi}(h)\hat{R}_{\chi}(t)$ ($\chi\in \hat{H}\setminus \{1\}$)
satisfies \ref{ss: swiden}(1) for
$h=h_\si$ (in fact, for any $h\in H$).
For this, fix a vertex $v_0$, and $g\in G_0$ with
$\varphi(g)=1$ as above. Take an arbitrary $x\in G_0$.
 Then we have to show that
\begin{equation*}
\frac{1}{|H|}\lim_{t\to 1} {\sum_{\chi\in \hat{H}}}'
\bar{\chi}(h)\hat{R}_{\chi}(t)(t^{w_{v_0}}\chi(g_{v_0})-1)(\chi_t(x)-1)=
-b_M(g_{v_0},p(x))\ (\mmod \ \Z).
\tag{8}
\end{equation*}
Via \ref{ss: fox}, the left hand side of (8) is
\begin{equation*}
\frac{1}{|H|}\lim_{t\to 1} {\sum_{\chi\in \hat{H}}}'
 \bar{\chi}(h)\cdot
\chi_t(\delta)\cdot \frac{\chi_t(x)-1}{\chi_t(g)-1}.
\tag{9}
\end{equation*}
Set  $a:=\varphi(x)$.   Since
$s(\delta x)=s(\delta)$, (9)  transforms as follows (use (6), (7) and
\ref{ss: fox}):
$$
\frac{1}{|H|}\lim_{t\to 1} {\sum_{\chi\in \hat{H}}}'
 \bar{\chi}(h)\cdot
\chi_t(\frac{\delta x-s(\delta x)}{g-1}-\frac{\delta-s(\delta)}{g-1})
=-\frac{1}{|H|}\lim_{t\to 1}1(h)\cdot
1_t(\frac{\delta x-s(\delta x)}{g-1}-\frac{\delta-s(\delta)}{g-1})$$
$$=-\frac{1}{|H|}\lim_{t\to 1} \frac{\Delta(t) t^a-|T|-\Delta(t)+|T|}{t-1}
=-\frac{1}{|H|}\Delta(1)a=-\frac{a}{m} \ (\mmod \ \Z).
$$
But the right hand side of (8), via \ref{ss: rbm}, is the same $-a/m\in \Q/\Z$.

\gtsubsection{Verification of \ref{ss: swiden}(2)}\label{(2)} Now
we will verify
\begin{equation*}
\frac{1}{|H|}\lim_{t\to 1} {\sum_{\chi\in \hat{H}}}'
\bar{\chi}(h_\si)\hat{R}_{\chi}(t)(t^{w_{v_0}}\chi(g_{v_0})-1)=
-\gq^c(\bar{\si})(g_{v_0})\ (\mmod \ \Z).
\tag{10}
\end{equation*}
The left  hand side is
\begin{equation*}
\frac{1}{|H|}\lim_{t\to 1} {\sum_{\chi\in \hat{H}}}'
 \bar{\chi}(h_\si)\frac{\chi_t(\delta)}{\chi_t(g)-1}.
%\tag{10}
\end{equation*}
The fraction in this expression can be written as (cf (7))
\begin{equation*}
\chi_t\Big(\frac{\delta-s(\delta)}{g-1}\Big)
+\frac{\chi(s(\delta))}{\chi_t(g)-1}.
\tag{11}
\end{equation*}
This sum-decomposition
provides two contributions. The first via (6), (7) and \ref{ss: fox} gives:
$$
\frac{1}{|H|}\lim_{t\to 1} {\sum_{\chi\in \hat{H}}}'
 \bar{\chi}(h_\si)\cdot
\chi_t(\frac{\delta-s(\delta)}{g-1})
%=-\frac{1}{|H|}\lim_{t\to 1}1(h_\si)\cdot 1_t(\frac{\delta-s(\delta)}{g-1})$$
%$$
=-\frac{1}{|H|}\lim_{t\to 1} \frac{\Delta(t)-\Delta(1)}{t-1}
=-\frac{1}{|H|}\Delta'(1) \ (\mmod \ \Z),
$$
where $\Delta'(t) $ denoted the
derivative of $\Delta$ with respect to $t$.
On the other hand, cf (5),
$$\frac{\Delta'(t)}{\Delta(t)}=\sum_v(\bar{\delta}_v-2)\frac{(
t^{w_v-1}+\cdots +t+1)'}{
t^{w_v-1}+\cdots +t+1},$$
hence
$$\frac{\Delta'(1)}{\Delta(1)}=\frac{1}{2}\sum_v(\bar{\delta}_v-2)(w_v-1).$$
Since $\Delta(1)=|T|=|H|/m$, the first contribution is
\begin{equation*}
\frac{1}{|H|}\lim_{t\to 1} {\sum_{\chi\in \hat{H}}}'
 \bar{\chi}(h_\si)\cdot
\chi_t(\frac{\delta-s(\delta)}{g-1})
=-\frac{1}{2m}\sum_v(\bar{\delta}_v-2)(w_v-1).
\end{equation*}
For the second contribution,
notice that $\chi(s(\delta))=\chi(\Sigma_T)$ is zero
unless $\chi$ is in the image of $\hat{r}\co \hat{\Z}_m\to \hat{H}$; if
$\chi$ is in this  image then $\chi(\Sigma_T)=|T|$. For any
$\chi\in \hat{\Z}_m$ we write $\chi(\hat{1})=\zeta$.
Assume that $r(h_{\si})=-\hat{a}$ (or equivalently, $
r(h_{\si})=-\frac{a}{m}\in \Q/\Z$).  Then
$$
\frac{1}{|H|}\lim_{t\to 1} {\sum_{\chi\in \hat{H}}}'
 \bar{\chi}(h_\si)\cdot
\frac{\chi(\Sigma_T)}{\chi_t(g)-1}=
\frac{1}{|H|}\lim_{t\to 1} {\sum_{\chi\in \hat{\Z}_m}}'
 \bar{\chi}(h_\si)\cdot
\frac{|T|}{\chi_t(g)-1}=\frac{1}{m}{\sum_{\zeta\in \Z_m}}'\zeta^a\cdot
\frac{1}{\zeta-1}.$$
Since $$\frac{1}{m}\sum_{\zeta\in \Z_m}
\frac{\zeta^a-1}{\zeta-1}= 0\ (\mmod \ \Z),$$
one gets that the second contribution is (cf \ref{eq: trig1}):
$$\frac{1}{m}{\sum_{\zeta\in \Z_m}}'\zeta^a\cdot
\frac{1}{\zeta-1}=-\frac{a}{m}+\frac{1}{m}{\sum_{\zeta\in \Z_m}}'
\frac{1}{\zeta-1}=-\frac{a}{m}-\frac{1}{2m}(m-1) \ (\mmod \ \Z).$$
Therefore, the left hand side of (10), modulo $\Z$,  is
$$-\frac{1}{2m}\sum_v(\bar{\delta}_v-2)(w_v-1)
-\frac{a}{m}-\frac{1}{2m}(m-1).$$
Notice that  $\sum_v(\bar{\delta}_v-2)=-1$.  Moreover
$w_v=-mI^{-1}_{vv_0}$, and
the coefficient $r_{v_0}$ of $Z_K$ equals
$1-\sum_v(\delta_v-2)I^{-1}_{vv_0}$ (cf \ref{ss: 5.2}),
   hence the above expression can be transformed into
$$-\frac{1}{2}-\frac{1}{2m}\sum_v(\bar{\delta}_v-2)w_v
-\frac{a}{m}=
-\frac{1}{2}r_{v_0}+\frac{1}{2}I^{-1}_{v_0v_0}-\frac{a}{m}.$$
Now, let us compute the right hand side of (10). Since $h_{\si}\cdot \si_{can}
=\si$ one has $2h_{\si}+Z_K=c(\si)$. Then the characteristic element
which provides $\bar{\si}$ is $-c(\si)=-2h_\si-Z_K$. Therefore
$$
-\gq^c(\bar{\si})(g_{v_0})=\frac{1}{2}(D_{v_0}-Z_K-2h_\si,D_{v_0})_\Q=
\frac{1}{2}(\sum_vI^{-1}_{vv_0}E_v-\sum_vr_vE_v,D_{v_0})_\Q-
(h_{\si},D_{v_0})_\Q$$
$$=\frac{1}{2}I^{-1}_{v_0v_0}-\frac{1}{2}r_{v_0}-(h_\si,D_{v_0})_\Q.$$
But, using \ref{ss: 2.2} and \ref{ss: rbm},
 $(h_\si,D_{v_0})_\Q=-b_M(h_\si,g_{v_0})=-r(h)=\frac{a}{m}$. This proves A.12(10). At this point we invoke the  following elementary fact.

\smallskip

\emph{Suppose  $q_1$, $q_2$ are two  quadratic functions  on the finite
abelian group $H$ associated with the bilinear forms $b_1, b_2$;
and $S\subset H$ is a    generating set  such that \,
$q_1(s)=q_2(s)$ \, and \, $b_1(s,h)=b_2(s,h)\;\;\mbox{for all}
\;\;s\in S$ and $ h\in H$. Then  $q_1(h)=q_2(h)$ for all $h\in H$}.

\smallskip

Using A.11(8),  A.12(10)  and the above fact  we obtain \ref{ss: swiden}(2),
for any $h$. The
identity \ref{ss: swiden}(2)  implies that $\bar{\chi}(h_\si)\hat{R}_\chi(t)
=\hat{\tT}_{M,\bar{\si}}(\chi)$. This concludes the proof of Theorem
\ref{ss: th1}.

(Notice that, in fact, we verified even more. First recall, cf
 \cite{Tu}, that  the {\em sign}   of the (sign-refined) torsion is
decided by universal rules.
In some cases its identification is rather involved.
The point is that the above verification also reassures us that in $(t)$
we have the right sign.)\endproof

\section{Basic facts concerning the Dedekind--Rademacher sums}
\setcounter{equation}{0}
In this Appendix we collected some facts about
(generalized)  Dedekind  sums  which  constitute a necessary minimum
in the concrete computation of the Seiberg--Witten invariants (and
in the understanding of the relationship between
Dedekind sums and Fourier analysis).
Let $\lfloor x\rfloor$ be  the integer part of $x$, and
$\{x\}:= x-\lfloor x\rfloor$ its fractional part.
 In the paper \cite{Ra},  Rademacher  introduces for every
pair of coprime  integers $h,k$ and any real numbers  $x,y$  the
following generalization of the classical Dedekind sum
\[
\bms(h,k;x,y)=\sum_{\mu=0}^{k-1}\LL \frac{\mu+y}{k}\RR \LL \frac{
h(\mu +y) }{k}+x \RR=-\bms(-h,k;-x,y),
\]
where  $((x))$ denotes the Dedekind symbol
\[
((x))=\left\{
\begin{array}{ccc}
\{x\} -1/2 & {\rm if} & x\in {\bR}\setminus {\bZ}\\
0 & {\rm if} & x\in {\bZ}.
\end{array}
\right.
\]
A simple computations shows that $\bms(h,k; x,y)$ depends only on
$x,y$ (mod 1).  Additionally
\[
\bms(h,k;x,y)=\bms(h-mk, k; x+my,y).
\]
Moreover, have the following result
\begin{equation}
\bms(1,k;0,y)=\left\{
\begin{array}{lr}
\frac{k}{12}+\frac{1}{6k}-\frac{1}{4}=\frac{(k-1)(k-2)}{12k} & y \in {\bZ} \\
                                     &   \\
\frac{k}{12} +\frac{1}{k}B_2(\{y\})  & y \in {\bR}\setminus {\bZ}\, ,
\end{array}
\right. \label{eq: quad}
\end{equation}
where $B_2(t)=t^2-t+1/6$ is the second Bernoulli polynomial.
If $x=y=0$ then we simply write $\bms(h,k)$.
Perhaps the  most important property of these Dedekind--Rademacher
sums is their reciprocity law which makes them    computationally
very friendly: their computational complexity is comparable with
the complexity of the classical Euclid's algorithm. To formulate
it we must distinguish two cases.

\noindent  $\bullet$\qua  Both $x$ and $y$ are integers.  Then
\begin{equation}
 \bms(h,k) +\bms(k,h)=-\frac{1}{4}+\frac{h^2+k^2+1}{12hk}.
\label{eq: rec1}
\end{equation}
$\bullet$\qua  $x$ and/or $y$ is not an integer. Then
\[
\bms(h,k;x,y)+\bms(k,h;y,x)
\]
\begin{equation}
=((x))\cdot ((y)) + \frac{h^2\psi_2(y) +\psi_2(h y+kx) +
k^2\psi_2(x)}{2hk} \label{eq: rec2}
\end{equation}
where $\psi_2(x):= B_2(\{x\})$. In particular, if $x,y\in{\bR}$
are not both  integers we deduce
\begin{equation}
\bms(1,m;x,y)=-((x))\cdot ((mx+y))+((x))((y))+\frac{\psi_2(y)
+\psi_2(y+mx) +m^2\psi_2(x)}{2m}. \label{eq: quad1}
\end{equation}
An important ingredient behind  the reciprocity law is the
following identity (\cite[Lemma 1]{Ra})
\begin{equation}
\sum_{\mu=0}^{k-1}\LL\frac{\mu+w}{k}\RR =((w))\;\;\mbox{for any } \ w\in
{\bR}. \label{eq: kubert}
\end{equation}
The various Fourier--Dedekind sums  we use in this paper  can be
expressed in terms of Dedekind--Rademacher sums.  This  follows
from the identity (\cite[page 170]{HZ})
\begin{equation}
\frac{1}{p}{\sum_{\zeta^p=1}}'\frac{\zeta^t}{1-\zeta}=
\LL\frac{2t-1}{2p}\RR,\;\;\mbox{for all} \
p,q\in {\bZ}, \; p>1. \label{eq: trig1}
\end{equation}
In other words, the function
\[
\{\zeta^p=1\}\subset {\bC}^*\ra
{\bC},\;\;\zeta\mapsto \left\{
\begin{array}{ccc}
0 & {\rm if} & \zeta=1\\
\frac{1}{1-\zeta} & {\rm if} & \zeta \neq 1
\end{array}
\right.
\]
is the Fourier transform of the function
\[
\Z_p\ra {\bC},\;\;  \hat{t}\mapsto
\LL\frac{2t-1}{2p}\RR.
\]
The identity (\ref{eq: trig1}) implies that
\[
\frac{1}{p}{\sum_{\zeta^p=1}}'
\frac{\zeta^{tq'}}{1-\zeta^q}=\LL\frac{2tq'-1}{2p}\RR,\ \mbox{for all}\
p,q\in {\bZ}, \; p>1,\; (p,q)=1,
\]
where $q'=q^{-1}\ (\mmod \ p)$.  Using the fact that   Fourier
transform  of the convolution product of two  functions
$\Z_p\ra {\bC}$ is the pointwise product of the Fourier
transforms of these functions we deduce after some simple
manipulations the following identity.
\begin{equation}
\frac{1}{p}{\sum_{\zeta^p=1}}'\frac{\zeta^t}{(\zeta-1)(\zeta^q-1)}=
-{\bms}\bigl(q,p;\frac{q+1-2t}{2p},-\frac{1}{2}\bigr).
\label{eq: trig2}
\end{equation}
If $t=0$ then by (\ref{eq: kubert}) (and a computation), or by
\cite[18a]{RG}, one has
\begin{equation}
\frac{1}{p}{\sum_{\zeta^p=1}}'\frac{1}{(\zeta-1)(\zeta^q-1)}=
-{\bms}(q,p)+\frac{p-1}{4p}.
\label{eq: trig2bis}
\end{equation}
By setting $q=-1$ and $t=0$ in the above equality we deduce
\begin{equation}
\frac{1}{p}{\sum_{\zeta^p=1}}'\frac{1}{|\zeta-1|^2}=-{\bms}(-1,p,0,-1/2)=
{\bms}(1,p,0,1/2)\stackrel{(\ref{eq:
quad})}{=}\frac{p}{12}-\frac{1}{12p}. \label{eq: trig3}
\end{equation}
The Fourier transform of the function
$ \mathfrak{d}_{p,q}:\Z_p\ra {\bC},\ \ \hat{t}\mapsto ((qt/p))
$
is the function (see \cite[Chapter 2, Section C]{RG})
\[
\{\zeta^p=1\}\ra {\bC},\;\;\zeta \mapsto \left\{
\begin{array}{ccc}
\frac{1}{2}-\frac{\zeta^q}{\zeta^q-1} & {\rm if} & \zeta\neq 1\\
0 & {\rm if} &\zeta=1.
\end{array}
\right.
\]
Then
\[
\bms(-q,p)=\sum_{t+s=0\, (\mbox{\tiny{mod}}\, p)}
\mathfrak{d}_{p,1}(t)\mathfrak{d}_{p,q}(s)=
(\mathfrak{d}_{p,1}\ast\mathfrak{d}_{p,q})(0)
\]
\[
=\frac{1}{p}{\sum_{\zeta^p=1}}'\Bigl(\frac{1}{2}-\frac{\zeta}{\zeta-1}\Bigr)\Bigl(\frac{1}{2}-\frac{\zeta^{q}}{\zeta^q-1}\Bigr).
\]
This implies (cf also with \cite{RG})
\begin{equation}
\frac{1}{p}{\sum_{\zeta^p=1}}'\frac{\zeta+1}{\zeta-1}\cdot\frac{\zeta^q+1}
{\zeta^q-1}=-4\bms(q,p).
\label{eq: trig4}
\end{equation}


\begin{thebibliography}

\bibitem{A'C} {\bf N A'Campo},  {\it La fonction zeta d'une monodromy}, Com. Math.
Helvetici, {50} (1975) 233--248

\bibitem{AGV} {\bf V\,I Arnold}, {\bf S\,M Gusein-Zade}, {\bf A\,N Varchenko},  {\it
Singularities of Differentiable Mappings}, Vol. 2 , Birkhauser, Boston (1988)

\bibitem{Artin} {\bf M Artin}, {\it Some numerical criteria for contractibility of
curves on algebraic surfaces},  Amer. J. of Math. {84} (1962) 485--496

\bibitem{Artin2} {\bf M Artin}, {\it On isolated rational singularities of
surfaces}, Amer. J. of Math. {88} (1966) 129--136

\bibitem{Bries} {\bf E Brieskorn}, {\it Beispiele zur Differentialtopologie von
Singularit\"aten}, Inventiones math. {2} (1969) 1--14

\bibitem{Chen1} {\bf W Chen}, {\it Casson invariant and Seiberg--Witten gauge
theory}, Turkish J. Math. {21} (1997) 61--81

\bibitem{CoS} {\bf O Collin}, {\bf N Saveliev}, {\it A geometric proof of the
Fintushel--Stern formula}, Adv. in Math. {147} (1999) 304--314

\bibitem{Co} {\bf O Collin}, {\it Equivariant Casson invariant for knots and the
Neumann--Wahl formula}, Osaka J. Math. {37} (2000) 57--71

\bibitem{Durfee} {\bf A Durfee}, {\it The Signature of Smoothings of Complex
Surface Singularities}, Math. Ann. {232} (1978) 85--98

\bibitem{FS} {\bf R Fintushel},  {\bf R Stern} {\it Instanton homology of Seifert
fibered homology three spheres}, Proc. London Math. Soc. {61} (1990)
109--137

\bibitem{FuS}  {\bf M Furuta}, {\bf B Steer}, {\it  Seifert fibered homology 3--spheres
and the Yang--Mills equations on Riemann surfaces with marked points}, Adv.
in Math. {96} (1992) 38--102

\bibitem{Grauert} {\bf H Grauert}, {\it \"Uber Modifikationen und
exceptionelle analytische Mengen}, Math. Ann. {146} (1962) 331--368

\bibitem{GSt} {\bf G-M  Greuel}, {\bf J\,H\,M Steenbrink}, {\it On the topology of
smoothable singularities}, Proc. Symp. of Pure Math. {\bf 40},
Part 1 (1983) 535--545

\bibitem{Gompf}{\bf R\,E Gompf}, {\it Handlebody construction of Stein surfaces},
Ann. of Math. 148 (1998)  619--693

\bibitem{GS} {\bf R\,E  Gompf}, {\bf A\,I Stipsicz}, {\it An Introduction to
$4$--Manifolds and Kirby Calculus}, Graduate Studies in Mathematics, vol. 20,
Amer. Math. Soc. (1999)

\bibitem{Hamm} {\bf H A Hamm}, {\it Exotische Sph\"aren als Umgebungsr\"nder in
speziellen komplexen R\"aumen}, Math. Ann, {197} (1972) 44--56

\bibitem{Hi} {\bf F Hirzebruch},  {\it Pontryagin classes of rational homology
manifolds and the signature of some affine hypersurfaces},
Proceedings of the Liverpool Singularities Symposium II (ed.
C.T.C. Wall) Lecture Notes in Math. {209}, Springer Verlag (971) 207--212

\bibitem{Ishii} {\bf S Ishii}, {\it The invariant $-K^2$ and continued fractions
for 2--dimensional cyclic quotient singularities}, preprint

\bibitem{JN} {\bf M Jankins},  {\bf W D Neumann}, {\it Lectures on Seifert
Manifolds}, Brandeis Lecture Notes (1983)

\bibitem{HZ} {\bf F Hirzebruch},  {D Zagier}, {\it The Atiyah--Singer Index
Theorem and Elementary Number Theory},  Math. Lect. Series {3}, Publish
or Perish Inc. Boston (1974)


\bibitem{Lauferme} {\bf H Laufer}, {\it On minimally elliptic singularities},
Amer. J. of Math. {99} (1977) 1257--1295

\bibitem{Laufermu} {\bf H Laufer}, {\it On $\mu$ for surface singularities},
Proc. of  Symp. in Pure Math. {30} (1977) 45--49

\bibitem{Lescop} {\bf C Lescop}, {\it Global Surgery Formula for the Casson--Walker
Invariant}, Annals of Math. Studies, vol. {140}, Princeton University
Press (1996)

\bibitem{Lim} {\bf Y Lim} {\it  Seiberg--Witten invariants for 3--manifolds in the
case $b_1=0$ or $1$}, Pacific J. of Math. {195} (2000) 179--204

\bibitem{LW}  {\bf E Looijenga},  {\bf J Wahl}, {\it  Quadratic functions and smoothing surface singularities}, Topology, {25} (1986) 261--291

\bibitem{MW}  {\bf M Marcolli}, {\bf B\,L Wang}, {\it  Seiberg--Witten invariant and the Casson--Walker invariant for rational homology 3--spheres},
{math.DG/0101127}, Geometri{\ae} Dedicata, to appear

\bibitem{MOY} {\bf T Mrowka}, {\bf  P Ozsvath}, {\bf B Yu}, {\it Seiberg--Witten monopoles on
Seifert fibered spaces}, Comm. Anal. and Geom. {5} (1997) 685--791

\bibitem{Neminv}{\bf  A N\'emethi}, {\it``Weakly'' elliptic Gorenstein
singularities of surfaces} ,  Invent. math. {137} (1999) 145--167

\bibitem{Nem1} {\bf  A N\'emethi}, {\it Five lectures on normal surface singularities},
{\it Proceedings of the  summer school}, Bolyai Society Mathematical Studies
 8, Low Dimensional Topology (1999)

\bibitem{nemsignat}  {\bf  A N\'emethi}, {\it The signature of $f(x,y)+z^n$},
Proceedings of Real and Complex Singularities , (C.T.C Wall's 60th
birthday meeting), Liverpool (England), August 1996; London Math.
Soc. Lecture Note Series, {263} (1999) 131--149

\bibitem{Nem2} {\bf  A N\'emethi}, {\it Dedekind sums and the
signature of $z^N+f(x,y)$}, Selecta Math. {4} (1998) 361--376

\bibitem{Neu} {\bf W Neumann}, {\it Abelian covers of quasihomogeneous surface
singularities}, {\it Proc. of Symposia in
Pure Mathematics}, vol. 40, Part 2, 233--244

\bibitem{NP} {\bf W Neumann},  {\it  A calculus for plumbing applied to the
topology of complex surface singularities and degenerating complex
curves}, Transactions of the AMS, {268} (1981) 299--344


\bibitem{NeR} {\bf W Neumann}, {\bf F Raymond},  {\it Seifert manifolds, plumbing,
 $\mu$--invariant and orientation reserving maps}, Algebraic and Geometric
Topology (Proceedings, Santa Barbara 1977), Lecture Notes in Math.
{664}, 161--196

\bibitem{NW} {\bf W Neumann}, {\bf  J Wahl}, {\it  Casson invariant of links
of singularities}, Comment. Math. Helvetici, {65} (1990) 58--78

\bibitem{NW2} {\bf W Neumann}, {\bf  J Wahl}, {\it Universal abelian covers of surface singularities},
 {\tt arXiv:math.AG/0110167}

\bibitem{Nico0} {\bf L\,I Nicolaescu}, {\it Adiabatic limits of the Seiberg--Witten
equations on Seifert manifolds}, Communication in Analysis and Geometry, {6} (1998) 331--392

\bibitem{Nico2} {\bf L\,I Nicolaescu}, {\it Finite energy Seiberg--Witten moduli spaces on
$4$--manifolds bounding Seifert fibrations},  Comm. Anal. Geom. {8} (2000)
 1027--1096

\bibitem{Nico3}{\bf L\,I Nicolaescu},  {\it Seiberg--Witten  invariants of lens spaces},
Canad J. of Math. {53} (2001) 780--808

\bibitem{N3} {\bf L\,I Nicolaescu}, {\it On the Reidemeister torsion of rational homology
spheres},  Int. J. of Math. and Math. Sci. {25} (2001) 11--17

\bibitem{Nico5}{\bf L\,I Nicolaescu}, \,{\it Seiberg--Witten invariants of rational homology
spheres}, {\tt arXiv:math.DG/0103020}

\bibitem{N2} {\bf L\,I Nicolaescu}, {\it Notes on  the Reidemeister Torsion},  Walter de
Gruyter, to appear

%\bibitem{Pinkh}  {\bf H Pinkham} {\it Normal surface singularities with ${\bC}^*$ action}, Math. Ann. {117} (1977) 183--193

\bibitem{Ra} {\bf H Rademacher}, {\it Some remarks on certain generalized
Dedekind sums}, Acta Arithmetica, {9} (1964) 97--105

\bibitem{RG} {\bf H Rademacher}, {\bf E Grosswald}, {\it Dedekind Sums},
The Carus Math. Monographs, MAA (1972)

\bibitem{Ratiu} {\bf A Ratiu}, {\it PhD Thesis},  Paris VII

\bibitem{Seade} {\bf J\,A Seade}, {\it A cobordism invariant for surface
singularities}, Proc. of Symp. in Pre Math. Vol. {40}, Part 2, 479--484
(1983)

\bibitem{Steenbrink}{\bf  J\,H\,M Steenbrink}, {\it Mixed Hodge structures
associated with isolated singularities}, Proc. Sumpos. Pure Math. {40},
Part 2, 513--536 (1983)

\bibitem{Tu5}{\bf V\,G  Turaev},  {\it  Torsion invariants of $Spin^c$--structures
on $3$--manifolds}, Math. Res. Letters, {4} (1997) 679--695

\bibitem{Tu} {\bf V\,G  Turaev}, {\it  Surgery formula for torsions and
Seiberg--Witten invariants of $3$--manifolds}, {\tt arXiv:math.GT/0101108}

\bibitem{Tu1} {\bf V\,G  Turaev}, {\it  Introduction to Combinatorial Torsions},
Lectures in Mathematics, ETH Zurich, Birkh\"{a}user (2001)

\bibitem{Blij} {\bf F Van der Blij}, {\it An invariant of quadratic
forms mod $8$}, Indag. Math. {21} (1959) 291--293

\bibitem{W} {\bf K Walker}, {\it An extension of Casson's invariant},
Annals of Math. Studies, vol. {126}, Princeton University
Press (1996)

\bibitem{Yau} {\bf S S-T Yau}, {\it On maximally elliptic singularities},
Transact. AMS, {257} (1980) 269--329
\end{thebibliography}
\end{document}